\documentstyle[amssymb,amsfonts]{amsart}

\def\R{{\mathbb{R}}}
\def\C{{\mathbb{C}}}
\def\Q{{\mathbb{Q}}}

\def\P{{\mathbb{P}}}
\def\E{{\mathbb{E}}}

\def\N{{\mathbb{N}}}
\def\L{{\mathbf{L}}}
\def\x{{\mathbf{x}}}
\def\B{{\mathcal{B}}}
\def\A{{\mathcal{A}}}
\def\D{{\mathcal{D}}}
 at 10 true pt

\def\tr{\operatorname{tr}}
\def\Re{\operatorname{Re}}
\def\Im{\operatorname{Im}}
\def\Z{{\mathbb{Z}}}

\def\eps{\varepsilon}
\newenvironment{proof}{\noindent {\bf Proof} }{\endprf\par}
\def \endprf{\hfill  {\vrule height6pt width6pt depth0pt}\medskip}
\def\emph#1{{\it #1}}
\def\textbf#1{{\bf #1}}

\theoremstyle{plain}
  \newtheorem{theorem}[subsection]{Theorem}

  \newtheorem{proposition}[subsection]{Proposition}
  
  \newtheorem{lemma}[subsection]{Lemma}
  \newtheorem{hypothesis}[subsection]{Hypothesis}
  \newtheorem{example}[subsection]{Example}
  \newtheorem{examples}[subsection]{Examples}
  \newtheorem{corollary}[subsection]{Corollary}

\theoremstyle{remark}
  \newtheorem{remark}[subsection]{Remark}
  \newtheorem{remarks}[subsection]{Remarks}

\theoremstyle{definition}
  \newtheorem{definition}[subsection]{Definition}

\include{psfig}

\begin{document}

\title[Constellations in the Gaussian primes]{The Gaussian primes contain arbitrarily shaped constellations}

\author{Terence Tao}
\address{Department of Mathematics, UCLA, Los Angeles CA 90095-1555}
\email{tao@@math.ucla.edu}

\begin{abstract}  We show that the Gaussian primes $P[i] \subseteq \Z[i]$ contain infinitely constellations of any prescribed shape and orientation.  More precisely, given any distinct Gaussian integers $v_0,\ldots,v_{k-1}$, we show that there are infinitely many sets $\{a+rv_0,\ldots,a+rv_{k-1}\}$, with $a \in \Z[i]$ and $r \in \Z \backslash \{0\}$, all of whose elements are Gaussian primes.  

The proof is modeled on that in \cite{gt-primes} and requires three ingredients.  The first is a hypergraph removal lemma of Gowers and R\"odl-Skokan,
or more precisely a slight strengthening of this lemma which can be found in \cite{tao:hyper}; this hypergraph removal lemma 
can be thought of as a generalization of the Szemer\'edi-Furstenberg-Katznelson theorem concerning multidimensional arithmetic
progressions.  The second ingredient is the transference argument from \cite{gt-primes}, which allows one to extend this hypergraph removal lemma to a relative version, weighted by a pseudorandom measure.  The third ingredient is a Goldston-Y{\i}ld{\i}r{\i}m type analysis for the Gaussian integers, similar to that in \cite{gt-primes}, which yields a pseudorandom measure which is concentrated on Gaussian ``almost primes''.
\end{abstract}

\maketitle

\section{Introduction}

A famous and deep theorem of Szemer\'edi \cite{szemeredi}
asserts that any set of integers of positive upper density contains arbitrarily long arithmetic progressions.  This theorem was extended by
Furstenberg and Katznelson \cite{fk} to higher dimensions, as follows.  If $Z$ is an additive group, we define a \emph{shape} in $Z$
to be a finite collection $(v_j)_{j \in J} \in Z^J$ of distinct elements in $Z$.  A \emph{constellation} in $Z$ with this shape is defined to be 
any $J$-tuple of the form $(a + rv_j)_{j \in J} \in Z^J$, where $a \in Z$ and $r \in \Z$, with all of the $a+rv_j$ being distinct.  
Note that we can define the product of an integer $r \in \Z$ with an additive group element $v \in Z$ in the usual manner.  Thus a constellation is nothing more than a homothetic copy of a given shape.

\begin{theorem}[Multidimensional Szemer\'edi's theorem, combinatorial version]\label{multiszemeredi}\cite{fk}  Let $d \geq 1$, and let $A$ be a subset of the lattice $\Z^d$ whose upper Banach density is strictly positive, thus
$$\limsup_{N \to \infty} \frac{|A \cap [-N,N]^d|}{|[-N,N]|^d} > 0,$$
where $[-N,N] := \{n \in\Z: -N \leq n \leq N\}$ and $|A|$ denotes the cardinality of $A$. Then for
any given shape $(v_j)_{j \in J}$ in $\Z^d$, the set $A$ contains infinitely many constellations $(a+rv_j)_{j \in J}$ with that shape.
\end{theorem}

Now consider the \emph{Gaussian primes} $P[i]$ in the \emph{Gaussian integers} $\Z[i] := \{ a+bi: a,b \in \Z\}$, defined
as those Gaussian integers $p \in \Z[i]$ which have no proper factors (other than \emph{units} $1, i, -1, -i$ and \emph{associates} $p, ip, -p, -ip$).  One can identify
$\Z[i]$ with $\Z^2$ in the obvious manner, however when one does so, the upper Banach density of $P[i]$ is zero and so Theorem \ref{multiszemeredi} does not directly apply.  Nevertheless, we are able to establish the following result, which is the main result of this paper.

\begin{theorem}[Constellations in the Gaussian primes]\label{main}  Let $(v_j)_{j \in J}$ be any shape in the Gaussian integers $\Z[i]$.
Then the Gaussian primes $P[i]$ contains infinitely many constellations with this shape.
\end{theorem}

Theorem \ref{main} can be thought of as the Gaussian counterpart of the recent result in \cite{gt-primes}
that the rational primes $P = \{2,3,5,\ldots\}$ contain arbitrarily long arithmetic progressions.  The latter result is connected to the
$d=1$ case of Theorem \ref{multiszemeredi}, whereas the results here are connected to the $d=2$ case.  It is likely that the method also extends
to cover some further results of this type, see Section \ref{discussion-sec}.  For instance, one can replace $P[i]$ in the above theorem by any subset of $P[i]$ of positive upper \emph{relative} Banach density, as in \cite{gt-primes}.  We remark that the scaling parameter $r$ can be chosen to be positive,
by the rather crude expedient of replacing the constellation $(v_j)_{j \in J}$ with the symmetrized constellation 
$(v_j)_{j \in J} \uplus (-v_j)_{j \in J}$.

Our approach to proving Theorem \ref{main} basically follows the strategy of \cite{gt-primes}.  A direct execution of that strategy would
proceed by somehow transferring Theorem \ref{multiszemeredi} to a relative version, weighted by a pseudorandom measure.  One would then construct
a pseudorandom measure concentrated on the Gaussian ``almost primes'' to conclude the argument.  It may well be possible to carry out this approach;
however we have proceeded by a slightly different route, not working with Theorem \ref{multiszemeredi} but a stronger result, which we call a
``strong hypergraph removal lemma'', which we shall discuss shortly.  (We will, however, still need to construct a pseudorandom measure concentrated in Gaussian almost primes.)

Theorem \ref{multiszemeredi} in the contrapositive, implies in particular that any subset of $\Z^d$ which contains only finitely many constellations
of a prescribed shape, must have density zero.  A more quantitative version of this assertion is as follows.  Given any finite non-empty set $Z$ and any function $f: Z \to \R$,
we use $\E(f) = \E(f|Z) = \E(f(x)|x \in Z) := \frac{1}{|Z|} \sum_{x \in Z} f(x)$ to denote the average value of $f$.  If $x,y_1,\ldots,y_n$ are parameters and $X > 0$ is a positive quantity,
we use $o_{x \to 0; y_1,\ldots,y_n}(X)$ to denote any quantity bounded in magnitude by $c(x,y_1,\ldots,y_n)X$, where $c$ is a function which goes to zero as $x \to 0$ for each fixed choice of $y_1,\ldots,y_n$.  Similarly we use $O_{y_1,\ldots,y_n}(X)$ to denote any quantity bounded in magnitude
by $C(y_1,\ldots,y_n)X$ for some quantity $C(y_1,\ldots,y_n) > 0$.

\begin{theorem}[Multidimensional Szemer\'edi's theorem, expectation version]\label{sz-multi-alt} Let $Z, Z'$ be two finite additive groups,  and let $(\phi_j)_{j \in J}$ be a finite collection of group homomorphisms $\phi_j: Z \to Z'$.  Let $A$ be a subset of $Z'$.  If we have
$$ \E\left( \prod_{j \in J} 1_A(x + \phi_j(r)) \big| x \in Z'; r \in Z\right) \leq \delta$$
for some $0 < \delta \leq 1$, then we have 
$$ \E( 1_A(x) | x \in Z') = o_{\delta \to 0;|J|}(1).$$
\end{theorem}

This particular result does not appear explicitly in the literature, but it follows from the work of Furstenberg and Katznelson
\cite{fk} in the cyclic case $Z = \Z/N\Z$, and from their later work \cite{fk2} on a density version of the Hales-Jewett theorem for the
general case.  It also follows from the hypergraph analysis of Gowers \cite{gowers-hyper}
and R\"odl-Skokan \cite{rodl}, \cite{rodl2}, or more precisely from Theorem \ref{main-2} below.
It is easy to see that Theorem \ref{sz-multi-alt} implies Theorem \ref{multiszemeredi}, by localizing the situation in Theorem \ref{multiszemeredi}
to a cyclic group such as $\Z_N^d$ for a large prime $N$, and then letting $N \to \infty$; we omit the standard details.  

The proof of Theorem \ref{sz-multi-alt} sketched above used methods from ergodic theory.  At first glance, it seems that the additive structures of the groups $Z$ and $Z'$ must play a key role; for instance, in the ergodic arguments of \cite{fk}, this structure is captured in the algebra of multiple commuting
shifts on a probability space.  However, it is a remarkable fact, observed by multiple authors, that Theorem \ref{sz-multi-alt} (and hence Theorem \ref{multiszemeredi}) can in fact be deduced from a stronger result - namely a ``hypergraph removal lemma'' - in which no additive structure is present.  We shall state this stronger result (or more precisely, a refinement of this result
in \cite{tao:hyper}) shortly, but first we need some notation.

\begin{definition}[Hypergraphs] If $J$ is a finite set and $d \geq 0$, we define ${J \choose d} := \{ e \subseteq J: |e| = d \}$ to be the set of all
subsets of $J$ of cardinality $d$.  A \emph{$d$-uniform hypergraph} on $J$ is then defined to be any subset $H \subseteq {J \choose d}$ of
${J \choose d}$.
\end{definition}

\begin{definition}[Hypergraph systems]  A \emph{hypergraph system} is a quadruplet $V = (J, (V_j)_{j \in J}, d, H)$, where $J$ is a finite set,
$(V_j)_{j \in J}$ is a collection of finite non-empty sets indexed by $J$, $d \geq 1$ is positive integer, and $H \subseteq {J \choose d}$ is a
$d$-uniform hypergraph.  For any $e \subseteq J$, we set $V_e := \prod_{j \in e} V_j$, and let $\pi_e: V_J \to V_e$ be the canonical projection map.  For each $e \in J$, let $\A_e$ be the $\sigma$-algebra on $V_J$ defined by $\A_e := \{ \pi_e^{-1}(E): E \subseteq V_e \}$.  
\end{definition}

\begin{remark}
Very roughly speaking, a hypergraph system corresponds to the notion of a \emph{measure-preserving system} in ergodic theory, though with the
notable difference that no analogue of the shift operator exists in a hypergraph system.  Indeed the $V_j$ are simply finite sets, and need not have any additive structure whatsoever.
\end{remark}

\begin{theorem}[Hypergraph removal lemma]\label{main-2}\cite{gowers-hyper}, \cite{nrs}, \cite{rodl}, \cite{rodl2}, \cite{tao:hyper}  
Let $V = (J, (V_j)_{j \in J}, d, H)$ be a hypergraph system.  For each $e \in H$, let $E_e$ be a set in $\A_e$ such that
\begin{equation}\label{Ee}
\E( \prod_{e \in H} 1_{E_e}(x) | x \in V_J) \leq \delta
\end{equation}
for some $0 < \delta < 1$.  Then for each $e \in H$ there exists a set $E'_e \in \A_e$ such that
$$ \bigcap_{e \in H} E'_e = \emptyset$$
and
$$ \E( 1_{E_e \backslash E'_e}(x) | x \in V_J) = o_{\delta \to 0; J}(1) \hbox{ for all } e \in H.$$
Furthermore, there exists sub-algebras $\B_{e'} \subseteq \A_{e'}$ whenever $e' \subset J$ and $|e'| < d$ obeying the complexity estimate
$$ |\B_{e'}| = O_{J, \delta}(1) \hbox{ whenever } e' \subseteq J \hbox{ and } |e'| < d$$
and
$$ E'_e \in \bigvee_{e' \subsetneq e} \B_{e'} \hbox{ for all } e \in H.$$
Here of course $\bigvee_{e' \subsetneq e} \B_{e'}$ is the smallest $\sigma$-algebra which contains $\B$.
\end{theorem}

\begin{remarks}
For this paper, we will only need this theorem in the special case when $d = |J|-1$ and $H$ is the simplex hypergraph $H = {J \choose |J|-1}$,
and when all the $V_j$ are equal to each other (in fact, they will all be set equal to a finite additive group $Z$).
On the other hand, this special case does not seem to be significantly easier to prove than the general case.  The hypothesis \eqref{Ee} asserts
that the sets $(E_e)_{e \in H}$ (which can be thought of as families of edges in a partite hypergraph) contain very few copies of $H$;
the hypergraph removal lemma then asserts that those copies of $H$ can be removed by replacing the edge sets $E_e$ with slightly different
edge sets $E'_e$ with bounded complexity.  For a more detailed discussion of this lemma, we refer the reader to the references given above.
\end{remarks}

At first glance, Theorem \ref{main-2} has nothing to do with Theorem \ref{sz-multi-alt}.  However, as observed in
\cite{rsz}, \cite{soly-roth}, \cite{soly-2}, \cite{frankl02}, \cite{gowers-hyper}, \cite{rodl2}, it is in fact relatively easy to deduce the former from the latter, and we include a proof below for the reader's convenience.

\begin{proof}[of Theorem \ref{sz-multi-alt} assuming Theorem \ref{main-2}]
Let us first make the ``ergodic'' hypothesis that the elements $\{ \phi_i(r) - \phi_j(r): i,j \in J; r \in Z \}$ generate $Z'$ as an additive group;
we will remove this hypothesis at the end of the argument.
Let $V = (J, (V_j)_{j \in J}, d, H)$ be the hypergraph system with $V_j := Z$, $d := |J|-1$, and $H := {J \choose d}$.  If
$e = J \backslash \{j\}$ is an element of $H$, we define the set $E_e \subseteq V_J = Z^J$ by
$$ E_e := \{ (x_i)_{i \in J} \in Z^J: \sum_{i \in J} \phi_i(x_i) - \phi_j(x_i) \in A \}.$$
Observe that the expression $\sum_{i \in J} \phi_i(x_i) - \phi_j(x_i)$ does not actually depend on $x_j$ and so $E_e \in \A_e$.

Now we compute the size of $\prod_{e \in H} 1_{E_e}$.  Let $\Phi: V_J \to Z' \times Z$ be the group homomorphism
$$ \Phi( (x_i)_{i \in J} ) := ( \sum_{j \in J} \phi_i(x_j), - \sum_{j \in J} x_j)$$
then we see from the definitions that
\begin{equation}\label{cape}
 \bigcap_{e \in H} E_e = \Phi^{-1}( \{ (a,r) \in Z' \times Z: a + \phi_j(r) \in A \hbox{ for all } j \} ).
\end{equation}
Consider the image of the group homomorphism $\Phi$.  This image contains all points of the form $(\phi_i(r)-\phi_j(r), 0)$ for $i,j \in J$ and $r \in Z$, and hence contains $Z' \times \{0\}$ by hypothesis.  It also contains all elements of the form $(-\phi_i(r), r)$ for any $r \in Z$ and $i \in J$.  Hence the image must be all of $Z' \times Z$; since $\Phi$ is a homomorphism, all the fibers $\Phi^{-1}(x,r)$ thus have the same cardinality.
We conclude
$$ \E( \prod_{e \in H} 1_{E_e}(x) | x \in V_J) = \E( \prod_{j \in J} 1_A(x + \phi_j(r)) | x \in Z'; r \in Z) \leq \delta$$
by hypothesis.  Applying Theorem \ref{main-2}, we can find $E'_e \in \A_e$ such that
$$ \bigcap_{e \in H} E'_e = \emptyset$$
and
$$ |E_e \backslash E'_e| = o_{\delta \to 0; |J|}(|V_J|) \hbox{ for all } e \in H.$$
We have additional information on the ``complexity'' of $E'_e$ but we will not need it for this argument.  

Next, from \eqref{cape} we see in particular that
$$ \Phi^{-1}(A \times \{0\}) \subseteq \bigcap_{e \in H} E_e;$$
since $\bigcap_{e \in H} E'_e = \emptyset$, we conclude that
$$ \Phi^{-1}(A \times \{0\}) \subseteq \bigcup_{e \in H} (E_e \backslash E'_e).$$
Thus by the pigeonhole principle there exists an $e = J \backslash \{j\}$ such that
$$ |(E_e \backslash E'_e) \cap \Phi^{-1}(A \times \{0\})| \geq \frac{|\Phi^{-1}(A \times \{0\})|}{|J|}.$$
The set $\Phi^{-1}(A \times \{0\})$ lives in the hyperplane $\{ (x_i)_{i \in J}: \sum_{i \in J} x_i = 0 \}$,
and in particular the projection map $\pi_e: V_J \to V_e$, which has multiplicity $|Z|$ everywhere,
is injective on $\phi^{-1}(A \times \{0\})$.  Hence we have
$$ |(E_e \backslash E'_e) \cap \Phi^{-1}(A \times \{0\})| \leq \frac{|E_e \backslash E_{e'}|}{|Z|}
= o_{\delta \to 0; |J|}(\frac{|V_J|}{|Z|}).$$
Since $\Phi$ is a surjective group homomorphism from $V_J$ to $Z' \times Z$, we have
$$ \frac{|\Phi^{-1}(A \times \{0\})|}{|V_J|} = \frac{|A|}{|Z' \times Z|} = \frac{1}{|Z|} \frac{|A|}{|Z'|}.$$
Combining these inequalities we obtain $|A| = o_{\delta \to 0; |J|}(|Z'|)$ as claimed.

To remove the ergodic hypothesis, we let $G$ be the subgroup of $Z'$ generated by
the elements $\{ \phi_i(r) - \phi_j(r): i,j \in J; r \in Z \}$.  We foliate $Z'$ into $|Z'|/|G|$ cosets of $G$.  An easy counting argument shows that
on all but $O( \sqrt{\delta} |Z'|/|G|)$ of these cosets $y+G$, we have
$$\E\left( \prod_{j \in J} 1_A(x + \phi_j(r)) \big| x \in y+G; r \in Z\right) \leq \sqrt{\delta}.$$
Applying the previous argument to each of these cosets, we conclude
$$ |A \cap (y+G)| = o_{\sqrt{\delta} \to 0; |J|}(|G|)$$
for each of these cosets.  Adding up the contributions for all of these cosets, as well as the $O( \sqrt{\delta} |Z'|/|G|)$ exceptional
cosets, we obtain $|A| = o_{\delta \to 0;|J|}(|Z'|)$ as claimed.
\end{proof}

\begin{remark}
Note in the above proof we did not need the complexity bounds on $E'_e$.
However, this fact will be important for us when we transfer this result to a weighted setting below.  The point is that the
pseudorandom weight function which we will introduce will be uniformly distributed with respect to lower order sets but not
with respect to arbitrary sets.
\end{remark}

Our proof of the number-theoretic results of this paper, and in particular Theorem \ref{main}, proceeds by a three-stage process
similar to that in \cite{gt-primes}.  Firstly, we apply the transference philosophy from \cite{gt-primes} to extend Theorem \ref{main-2}
to a relative version of that theorem, weighted by a pseudorandom system $(\nu_e)_{e \in H}$ of measures; this shall be done by following the arguments in \cite{gt-primes} closely, the main observation being that those arguments did not significantly rely on any additive structure in the underlying system and thus generalize from the ergodic system $\Z_N$ to an arbitrary hypergraph system without any fundamental new difficulties.  Next, 
by repeating the deduction of Theorem \ref{sz-multi-alt} from Theorem \ref{main-2}, we obtain a relative version of Theorem \ref{sz-multi-alt}, in which the set $A$ is measured with respect to a pseudorandom measure $\nu$; this step of the argument is quite easy.  Finally, we apply
this relative version of Theorem \ref{sz-multi-alt} to the Gaussian primes by constructing a psuedorandom majorant for these primes in the spirit
of the work of Goldston and Y{\i}ld{\i}r{\i}m (with some additional simplifications introduced in \cite{tao-gy}).

One additional technical complication which appears in this work is that the Gaussian primes (or almost primes) contain certain correlations which are not present in the rational case.  In particular, the Gaussian (almost) primes have a different density on lines such as the real line, than they do on
all of $\Z[i]$.  Also, there is an obvious correlation between $p$ being a Gaussian (almost) prime and $\overline{p}$ being a Gaussian (almost) prime.  We shall eliminate the first type of correlation by excluding the ``exceptional'' Gaussian primes whose norm is not a rational prime.  The second type of correlation cannot be eliminated so easily, but fortunately its contributions to the error terms are ultimately manageable.

The author is supported by a grant from the Packard foundation.  The author also thanks Timothy Gowers and Ben Green for some helpful
conversations, and Lilian Matthiesen for pointing out the need for a self-incommensurability hypothesis.

\section{Pseudorandomness}

Before we can state our relative versions of Theorem \ref{main-2} and Theorem \ref{sz-multi-alt}, we must introduce the notion of
a \emph{pseudorandom system of measures} $(\nu_e)_{e \in H}$ on a hypergraph system $V = (J, (V_j)_{j \in J}, d, H)$.  Strictly speaking,
the concept of pseudorandomness will not be associated with a single system of measures on a hypergraph system, but rather on
a one-parameter family of measures $(\nu_e)_{e \in H} = (\nu^{(N)}_e)_{e \in H}$ on a hypergraph system $V = V^{(N)}$, where $N$
ranges over a sequence of numbers tending to infinity (e.g. $N$ could range over the primes).  This is in order to make sense of error terms such as
$o_{N \to \infty}(1)$.  However we will usually suppress the explicit dependence of our objects on $N$, as we shall work almost exclusively with a single fixed (large) value of $N$. Indeed our notation (particularly the expectation notation) is deliberately designed to hide all factors of $N$, in order to work easily in the asymptotic regime $N \to \infty$.  The concept of a pseudorandom system is closely analogous to that of a \emph{pseudorandom measure} in
\cite{gt-primes}, where the hypergraph system was replaced by the ergodic system $\Z/N\Z$.

In the rest of this paper, we fix the finite set $J$ and the index $d$, as well as the hypergraph $H \subseteq {J \choose d}$; in particular,
these objects will not depend on the parameter $N$.  We will
allow all implicit constants in the $O()$ and $o()$ notation to depend on $J$, $d$, and $H$; indeed, since for any fixed $J$ there are only
finitely many possible values of $d$ and $H$, this is the same as requiring all implicit constants to depend on $J$.

\begin{definition}[System of measures]  We define a \emph{system of measures} $(\nu_e)_{e \in H}$ 
to be a hypergraph system $V = V^{(N)} = (J, (V_j^{(N)})_{j \in J}, d, H)$ depending on a parameter $N$ (ranging over a sequence of
numbers tending to infinity),
together with a collection of non-negative functions $\nu_e = \nu_e^{(N)}: V_e^{(N)} \to \R^+$, obeying the normalization condition
\begin{equation}\label{nue-mean}
\E( \nu_e(x_e) | x_e \in V_e ) = 1 + o_{N \to \infty}(1).
\end{equation}
We will usually suppress the dependence of $V$ and $(\nu_e)_{e \in H}$ on the parameter $N$.
\end{definition}

\begin{example}  One could set $V_j^{(N)} = \Z/N\Z$, and let $\nu_e: (\Z/N\Z)^e \to \R^+$ be a random function such that
for each $x \in (\Z/N\Z)^e$, $\nu_e(x) = \log N$ with independent probability $1/\log N$, and $\nu_e(x) = 0$ otherwise.
Then with high probability, $(\nu_e)_{e \in H}$ will be a system of measures, and it will also with high probability satisfy
the pseudorandomness conditions we shall give shortly.  For a more sophisticated example, see Example \ref{model} below.
\end{example}

\begin{remark}
Note we do not require that $\nu_e$ be bounded, or even that it obey any sort of $L^p$ type moment condition (for instance,
$\E(\nu_e(x_e)^2 | x_e \in V_e)$ need not be bounded). Indeed, for applications to number theory (or indeed to any application involving sets of Banach density zero) it is vital that we allow
these moments to be unbounded.  However, we shall shortly require that various \emph{correlations} involving the $\nu_e$ be bounded.
\end{remark}

The condition \eqref{nue-mean} is not strong enough by itself for our applications, and we must supplement it with three conditions,
the \emph{dual function condition},
the \emph{linear forms condition} and the \emph{correlation condition}.  These closely mimic the conditions of the same name in \cite{gt-primes},
(where the dual function condition and linear forms condition were combined into a single (affine-)linear forms condition), 
though there are some minor technical differences.

\begin{definition}[Discrete cube]  If $e$ is a finite set, we let $\{0,1\}^e$ be the set of all binary $e$-tuples $\omega = (\omega_j)_{j \in e}$
where each $\omega_j$ is either 0 or 1.  Observe that $\{0,1\}^e$ contains in particular the zero $e$-tuple $0^e := (0)_{j \in e}$
and the one $e$-tuple $1^e := (1)_{j \in e}$.  If $x^{(0)}_J = (x^{(0)}_j)_{j \in J}$ and
$x^{(1)}_J = (x^{(1)}_j)_{j \in J}$ are two elements of $V_J$, $e$ is a subset of $J$, and $\omega \in \{0,1\}^e$ is a binary $e$-tuple, we define $x^{(\omega)}_e \in V_e$ to be the element
$$ x^{(\omega)}_e := (x^{(\omega_j)}_j)_{j \in e}.$$
We abbreviate $x^{(0^e)}_e$ as $x^{(0)}_e$, thus
$$ x^{(0)}_e := (x^{(0)}_j)_{j \in e}$$
and define $x^{(1)}_e$ similarly.
\end{definition}

\begin{definition}[Dual function]\label{dual-def} Let $V = (J, (V_j)_{j \in J}, d, H)$ be a hypergraph system, and let $e \in H$.  If $f: V_e \to \R$ is a function,
we define its \emph{dual function} $\D_e f: V_e \to \R$ by the formula
\begin{equation}\label{defdef}
\D_e f(x^{(0)}_e) := \E( \prod_{\omega \in \{0,1\}^e: \omega \neq 0^e} f(x^{(\omega)}_e) | x^{(1)}_e \in V_e )
\end{equation}
for all $x^{(0)}_e \in V_e$.
\end{definition}

\begin{example}  If $e = \{1,2\}$, then
$$ \D_{\{1,2\}}(f)(x_1,x_2) = \E( f(x_1,x'_2) f(x'_1, x_2) f(x'_1, x'_2) | x'_1 \in V_1, x'_2 \in V_2 ).$$
The dual functions will be an indispensable tool in our analysis of the Gowers cube norms $\|f_e\|_{\Box^e}$, which we shall introduce later and which
will play a pivotal role in our arguments.
\end{example}

\begin{definition}[Dual function condition]\label{dual-c} A system of measures $(\nu_e)_{e \in H}$ on the hypergraph system $V = (J, (V_j)_{j \in J}, d, H)$ is said to obey the \emph{dual function condition} if one has the pointwise estimate
$$ \D_e(\nu_e+1)(x^{(0)}_e) = O(1)$$
for all $e \in H$ and $x^{(0)}_e \in V_e$.
\end{definition}

\begin{definition}[Linear forms condition]\label{lfc}  A system of measures $(\nu_e)_{e \in H}$ on the hypergraph system $V = (J, (V_j)_{j \in J}, d, H)$
is said to obey the \emph{linear forms condition} if one has
\begin{equation}\label{linform-special}
\E( \prod_{e \in H} \prod_{\omega \in \{0,1\}^{e}} \nu_e( x_e^{(\omega)} )^{n_{e,\omega}} | x_J^{(0)}, x_J^{(1)} \in V_J ) = 1 + o_{N \to \infty}(1)
\end{equation}
for any choice of exponents $n_{e,\omega} \in \{0,1\}$.  
\end{definition}

\begin{example}  If $J = \{1,2,3\}$, $d=2$, and $H = {J \choose 2}$, then \eqref{linform-special} asserts that 
\begin{align*}
\E( &\prod_{ij = 12,23,31} \nu_{ij}(x_i,x_j) \nu_{ij}(x_i,x'_j) \nu_{ij}(x'_i,x_j) \nu_{ij}(x'_i,x'_j) \\
&\quad | x_1,x'_1 \in V_1, x_2,x'_2 \in V_2,
x_3,x'_3 \in V_3 ) = 1 + o_{N \to \infty}(1),
\end{align*}
and similarly if one or more of the twelve factors of $\nu$ in the expectation is deleted.  
\end{example}

\begin{remark}
The condition \eqref{linform-special} can be viewed
as a fairly strong assertion of independence between the quantities 
$\nu_e( x_e^{(\omega)} )$; they in particular imply that each weight $\nu_e$ is pseudorandom in the sense of \cite{krs-hyper}, 
\cite{gowers-hyper} but are significantly stronger than those bounds alone. Note that most instances of \eqref{linform-special} are coupled expressions which involve several of the $\nu_e$ in some entangled way; it may be possible to use multiple applications of the Cauchy-Schwarz inequality to 
replace these conditions by ``pure'' pseudorandomness conditions involving each of the $\nu_e$ separately, but we have not sought to do so here.
\end{remark}

\begin{remark}
Note that the linear forms condition \eqref{linform-special} implies \eqref{nue-mean} as a special case (when all but one of the exponents $n_{e,\omega}$ is set to zero).  
However we have chosen to isolate \eqref{nue-mean} for expository reasons, to emphasize the normalized nature of the $\nu_e$.
\end{remark}

\begin{example}\label{model} A model instance of a pseudorandom system of measures, of relevance to number theory, is as follows.  Let $J = \{1,\ldots,k\}$,
let $d=k-1$, and $H = {J \choose d}$.  Let $N$ be a very large integer, and let $w = w(N)$ be a moderately large integer growing slowly with $N$
(so $1/w = o_{N \to \infty}(1)$).  Let $W = \prod_{p < w} p$ be the product of the rational primes less than $w$, and let $b_1, \ldots, b_k$
be integers in $\{0,\ldots,W-1\}$ such that $\sum_{1 \leq i \leq k}(i-j) b_i$ is coprime to $W$ for each $1 \leq j \leq k$.  For each $j \in J$, let $V_j$ be the set
$$ V_j := \{ Wn + b_j: 1 \leq n \leq N \}$$
and for each $e = J \backslash \{j\}$, let $\nu_e: V_e \to \R^+$ be the function
$$ \nu_e( (x_j)_{j \in e} ) :=  \frac{\phi(W)}{W} \Lambda( \sum_{i \in e} (i-j) x_i )$$
where $\phi(W)$ is the Euler totient function of $W$ and $\Lambda$ is the von Mangoldt function.  Then, assuming a certain strong 
form of the Hardy-Littlewood prime tuples conjecture,
this system of measures will obey the linear forms condition if $w$ is a sufficiently slowly growing function of $N$.  Of course, to verify
the prime tuples conjecture is considered to be impossible by current technology; however, by modifying the arguments in \cite{gt-primes}
one can replace the normalized von Mangoldt function $\frac{\phi(W)}{W} \Lambda$ by a slightly larger 
pseudorandom function $\nu$ (essentially a truncated divisor sum
of Goldston-Y{\i}ld{\i}r{\i}m type) for which these types of conditions can be much more easily verified.  See \cite{gt-primes}.
\end{example}

In addition to controlling dual functions and linear form expectations, we will also need to control correlations (involving only a 
single measure $\nu_e$) in which both vertices $x_i^{(0)}, x_i^{(1)}$ from a vertex set $V_i$ are fixed; this quantity then 
measures some sort of pair correlation 
between $x_{e \backslash \{j\}}^{(0)}$ and $x_{e \backslash \{j\}}^{(1)}$.  For such expressions one cannot
expect a uniform bound such as $1 + o_{N \to \infty}(1)$ or even $O(1)$, because the diagonal case $x_{e \backslash \{j\}}^{(0)} = x_{e \backslash \{j\}}^{(1)}$ will almost certainly
have an abnormally large (and unbounded) correlation.  In number theoretic applications (such as Example \ref{model}), there are a few other cases where the correlation
is expected to be abnormally large, notably when $\sum_{i \in e \backslash \{j\}} x_i^{(0)} - x_i^{(1)}$ has an extremely large number of small prime factors (e.g. if it is a ``smooth'' number).  These correlations
can become unbounded (thanks to the divergence of the Euler product $\prod_p (1 - \frac{1}{p})^{-1}$, which diverges both for rational and for Gaussian primes).  However, the correlations will still be
bounded on the average, and even have bounded moments of any given order.  More precisely, we have

\begin{definition}[Correlation condition]\label{corr-c}  A system of measures $(\nu_e)_{e \in H}$ on the hypergraph system $V = (J, (V_j)_{j \in J}, d, H)$
is said to obey the \emph{correlation condition} if we have
\begin{equation}\label{correlation}
\E\left( \E( \prod_{\omega \in \{0,1\}^e} \nu_e( x_e^{(\omega)} )^{n_{e,\omega}} | x_{j}^{(0)}, x_{j}^{(1)} \in V_j )^K \bigg|
x_{e \backslash \{j\}}^{(0)}, x_{e \backslash \{j\}}^{(1)} \in V_{e \backslash \{j\}} \right) = O_K(1)
\end{equation}
for every $e \in H$, $j \in e$, any choice of exponents $n_{e,\omega} \in \{0,1\}$, and any integer $K \geq 0$.
\end{definition}

\begin{example}  If $J = \{1,2\}$, $d=2$, and $H = {J \choose 2}$, then \eqref{correlation} with $e = \{1,2\}$ and $j=1$ asserts that 
$$
\E \left( \E( \nu_{12}(x_1,x_2) \nu_{12}(x'_1,x_2) | x_2 \in V_2 )^K | x_1, x'_1 \in V_1 \right) = O_K(1)
$$
for any $K \geq 0$, and similarly if one or both of the $\nu_{12}$ factors are deleted.  Thus the pair correlations of $\nu_{12}$ are bounded in $L^K$ for any $K$.
\end{example}

\begin{definition}[Pseudorandom system]  A system of measures $(\nu_e)_{e \in H}$ on the hypergraph system $V = (J, (V_j)_{j \in J}, d, H)$
is said to be \emph{pseudorandom} if it obeys the dual function condition, the linear forms condition and the correlation condition.
\end{definition}

The system $(1)_{e \in H}$ is a rather trivial example of a pseudorandom system of measures.  More generally, we have the following simple but handy lemma that says that the arithmetic mean of a pseudorandom system $(\nu_e)_{e \in H}$ with $(1)_{e \in H}$ is also pseudorandom:

\begin{lemma}\label{half} Let $V = (J, (V_j)_{j \in J}, d, H)$ be a hypergraph system, and let $(\nu_e)_{e \in H}$ be a system of pseudorandom measures.  Then $(\frac{1}{2} + \frac{1}{2}\nu_e)_{e \in H}$ is also a system of pseudorandom measures (perhaps with slightly different constants
in the $o()$ and $O()$ notations).
\end{lemma}

\begin{proof}  This is a reprise of \cite[Lemma 5.2]{gt-primes}.  The dual function condition follows from the pointwise estimate
$$ \D_e( \frac{1}{2} + \frac{1}{2}\nu_e + 1 ) \leq \D_e( \frac{3}{2} (\nu_e + 1 ) ) = (\frac{3}{2})^{2^d-1} \D_e( \nu_e + 1).$$
As for the linear forms and correlation conditions, from the binomial formula we have
\begin{align*}
&\prod_{e \in H} \prod_{\omega \in \{0,1\}^e} (\frac{1}{2} + \frac{1}{2} \nu_e( x_e^{(\omega)} ))^{n_{e,\omega}}\\
&= \E( \prod_{e \in H} \prod_{\omega \in \{0,1\}^e} \nu_e( x_e^{(\omega)})^{n_{e,\omega} m_{e,\omega}} | m_{e,\omega} \in \{0,1\} 
\hbox{ for all } e \in H, \omega \in \{0,1\}^e )
\end{align*}
and the claim follows by linearity of expectation.
\end{proof}

In \cite{gt-primes}, Szemer\'edi's theorem was extended via a ``transference principle''
to a relative version, weighted with a pseudorandom measure.  In this paper we shall apply the same transference principle
to extend Theorem \ref{main-2} to a relative version, which we state as follows.

\begin{theorem}[Relative hypergraph removal lemma]\label{transference}
Let $V = (J, (V_j)_{j \in J}, d, H)$ be a hypergraph system, and let $(\nu_e)_{e \in H}$ be a system of pseudorandom measures, 
For each $e \in H$, let $E_e$ be a set in $\A_e$ such that
\begin{equation}\label{nupe}
 \E( \prod_{e \in H} 1_{E_e}(x) \nu_e(\pi_e(x)) | x \in V_J ) \leq \delta
 \end{equation}
for some $0 < \delta < 1$.  Then, if $N$ is sufficiently large depending on $\delta$ and $J$,
for each $e \in H$ there exists a set $E'_e \in \A_e$ such that
$$ \bigcap_{e \in H} E'_e = \emptyset$$
and
\begin{equation}\label{eex}
 \E( 1_{E_e \backslash E'_e}(x) \nu_e(\pi_e(x)) | x \in V_J ) = o_{\delta \to 0}(1) + o_{N \to \infty; \delta}(1) \hbox{ for all } e \in H.
 \end{equation}
 Recall that all constants are allowed to depend on $J$.
Furthermore, there exists a $\sigma$-algebra $\B_{e'} \subseteq \A_{e'}$ for all $e'$ with $|e'| < d$ such that
$$ |\B_{e'}| = O_{\delta}(1) \hbox{ whenever } e' \subseteq J \hbox{ and } |e'| < d$$
and
$$ E'_e \in \bigvee_{e' \subsetneq e} \B_{e'} \hbox{ for all } e \in H.$$
\end{theorem}

The proof of Theorem \ref{transference} is lengthy and shall occupy Sections \ref{cube-sec}-\ref{structure-sec}.
Just as Theorem \ref{main-2} implies Theorem \ref{sz-multi-alt}, Theorem \ref{transference} implies the following relative version
of Theorem \ref{sz-multi-alt}.

\begin{theorem}[Relative multidimensional Szemer\'edi's theorem]\label{sz-multi-alt-relative} Let $Z, Z'$ be two finite additive groups,  and let $(\phi_j)_{j \in J}$ be a finite collection of group homomorphisms $\phi_j: Z \to Z'$ be any group homomorphisms from $Z$ to $Z'$.  We assume the ergodic hypothesis that the elements $\{ \phi_i(r) - \phi_j(r): i,j \in J, r \in Z\}$ generate $Z'$ as an abelian group.  Let $\nu: Z' \to \R^+$ be
a non-negative function, with the property that in the hypergraph system $(J, (V_j)_{j \in J}, d,H)$ with $V_j := Z$, $d := |J|-1$, $H := {J \choose d}$, the
collection $(\nu_e)_{e \in H}$ defined by
$$ \nu_{J \backslash \{j\}}((x_i)_{i \in J}) := \nu(\sum_{i \in J} \phi_i(x_i) - \phi_j(x_i))$$
is a pseudorandom family of measures.  Then if $A$ is a subset of $Z'$ such that
\begin{equation}\label{anu}
\E( \prod_{j \in J} 1_A(x + \phi_j(r)) \nu(x + \phi_j(r)) | x \in Z'; r \in Z) \leq \delta
\end{equation}
for some $0 < \delta \leq 1$, then we have 
$$ \E( 1_A(x) \nu(x) | x \in Z') = o_{\delta \to 0;|J|}(1) + o_{N \to \infty; \delta}(1).$$
\end{theorem}

\begin{proof}[of Theorem \ref{sz-multi-alt-relative} assuming Theorem \ref{transference}]  This shall be a reprise of the proof of
Theorem \ref{sz-multi-alt}.  We may assume $N$ is large since the claim is trivial otherwise.
As in that proof, we define the set $E_e \in \A_e$ for each element
$e = J \backslash \{j\}$ of $H$ as
$$ E_e := \{ (x_i)_{i \in J} \in Z^J: \sum_{i \in J} \phi_i(x_i) - \phi_j(x_i) \in A \},$$
and we recall the group homomorphism $\Phi: V_J \to Z' \times Z$ defined by
$$ \Phi( (x_i)_{i \in J} ) := ( \sum_{j \in J} \phi_i(x_j), - \sum_{j \in J} x_j).$$
Then we have
$$
 \prod_{e \in H} 1_{E_e}(x) \nu_e(\pi_e(x)) = \prod_{j \in J} 1_A(a + \phi_j(r)) \nu(a + \phi_j(r))
$$
for all $x \in V_J$, where $(a,r) := \Phi(x)$.  From the ergodic hypothesis, $\Phi$ is surjective, and hence all the fibers
$\Phi^{-1}(a,r)$ have the same cardinality.  Thus
$$ \E( \prod_{e \in H} 1_{E_e}(x) \nu_e(\pi_e(x)) | x \in V_J ) = \E( \prod_{j \in J} 1_A(a + \phi_j(r)) \nu(a + \phi_j(r)) ) \leq \delta$$
by hypothesis.  Applying Theorem \ref{transference} (for $N$ large enough), we can find $E'_e \in \A_e$ such that
\begin{equation}\label{ee-empty}
 \bigcap_{e \in H} E'_e = \emptyset
\end{equation}
and
\begin{equation}\label{ee-rare}
\E( 1_{E_e \backslash E'_e}(x) \nu_e(\pi_e(x)) | x \in V_J ) = o_{\delta \to 0; |J|}(1) + o_{N \to \infty; \delta,|J|}(1) \hbox{ for all } e \in H.
\end{equation}
Once again, we will not need to use the additional complexity information on $E'_e$.

Next, we observe from the definition of $E_e$ that
$$ 1_A(a) 1_{r=0} = 1_{r=0} \prod_{e \in H} 1_{E_e}(x)$$
 for all $x \in V_J$, where $(a,r) = \Phi(x)$ as before.  From \eqref{ee-empty} we conclude that
$$ 1_A(a) 1_{r=0} \leq \sum_{e \in H} 1_{r=0} 1_{E_e \backslash E'_e}(x).$$
Multiplying by $\nu(a)$, averaging in $x$, and then applying the pigeonhole principle, there exists
an $e = J \backslash \{j\}$ in $H$ such that
$$ \frac{1}{|J|} \E( 1_A(a) 1_{r=0} \nu(a) | x \in V_J ) \leq \E( 1_{r=0} \nu(a) 1_{E_e \backslash E'_e}(x) | x \in V_J ).$$
Observe that $\nu(a) = \nu_e(\pi_e(x))$.  Also, recall that the fibers $\Phi^{-1}(a,r)$ all have equal cardinality.  Thus we have
$$ \frac{1}{|J|} \E( 1_A(a) 1_{r=0} \nu(a) | (a,r) \in Z' \times Z ) \leq \E( 1_{r=0} \nu_e(\pi_e(x)) 1_{E_e \backslash E'_e}(x) | x \in V_J ).$$
Since the function $\nu_e(\pi_e(x)) 1_{E_e \backslash E'_e}(x)$ does not depend on the $x_j$ variable, and that the constraint $r=0$ forces $x_j$ to be determined by all the other variables, we have
$$ \E( 1_{r=0} \nu_e(\pi_e(x)) 1_{E_e \backslash E'_e}(x) | x \in V_J ) = \frac{1}{|Z|} \E( \nu_e(\pi_e(x)) 1_{E_e \backslash E'_e}(x) | x \in V_J ).$$
Also, we have $\E( 1_A(a) 1_{r=0} \nu(a) | (a,r) \in Z' \times Z ) = \frac{1}{|Z|} \E( 1_A(a) \nu(a) | a \in Z' )$.  Thus
$$ \E( 1_A(a) \nu(a) | a \in Z' ) \leq |J| \E( \nu_e(\pi_e(x)) 1_{E_e \backslash E'_e}(x) | x \in V_J )$$ 
and the claim follows from \eqref{ee-rare}.
\end{proof}

\begin{remarks} The complexity bound was not used in this argument, however we will need the complexity bound from Theorem \ref{main-2} in order to successfully transfer that theorem to the relative setting.  The ergodic hypothesis can be dropped by foliating $Z'$ into cosets as in the proof of Theorem \ref{sz-multi-alt}, but one has to modify the pseudorandomness hypotheses on $\nu$ accordingly; we omit the details.
\end{remarks}

\begin{remark} Theorem \ref{sz-multi-alt-relative} can be used to prove a slight variant of the relative Szemer\'edi theorem in
\cite[Theorem 3.5]{gt-primes} (with some minor variations in the linear forms and correlation condition).  
This is unsurprising given that the proof of Theorem \ref{sz-multi-alt-relative} given here closely follows
the proof of that theorem in \cite{gt-primes}.
\end{remark}

In the next few sections we shall prove Theorem \ref{transference}, and hence Theorem \ref{sz-multi-alt-relative}.  In the second half of the
paper (from Section \ref{Gaussian-sec} onwards) we shall apply Theorem \ref{sz-multi-alt-relative} to questions concerning the primes and Gaussian primes.

\section{The Gowers cube norm, and overview of proof of Theorem \ref{transference}}\label{cube-sec}

In this section we shall recall the \emph{Gowers cube norm} $\|f\|_{\Box^e}$, which shall be a fundamental tool in our proof of
Theorem \ref{transference}, playing a role closely analogous to that of the Gowers uniformity norm $\|f\|_{U^d}$ in \cite{gt-primes}.
We will then use this norm to split the proof of Theorem \ref{transference} into four components.  One component is a weighted version (Theorem \ref{weight-main}) of the hypergraph removal lemma, which is a minor generalization of Theorem \ref{main-2}.  Another component will be
a \emph{generalized von Neumann theorem} (Theorem \ref{gvn}), which essentially asserts that functions with small cube norm have a negligible impact on
the quantity \eqref{nupe}.  A third component is a \emph{structure theorem}, which decomposes the function $1_{E_e} \nu_e$ into
a bounded non-negative function (which can be dealt with using Theorem \ref{main-2}) and a remainder with small cube norm (which can
be dealt with using Theorem \ref{gvn}), plus a negligible error.  Finally (and this is where we need the complexity information from Theorem \ref{main-2}), we need a simple result (Corollary \ref{loword-neg}) which asserts that functions with small cube norm are uniformly distributed with respect to lower order sets.

We now turn to the details.  We begin by defining the Gowers cube norm.

\begin{definition}[Gowers cube norm]  Let $V = (J, (V_j)_{j \in J}, d, H)$ be a hypergraph system, let $e$ be an element of $H$, and let
$f: V_e \to \R$ be a function.  We define the \emph{Gowers cube norm} $\|f\|_{\Box^e}$ of $f$ to be the quantity
$$ \|f\|_{\Box^e} := \E( \prod_{\omega \in \{0,1\}^e} f(x_e^{(\omega)}) | x_e^{(0)}, x_e^{(1)} \in V_e )^{1/2^{|e|}}.$$
\end{definition}

\begin{examples}  If $e$ is empty, $e = \emptyset$, then $V_e$ is a singleton set, and $\|f\|_{\Box^e}$ is simply equal to the single value of $f$ on $V_e$; in particular $\|f\|_{\Box^e}$ can be negative in this case.  If $e$ is a point, thus $e = \{j\}$, then
$$ \|f\|_{\Box^e} = \E( f(x_e^{(0)}) f(x_e^{(1)}) | x_e^{(0)}, x_e^{(1)} \in V_e)^{1/2} = |\E( f(x) | x \in V_j )|.$$
In particular, the $\Box^e$ ``norm'' is only a semi-norm in this case.
If $e$ consists of two points, thus $e = \{i,j\}$, then
\begin{align*}
\|f\|_{\Box^e} 
&= \E( f(x_e^{(0,0)}) f(x_e^{(0,1)}) f(x_e^{(1,0)}) f(x_e^{(1,1)}) | x^{(0)}_e, x^{(1)}_e \in V_e )^{1/4} \\
&= \E( f(x_i, x_j) f(x_i, x'_j), f(x'_i, x_j) f(x'_i, x'_j) | x_i, x'_i \in V_i; x_j, x'_j \in V_j)^{1/4} \\
&= \E( \E( f(x_i, x_j) f(x'_i, x_j) | x_j \in V_j)^2 | x_i, x'_i \in V_i )^{1/4}.
\end{align*}
Thus $\|f\|_{\Box^e}$ is non-negative (and one can easily verify that it vanishes if and only if $f$ is identically zero).
In this case one can view $f$ as the kernel of a linear operator $T$ from $V_i$ to $V_j$, and $\|f\|_{\Box^e}$ can be viewed as square root of the normalized Hilbert-Schmidt norm of $T^* T$, or as the 4-Schatten norm $\tr( TT^* TT^*)^{1/4}$.  Alternatively, one can view $f$ as a weighted bipartite graph from $V_i$ to $V_j$, and then $\|f\|_{\Box^e}$ is a normalized count of the 4-cycles in this graph, weighted by $f$.
\end{examples}

\begin{example}  Suppose $V_j = Z$ for some abelian group, $e \in H$, and $f: V_e \to \R$ has the special form
$$ f( (x_j)_{j \in e} ) = F( \sum_{j \in e} x_j )$$
for some function $F: Z \to \R$.  Then $\|f\|_{\Box^e} = \|F\|_{U^d(Z)}$, where $d = |e|$ and the $U^d(Z)$ norm is the Gowers uniformity
norm, defined for instance in \cite{gowers}, \cite{gt-primes}, \cite{tao:ergodic}.
\end{example}

\begin{remark}
The $\Box^e$ norm is closely related to the concept of a dual function $\D_e(f)$, see \eqref{fdef} below.  Indeed, just as in \cite{gt-primes},
the complementarity between Gowers uniform functions - that is, functions with small $\Box^e$ norm - and Gowers anti-uniform functions (specifically,
functions generated by dual functions) will lie at the heart of the transference principle that underlies this paper.
\end{remark}

If $|e| \geq 1$, and we split $e = e' \cup \{j\}$ for an arbitrary $j \in e$, where $e' := e' \backslash \{j\}$, then one can verify the identity
\begin{equation}\label{fbe}
\| f \|_{\Box^e} = \E( (\prod_{\omega \in \{0,1\}^{e'}} f(x_{e'}^{(\omega)}, x_j) | x_j \in V_j)^2
 | x_{e'}^{(0)}, x_{e'}^{(1)} \in V_{e'} )^{1/2^{|e|}}
\end{equation}
and thus $\|f\|_{\Box^e}$ is non-negative.  One can also verify that $\Box^e$ obeys the triangle inequality when $|e| \geq 1$ (see e.g. \cite{tao:survey}) but we will not need this fact here.  One further consequence of the identity \eqref{fbe} is that
$$ \|fg\|_{\Box^e} \leq \|f\|_{\Box^e}$$
whenever $g: V_e \to [-1,1]$ is a bounded function which is independent of the $x_j$ variable for some $j \in e$.  In particular,
$g$ can be a indicator function.  Iterating this claim, we obtain

\begin{corollary}\label{loword}  Let $V = (J, (V_j)_{j \in J},d,H)$ be a hypergraph system, and let $e \in H$.
Let $f: V_e \to \R$ be a function, and for each $e' \subsetneq e$ let $E_{e'}$ be a subset of $V_{e'}$.  Then we have
$$ |\E( f( x_e ) \prod_{e' \subsetneq e} 1_{E_{e'}}( x_{e'} ) | x_e \in V_e )| \leq \|f\|_{\Box^e},$$
where $x_{e'}$ is the restriction of $x_e$ to $V_{e'}$ (thus if $x_e = (x_j)_{j \in e}$ then $x_{e'} = (x_j)_{j \in e'}$).
\end{corollary}

In particular, we have the following result, which is one of four ingredients necessary to prove Theorem \ref{transference}.  It asserts
that Gowers uniform functions are uniformly distributed across lower order sets - sets which arise from the $\sigma$-algebras $\A_{e'}$
with $e'$ strictly smaller than $e$.

\begin{corollary}[Gowers uniform functions are orthogonal to lower order sets]\label{loword-neg}  
Let $V = (J, (V_j)_{j \in J}, d, H)$ be a hypergraph system, and let $(\nu_e)_{e \in H}$ be a system of pseudorandom measures.
Suppose there exists sub-algebras $\B_{e'} \subseteq \A_{e'}$ whenever $e' \subset J$ and $|e'| < d$ obeying the complexity estimate
$$ |\B_{e'}| \leq M \hbox{ whenever } e' \subseteq J \hbox{ and } |e'| < d$$
for some $M$.  For each $e \in H$, let $E'_e$ be a set in $\bigvee_{e' \subsetneq e} \B_{e'}$.  Then we have
$$ \E( 1_{E'_e}(x) f(\pi_e(x)) | x \in V_J ) = O_M( \|f\|_{\Box^e} )$$
for any $f: V_e \to \R$.
\end{corollary}

\begin{proof}  We can decompose $E'_e$ as the union of $O_M(1)$ atoms of $\bigvee_{e' \subsetneq e} \B_{e'}$, each of which are in turn the intersection of atoms from $\B_{e'}$.  By the triangle inequality, it thus suffices to show that
$$ |\E( f(\pi_e(x)) \prod_{e' \subseteq e} 1_{F_{e'}}(x) | x \in V_J )| \leq \|f\|_{\Box^e} $$
whenever $F_{e'} \in \B_{e'}$.  But this follows from Corollary \ref{loword} after eliminating the redundant averaging over
those variables $x_j$ for which $j \in J \backslash e$.
\end{proof}

The second ingredient we need to prove Theorem \ref{transference} is the following minor generalization of Theorem \ref{main-2}, which does not
involve a pseudorandom system of measures, but replaces the sets $E_e$ by bounded weight functions.

\begin{theorem}[Weighted hypergraph removal lemma]\label{weight-main}
Let $V = (J, (V_j)_{j \in J}, d, H)$ be a hypergraph system.  For each $e \in H$, let $f_e: V_e \to [0,1]$ be a bounded non-negative function
\begin{equation}\label{fide}
 \E( \prod_{e \in H} f_e(\pi_e(x)) | x \in V_J)| \leq \delta
\end{equation}
for some $0 < \delta < 1$.  Then for each $e \in H$ there exists a set $E'_e \in \A_e$ such that
\begin{equation}\label{ehe}
\bigcap_{e \in H} E'_e = \emptyset
\end{equation}
and
\begin{equation}\label{fepe}
 \E( f_e(\pi_e(x)) 1_{V_J \backslash E'_e}(x) | x \in V_J) = o_{\delta \to 0}(1) \hbox{ for all } e \in H.
\end{equation}
Furthermore, there exists sub-algebras $\B_{e'} \subseteq \A_{e'}$ whenever $e' \subset J$ and $|e'| < d$ obeying the complexity estimate
$$ |\B_{e'}| = O_{\delta}(1) \hbox{ whenever } e' \subseteq J \hbox{ and } |e'| < d$$
and
$$ E'_e \in \bigvee_{e' \subsetneq e} \B_{e'} \hbox{ for all } e \in H.$$
\end{theorem}

Note that Theorem \ref{main-2} is the special case of Theorem \ref{weight-main} in the case when the $f_e$ are indicator functions.

\begin{proof}  For each $e \in H$, let $E_e \subseteq V_J$ be the set
$$ E_e := \{ x \in V_J: f_e(\pi_e(x)) \geq \delta^{\frac{1}{2|H|}} \}.$$
Clearly, $E_e \in \A_e$.  From \eqref{fide} we see that
$$ \E( \prod_{e \in H} 1_{E_e}(x) | x \in V_j)| \leq \delta^{1/2}.$$
Applying Theorem \ref{main-2}, we obtain a set $E'_e \in \A_e$ for each $e \in H$ obeying \eqref{ehe} and the desired complexity bounds,
and such that
$$ \E( 1_{E_e}(x) 1_{V_J \backslash E'_e}(x) | x \in V_J) = o_{\delta \to 0}(1) \hbox{ for all } e \in H.$$
Using the pointwise estimate $f_e(\pi_e(x)) \leq 1_{E_e}(x) + \delta^{\frac{1}{2|H|}}$, we obtain \eqref{fepe}, and the claim follows.
\end{proof}

The third ingredient of the proof of Theorem \ref{transference} is the following generalized von Neumann theorem, which we prove in Section \ref{gvn-sec}.  It asserts that Gowers uniform functions have a negligible impact on averages such as those appearing in \eqref{anu}, even
when such functions are bounded by a pseudorandom system of measures rather than by 1.

\begin{theorem}[Generalized von Neumann theorem]\label{gvn}  Let $V = (J, (V_j)_{j \in J}, d, H)$ be a hypergraph system, and let $(\nu_e)_{e \in H}$ be a system of
pseudorandom measures on $V$.  For every $e \in H$, let $f_e: V_e \to \R$ be a function such that we have the pointwise estimates
\begin{equation}\label{fne}
 |f_e(x_e)| \leq \nu_e(x_e) \hbox{ for all } x_e \in V_e \hbox{ and } e \in H.
 \end{equation}
Then we have
$$ \E( \prod_{e \in H} f_e(\pi_e(x)) | x \in V_J ) = O( \inf_{e \in H} \| f_e \|_{\Box^e} ) + o_{N \to \infty}(1).$$
\end{theorem}

This theorem will follow from multiple applications of the Cauchy-Schwarz inequality; the main difficulty is that of setting up a notational system which is not too cumbersome in order to track all the variables.  It is the analogue of \cite[Proposition 5.3]{gt-primes}.

The final ingredient in the proof of Theorem \ref{transference} is the following structure theorem, which is the analogue of \cite[Proposition 8.1]{gt-primes}.  It splits an arbitrary system of functions (bounded by a pseudorandom system) into a bounded component, plus a Gowers uniform component, outside of a set of negligible measure.

\begin{theorem}[Structure theorem]\label{structure}
Let $V = (J, (V_j)_{j \in J}, d, H)$ be a hypergraph system, and let $(\nu_e)_{e \in H}$ be a system of
pseudorandom measures on $V$.  Let $e \in H$, and let $f_e: V_e \to \R^+$ be a non-negative function such that we have the pointwise estimate
\begin{equation}\label{fnue}
 0 \leq f_e(x_e) \leq \nu_e(x_e) \hbox{ for all } x_e \in V_e.
 \end{equation}
Let $0 < \eps \ll 1$ be a small parameter, and assume $N$ sufficiently large depending on $\eps$.  Then there exists a $\sigma$-algebra $\B_e$ on $V_e$ and	an exceptional set $\Omega_e \in \B_e$ obeying the smallness condition
\begin{equation}\label{omegasmall}
\E( 1_{\Omega_e}(x_e) \nu_e(x_e) | x_e \in V_e ) = o_{N \to \infty; \eps}(1)
\end{equation}
and such that $\nu_e$ is uniformly distributed outside of $\Omega_e$:
\begin{equation}\label{nu-distributed}
\E(\nu_e|\B_e)(x_e) = 1 + o_{N \to \infty; \eps}(1) \hbox{ for all } x \in V_e \backslash \Omega_e.
\end{equation}
Furthermore, we have the uniformity estimate
\begin{equation}\label{nu-cube}
\| (1 - 1_{\Omega_e}) (f_e - \E(f_e | \B_e)) \|_{\Box^e} \leq \eps^{1/2^{|J|}}.
\end{equation}
\end{theorem}

The proof of this theorem is somewhat lengthy and will occupy Sections \ref{dual-sec}-\ref{structure-sec}.  Assuming both Theorem \ref{gvn} and Theorem \ref{structure}, we can now
combine all the above ingredients to prove Theorem \ref{transference} (and hence Theorem \ref{sz-multi-alt-relative}.

\begin{proof}[of Theorem \ref{transference} assuming Theorems \ref{gvn}, \ref{structure}]
Let $V$, $(\nu_e)_{e \in H}$, $(E_e)_{e \in H}$ be as in Theorem \ref{transference}.  Since $E_e \in \A_e$, we can write $E_e = \pi_e^{-1}(F_e)$ for some set
$F_e \subseteq V_e$.  Let $0 < \eps < \delta^{2^{|J|}}$ be a small parameter (depending on $\delta$, of course) to be chosen later.  We may assume that $N$
is large depending on $\eps$ and $\delta$ as the claim is trivial otherwise.
Applying Theorem \ref{structure} once for each $e \in H$ with $f_e := 1_{F_e} \nu_e$, we can find $\sigma$-algebras $\B_e$ on $V_e$
and sets $\Omega_e \in \B_e$ obeying \eqref{omegasmall}, \eqref{nu-distributed}, \eqref{nu-cube}.  

Now write $f_{e,\Box} := (1 - 1_{\Omega_e}) (f_e - \E(f_e | \B_e))$ and $f_{e,\Box^\perp} := (1 - 1_{\Omega_e}) \E(f_e | \B_e)$, thus
$f_{e,\Box}$ and $f_{e,\Box^\perp}$ are real-valued functions on $V_e$ which add up to $(1 - 1_{\Omega_e}) f_e$, which is of course
bounded by $f_e = 1_{F_e} \nu_e$.
From \eqref{nu-distributed}, \eqref{nu-cube} we have the estimates
\begin{align*}
\| f_{e,\Box} \|_{\Box^e} &\leq \eps^{1/2^{|J|}} \\
0 \leq f_{e,\Box^\perp}(x_e) &\leq 1 + o_{N \to \infty; \eps}(1) \hbox{ for all } x_e \in V_e\\
|f_{e,\Box}(x_e)| &\leq \nu_e(x_e) + 1 + o_{N \to \infty; \eps}(1) \hbox{ for all } x_e \in V_e \\
0 \leq f_{e,\Box}(\pi_e(x)) + f_{e,\Box^\perp}(\pi_e(x)) &\leq 1_{E_e}(x) \nu_e(\pi_e(x)) \hbox{ for all } x \in V_J.
\end{align*}
Thus we have split $f_e$ (modulo a negligible error) into a bounded component $f_{e,\Box^\perp}$, and a component $f_{e,\Box}$ with small
$\Box^e$ norm.  From the latter estimate and \eqref{nupe} we have
$$  \E( \prod_{e \in H} (f_{e,\Box}(\pi_e(x)) + f_{e,\Box^\perp}(\pi_e(x))) | x \in V_J ) \leq \delta.$$
We split the left-hand side into $2^{|H|} = O(1)$ terms in the obvious manner.  All but one of these terms involves at least
one function $f_{e,\Box}$.  Applying Theorem \ref{gvn} (using Lemma \ref{half} to replace $\nu_e$ by $\frac{1}{2} + \frac{1}{2} \nu_e$, and
scaling by the harmless factor $2 + o_{N \to \infty; \eps}(1)$) and
the above estimates, we see that the contribution of each such term is $O(\eps^{1/2^{|J|}})$ (if $N$ is sufficiently large depending on $\eps$).
By the triangle inequality, we thus conclude
$$  \E( \prod_{e \in H} f_{e,\Box^\perp}(\pi_e(x)) | x \in V_J ) \leq \delta + O(\eps^{1/2^{|J|}}) = O(\delta)$$
since we are taking $0< \eps < \delta^{2^{|J|}}$.  We can now apply Theorem \ref{weight-main} (with $f_e$ replaced by $f_{e,\Box^\perp}$)
to obtain sets $E'_e \in \A_e$ for each $e \in H$ obeying \eqref{ehe} and
\begin{equation}\label{fepeu}
 \E( f_{e,\Box^\perp}(\pi_e(x)) 1_{V_J \backslash E'_e}(x) | x \in V_J) = o_{\delta \to 0}(1) \hbox{ for all } e \in H.
 \end{equation}
Furthermore, there exists sub-algebras $\B_{e'} \subseteq \A_{e'}$ whenever $e' \subset J$ and $|e'| < d$ obeying the complexity estimate
$$ |\B_{e'}| = O_{\delta}(1) \hbox{ whenever } e' \subseteq J \hbox{ and } |e'| < d$$
and
$$ E'_e \in \bigvee_{e' \subsetneq e} \B_{e'} \hbox{ for all } e \in H.$$
The only remaining thing to establish is \eqref{eex}. Applying Corollary \ref{loword-neg} we obtain
$$  \E( f_{e,\Box}(\pi_e(x)) 1_{V_J \backslash E'_e}(x) | x \in V_J) = O_\delta(\|f_{e,\Box}\|_{\Box^e}) = O_\delta(\eps) \hbox{ for all } e \in H.$$
Adding this to \eqref{fepeu} we conclude
$$
 \E\left( (1 - 1_{\Omega_e}(\pi_e(x))) f_e(\pi_e(x)) 1_{V_J \backslash E'_e}(x) \bigg| x \in V_J\right) = o_{\delta \to 0}(1)\hbox{ for all } e \in H
$$
if $\eps$ is sufficiently small depending on $\delta$ (and $N$ is sufficiently large depending on $\eps$).  Thus we have
$$
 \E\left( (1 - 1_{\Omega_e}(\pi_e(x)) \nu_e(\pi_e(x)) 1_{E_e \backslash E'_e}(x) \bigg| x \in V_J\right) = o_{\delta \to 0}(1) \hbox{ for all } e \in H.
$$
From this and \eqref{omegasmall} we have \eqref{eex} as desired.
\end{proof}

It now remains to prove Theorem \ref{gvn} and Theorem \ref{structure}, which we shall do in the next few sections.  The proofs of these theorems
can be read independently of each other.

\section{A generalized von Neumann theorem}\label{gvn-sec}

The purpose of this section is to prove Theorem \ref{gvn}.  We shall follow the proof of \cite[Proposition 5.3]{gt-primes} closely.  
The basic idea is to repeatedly use the Cauchy-Schwarz inequality to
replace each of the $f_e$ factors by a $\nu_e$ in turn, until only one function $f_e$ remains.
The key estimate for doing so is the following:

\begin{proposition}[Cauchy-Schwarz]\label{ics}  
Let $V = (J, (V_j)_{j \in J}, d, H)$ be a hypergraph system, and let $(\nu_e)_{e \in H}$ be a system of
pseudorandom measures on $V$.  For every $e \in H$, let $f_e: V_e \to \R$ be a function such that we have the pointwise estimates
\eqref{fne}.  For any set $J' \subseteq J$, let $Q_{J'}$ denote the quantity
\begin{align*}
Q_{J'} &:= \E\biggl( \prod_{e \in H: J' \subseteq e} \prod_{\omega \in \{0,1\}^{J'}} f_e( x^{(\omega)}_e ) \\
&\quad \prod_{e \in H: J' \not \subseteq e} 
\prod_{\omega \in \{0,1\}^{e \cap J'}} \nu_e( x^{(\omega)}_e )
| x^{(0)}_J, x^{(1)}_J \in V_J; x^{(0)}_{J \backslash J'} = x^{(1)}_{J \backslash J'} \biggr)
\end{align*}
where we extend $\omega$ arbitrarily from $J'$ or $e \cap J'$ to $e$ (the exact choice of extension is unimportant since
$x^{(0)}_{J \backslash J'} = x^{(1)}_{J \backslash J'}$).  Then for any $J' \subsetneq J$ and $j_0 \in J \backslash J'$, we have
$$ |Q_{J'}| \leq (1 + o_{N \to \infty}(1)) |Q_{J' \cup \{j_0\}}|^{1/2}.$$
\end{proposition}

\begin{example} If $J = \{1,2,3\}$ and $H = {J \choose 2}$, then
\begin{align*}
Q_\emptyset &= \E( f_{\{1,2\}}(x_1,x_2) f_{\{2,3\}}(x_2,x_3) f_{\{3,1\}}(x_3,x_1) | x_1 \in V_1, x_2 \in V_2, x_3 \in V_3) \\
Q_{\{1\}} &= \E( f_{\{1,2\}}(x_1,x_2) f_{\{1,2\}}(x'_1,x_2) \nu_{\{2,3\}}(x_2,x_3) f_{\{3,1\}}(x_3,x_1) f_{\{3,1\}}(x_3,x'_1) \\
&\quad |
x_1,x'_1 \in V_1, x_2 \in V_2, x_3 \in V_3 ) \\
Q_{\{1,2\}} &= \E( f_{\{1,2\}}(x_1,x_2) f_{\{1,2\}}(x'_1,x_2) f_{\{1,2\}}(x_1,x'_2) f_{\{1,2\}}(x'_1,x'_2) \\
&\quad \nu_{\{2,3\}}(x_2,x_3) \nu_{\{2,3\}}(x'_2,x_3)
\nu_{\{3,1\}}(x_3,x_1) \nu_{\{3,1\}}(x_3,x'_1) \\
&\quad | x_1, x'_1 \in V_1, x_2,x'_2 \in V_2, x_3 \in V_3 ).
\end{align*}
\end{example}

\begin{proof}  
For all pairs $(x^{(0)}_J, x^{(1)}_J) \in V_J \times V_J$ with $x^{(0)}_{J \backslash J'} = x^{(1)}_{J \backslash J'}$, 
let us define the functions
\begin{align*}
F(x^{(0)}_J, x^{(1)}_J) &:= \prod_{e \in H: J' \subseteq e; j_0 \in e} \prod_{\omega \in \{0,1\}^{J'}} f_e( x^{(\omega)}_e )  \\
G(x^{(0)}_J, x^{(1)}_J) &:= \prod_{e \in H: J' \subseteq e; j_0 \not \in e} \prod_{\omega \in \{0,1\}^{J'}} f_e( x^{(\omega)}_e )  \\
K(x^{(0)}_J, x^{(1)}_J) &:= 
\prod_{e \in H: J' \not \subseteq e; j_0 \in e} \prod_{\omega \in \{0,1\}^{e \cap J'}} \nu_e( x^{(\omega)}_e ) \\
L(x^{(0)}_J, x^{(1)}_J) &:= 
\prod_{e \in H: J' \not \subseteq e; j_0 \not \in e} \prod_{\omega \in \{0,1\}^{e \cap J'}} \nu_e( x^{(\omega)}_e ) \\
M(x^{(0)}_J, x^{(1)}_J) &:= 
\prod_{e \in H: j_0 \not \in e} \prod_{\omega \in \{0,1\}^{e \cap J'}} \nu_e( x^{(\omega)}_e ), 
\end{align*}
then we can write
$$ Q_{J'} = \E( FGKL(x^{(0)}_J, x^{(1)}_J) | x^{(0)}_J, x^{(1)}_J \in V_J;
x^{(0)}_{J \backslash J'} = x^{(1)}_{J \backslash J'} ).$$
Write $J^* := J \backslash \{j_0\}$.  Currently, we are averaging over a pair $(x^{(0)}_J, x^{(1)}_J)$ in $V_J \times V_J$
with $x^{(0)}_{J \backslash J'} = x^{(1)}_{J \backslash J'}$.  But this is equivalent to averaging over a pair
$(x^{(0)}_{J^*}, x^{(1)}_{J^*})$ in $V_{J^*} \times V_{J^*}$ with $x^{(0)}_{J^* \backslash J'} = x^{(1)}_{J^* \backslash J'}$, together
with an element $x_{j_0} \in V_{j_0}$, with the understanding that $x^{(0)}_{j_0} = x^{(1)}_{j_0} = x_{j_0}$.  If one performs this change of
variables, then the functions $G$ and $L$ become independent of $x_{j_0}$.  Thus we can write $Q_{J'}$
(with a slight abuse of notation) as
$$ \E( \E( FK(x^{(0)}_{J^*}, x^{(1)}_{J^*}, x_{j_0}) | x_{j_0} \in V_{j_0} ) GL(x^{(0)}_{J^*}, x^{(1)}_{J^*}) | x^{(0)}_{J^*}, x^{(1)}_{J^*}
\in V_{J^*}; x^{(0)}_{J^* \backslash J'} = x^{(1)}_{J^* \backslash J'}).$$
By the hypothesis \eqref{fne}, we have $|G(q^*)| L(q^*) \leq M(q^*)$.
Applying Cauchy-Schwarz, we thus have
$$ |Q_{J'}| \leq X^{1/2} Y^{1/2} $$
where
$$ X := \E( |\E( FK(x^{(0)}_{J^*}, x^{(1)}_{J^*}, x_{j_0}) | x_{j_0}\in V_{j_0})|^2 M(x^{(0)}_{J^*}, x^{(1)}_{J^*}) | x^{(0)}_{J^*}, x^{(1)}_{J^*}
\in V_{J^*}; x^{(0)}_{J^* \backslash J'} = x^{(1)}_{J^* \backslash J'} )$$
and
$$Y := \E( M(x^{(0)}_{J^*}, x^{(1)}_{J^*}) | x^{(0)}_{J^*}, x^{(1)}_{J^*} \in V_{J^*}; x^{(0)}_{J^* \backslash J'} = x^{(1)}_{J^* \backslash J'} ).$$
From the linear forms condition \eqref{linform-special} we have
$$ Y = 1 + o_{N \to \infty}(1).$$
On the other hand, we can expand $X$ as
\begin{align*}
 \E\bigl( &FK(x^{(0)}_{J^*}, x^{(1)}_{J^*}, x_{j_0}^{(0)}) FK(x^{(0)}_{J^*}, x^{(1)}_{J^*}, x_{j_0}^{(1)}) 
M(x^{(0)}_{J^*}, x^{(1)}_{J^*}) \\
&| x^{(0)}_{J^*}, x^{(1)}_{J^*} \in V_{J^*}; x^{(0)}_{J^* \backslash J'} = x^{(1)}_{J^* \backslash J'}; 
x_{j_0}^{(0)}, x_{j_0}^{(1)} \in V_{j_0} \bigr).
\end{align*}
Re-inserting the definitions of $F, K, M$ and comparing this against $Q_{J' \cup \{j_0\}}$, we 
conclude that $X = Q_{J' \cup \{j_0\}}$.
\end{proof}

Now we prove Theorem \ref{gvn}.  

\begin{proof}[of Theorem \ref{gvn}]
Pick any $e_0 \in H$.  It suffices to show that
$$ \E( \prod_{e \in H} f_e(\pi_e(x)) | x \in V_J ) = O( \| f_{e_0} \|_{\Box^{e_0}} ) + o_{N \to \infty}(1).$$

Applying Proposition \ref{ics} repeatedly, we see that
$$ |Q_\emptyset| \leq (1 + o_{N \to \infty}(1)) |Q_{e_0}|^{1/2^d}.$$
On the other hand, direct computation shows that
$$ Q_\emptyset = \E( \prod_{e \in H} f_e(\pi_e(x)) | x \in V_J ).$$
Thus it suffices to show that
$$ Q_{e_0} = O( \|f_{e_0} \|_{\Box^{e_0}}^{2^d} ) + o_{N \to \infty}(1).$$
We may expand
\begin{align*}
Q_{e_0} &= \E( \prod_{\omega \in \{0,1\}^{e_0}} f_{e_0}( x^{(\omega)}_{e_0} )
 \prod_{e \in H \backslash \{e_0\}} \prod_{\omega \in \{0,1\}^e: \omega_j = 0 \hbox{ for all } j \in e \backslash e_0} \nu_e( x^{(\omega)}_e ) \\
 &\quad | x^{(0)}_J, x^{(1)}_J \in V_J; x^{(0)}_{J \backslash e_0} = x^{(1)}_{J \backslash e_0} ) \\
 & = \E( W(x^{(0)}_{e_0}, x^{(1)}_{e_0}) \prod_{\omega \in \{0,1\}^{e_0}} f_{e_0}( x^{(\omega)}_{e_0} )
| x^{(0)}_{e_0}, x^{(1)}_{e_0} \in V_{e_0} ),
\end{align*}
where $W(x^{(0)}_{e_0}, x^{(1)}_{e_0} )$ is the cube counting function
$$ W(x^{(0)}_{e_0}, x^{(1)}_{e_0}) :=
\E(
\prod_{e \in H \backslash \{e_0\}} \prod_{\omega \in \{0,1\}^{e \cap e_0}} 
\nu_e( x^{(\omega)}_e ) | x^{(0)}_{J \backslash e_0} = x^{(1)}_{J \backslash e_0}
\in V_{J \backslash e_0} ).$$
On the other hand, by definition of the $\Box^{e_0}$ norm we have
$$ \E( \prod_{\omega \in \{0,1\}^{e_0}} f_{e_0}( x^{(\omega)}_{e_0} )
| x^{(0)}_{e_0}, x^{(1)}_{e_0} \in V_{e_0} )
= \| f_{e_0} \|_{\Box^{e_0}}^{2^d}.$$
Thus by the triangle inequality, it will suffice to show that
$$ \E( (W(x^{(0)}_{e_0}, x^{(1)}_{e_0})-1) \prod_{\omega \in \{0,1\}^{e_0}} f_{e_0}( x^{(\omega)}_{e_0} )
| x^{(0)}_{e_0}, x^{(1)}_{e_0} \in V_{e_0} )
= o_{N \to \infty}(1).
$$
Applying \eqref{fne} and Cauchy-Schwarz, it suffices to show that
$$ \E( |W(x^{(0)}_{e_0}, x^{(1)}_{e_0})-1|^n \prod_{\omega \in \{0,1\}^{e_0}} \nu_{e_0}( x^{(\omega)}_{e_0} )
| x^{(0)}_{e_0}, x^{(1)}_{e_0} \in V_{e_0} )
= o_{N \to \infty}(1)$$
for $n=0,2$.  Expanding this out, it suffices to show that
$$ \E( W(x^{(0)}_{e_0}, x^{(1)}_{e_0})^n \prod_{\omega \in  \{0,1\}^{e_0}} \nu_{e_0}( x^{(\omega)}_{e_0} )
| x^{(0)}_{e_0}, x^{(1)}_{e_0} \in V_{e_0} )
= 1+o_{N \to \infty}(1)$$
for $n=0,1,2$.  But the left-hand side can be rewritten as
$$ \E( 
[\prod_{e \in H \backslash \{e_0\}} \prod_{i=1}^n \prod_{\omega \in \{0,1\}^{e}: \omega_j = i \hbox{ for all } j \in e \backslash e_0} 
\nu_e( x^{(\omega)}_e )] \prod_{\omega \in  \{0,1\}^{e_0}} \nu_{e_0}( x^{(\omega)}_{e_0} )
| x^{(0)}_J, x^{(1)}_J \in V_J )
$$
and the claim thus follows from \eqref{linform-special}.
\end{proof}

\section{Dual functions and a uniform distribution property}\label{dual-sec}

We now turn to the proof of Theorem \ref{structure}.
As with \cite{gt-primes}, a key tool will be the notion of \emph{dual function} introduced in Definition \ref{dual-def}.
By definition of $\D_e$ and of the $\Box^e$ norm we observe the identity
\begin{equation}\label{fdef}
\E( f \D_e f ) = \|f\|_{\Box^e}^{2^d}
\end{equation}
for all $f: V_e \to \R$.  Thus if $f$ is not Gowers uniform in the sense that $\|f\|_{\Box^e}$ is large, then $f$ will have a large
correlation with its dual function.

The next important observation, which is a direct consequence of the dual function condition (Definition
\ref{dual-c}) is that if $f$ is bounded pointwise by $\nu_e + 1$, then the dual function is uniformly bounded:
\begin{equation}\label{def}
\D_e f(x^{(0)}_e) = O(1) \hbox{ for all } x^{(0)}_e \in V_e.
\end{equation}

We now come to a deeper property of dual functions, namely that a pseudorandom measure $\nu_e$ is uniformly distributed with respect to
arbitrary polynomial combinations of these functions.

\begin{proposition}[Uniform distribution property]\label{udp}  Let $V = (J, (V_j)_{j \in J}, d, H)$ be a hypergraph system, let $(\nu_e)_{e \in H}$ be a system of pseudorandom measures, and let $e \in H$.  Let $K$ be a finite set, and for each $k \in K$ let $f_k: V_e \to \R$ be a function such that
\begin{equation}\label{fke}
|f_k(x_e)| \leq \nu_e(x_e) + 1 \hbox{ for all } x_e \in V_e.
\end{equation}
Then we have
\begin{equation}\label{nude}
 \left|\E\left( (\nu_e(x_e) - 1) \prod_{k \in K} \D_e f_k(x_e) \big| x_e \in V_e \right)\right| = o_{N \to \infty; K}(1).
 \end{equation}
\end{proposition}

As in \cite[Lemma 6.3]{gt-primes}, the key feature here is that $K$ is allowed to be arbitrarily large.

\begin{proof}  We may use the trick of
using Lemma \ref{half} (conceding a factor of $2^{|K|}$) to replace the hypothesis \eqref{fke} by the stronger hypothesis
\begin{equation}\label{fke-better}
|f_k(x_e)| \leq \nu_e(x_e) \hbox{ for all } x_e \in V_e.
\end{equation}
Let us write $g := \nu_e - 1$.  By relabeling we may assume that $0,1 \not \in K$.  For any $e' \subseteq e$, we introduce the quantity
$Q_{e'}$, defined as
\begin{align*}
Q_{e'} := &\E\biggl( \prod_{\omega \in \{0,1\}^e: \omega_j = 0 \hbox{ for all } j \in e \backslash e'} g(x^{(\omega)}_e)\\
&\quad
\prod_{k \in K} \prod_{\omega \in (\{0,k\}^{e \backslash e'} \backslash 0^{e \backslash e'}) \times \{0,1\}^{e'}}
f_k( x^{(\omega)}_e ) \\
&\quad\quad | x^{(0)}_e \in V_e; x^{(1)}_{e'} \in V_{e'};
x^{(k)}_{e \backslash e'} \in V_{e \backslash e'} \hbox{ for all } k \in K \biggr).
\end{align*}

\begin{example}  If $e = \{1,2\}$, then
\begin{align*}
Q_\emptyset &= \E( g(x_1,x_2) \prod_{k \in K} f_k(x_1,x^{(k)}_2) f_k(x^{(k)}_1,x_2) f_k(x^{(k)}_1,x^{(k)}_2)\\
&\quad | x_1,x^{(k)}_1, \in V_1, x_2,x^{(k)}_2 \in V_2 \hbox{ for } k \in K) \\
Q_{\{1\}} &= \E( g(x_1,x_2) g(x'_1,x_2) \prod_{k \in K} f_k(x_1,x^{(k)}_2) f_k(x'_1,x^{(k)}_2) \\
&| x_1,x'_1 \in V_1, x_2,x^{(k)}_2 \in V_2 \hbox{ for } k \in K) \\
Q_{\{1,2\}} &= \E( g(x_1,x_2) g(x'_1,x_2) g(x_1,x'_2) g(x'_1,x'_2) | x_1,x'_1 \in V_1, x_2,x'_2 \in V_2)
\end{align*}
\end{example}

We claim the following analogue of Proposition \ref{ics}.

\begin{proposition}[Cauchy-Schwarz]\label{ics-2}  Let the notation and assumptions be as above.  Then for
any $e' \subsetneq e$ and $j \in e \backslash e'$, we have
$$ |Q_{e'}| \leq O_K( |Q_{e' \cup \{j\}}|^{1/2} ).$$
\end{proposition}

If we assume this proposition, then by iterating it we obtain
\begin{align*}
|\E( (\nu_e(x_e) - 1) \prod_{k \in K} \D_e f_k(x_e) | x_e \in V_e )| &= |Q_\emptyset| \\
&= O_K( |Q_e|^{1/2^d} ) \\
&= O_K( |\E( \prod_{\omega \in \{0,1\}^{e}} g(x^{(\omega)}_e) | x^{(0)}_e \in V_e; x^{(1)}_{e} \in V_{e})|^{1/2^{d}})\\
&= O_K( |\sum_{A \subseteq \{0,1\}^e} (-1)^A \E( \nu_e(x^{(\omega)}_e) | x^{(0)}_e \in V_e; x^{(1)}_{e} \in V_{e}) |^{1/2^{d}})\\
&= O_K( |\sum_{A \subseteq \{0,1\}^e} (-1)^A (1 + o_{N \to \infty}(1)) |^{1/2^{d}})\\
&= o_{N \to \infty;K}(1)
\end{align*}
as desired, where we have used \eqref{linform-special} and the binomial formula $\sum_{A \subseteq \{0,1\}^e} (-1)^A = (1-1)^{|\{0,1\}^e|} = 0$.
Thus it remains to prove the proposition.
To control $Q_{e'}$, we organize the variables $x^{(0)}_e$, $x^{(1)}_{e'}$, $x^{(k)}_{e \backslash e'}$ into three groups $\vec x, \vec y, \vec z$, where
\begin{align*}
\vec x &:= (x^{(0)}_{e \backslash \{j\}}, x^{(1)}_{e'}, (x^{(k)}_{e \backslash (e' \cup \{j\})})_{k \in K}) \in 
X := V_{e \backslash \{j\}} \times V_{e'} \times V_{e \backslash (e' \cup \{j\})}^K\\
\vec y &:= x^{(0)}_{j} \in Y := V_{j} \\
\vec z &:= (x^{(k)}_{j})_{k \in K} \in Z := V_j^K.
\end{align*}
We can then factorize 
$$ Q_{e'} = \E\left( \E( F(\vec x, \vec y) | \vec y \in Y) \E( G(\vec x, \vec z) | \vec z \in Z ) | \vec x \in X \right) $$
where
$$ F(\vec x, \vec y) := 
\prod_{\omega \in \{0\}^{e \backslash e'} \times \{0,1\}^{e'}} g(x^{(\omega)}_e)
\prod_{k \in K} \prod_{\omega \in (\{0,k\}^{e \backslash (e' \cup \{j\})} \backslash 0^{e \backslash (e' \cup \{j\})}) \times \{0\}^{\{j\}} \times \{0,1\}^{e'}}
f_k( x^{(\omega)}_e )$$
and
$$ G(\vec x, \vec z) :=
\prod_{k \in K} \prod_{\omega \in \{0,k\}^{e \backslash (e' \cup \{j\})} \times \{k\}^{\{j\}} \times \{0,1\}^{e'}}
f_k( x^{(\omega)}_e ).$$
Applying Cauchy-Schwarz, we then have
$$ |Q_{e'}| \leq \E\left( |\E( F(\vec x, \vec y) | \vec y \in Y)|^2 | \vec x \in X\right)^{1/2}
\E\left( |\E( G(\vec x, \vec z) | \vec z \in Z )|^2 | \vec x \in X \right)^{1/2}.$$
By using the definition of $F$, we have
$$ \E( |\E( F(\vec x, \vec y) | \vec y \in Y)|^2 | \vec x \in X) = Q_{e' \cup \{j\}}.$$
Thus it will suffice to show that
\begin{equation}\label{GB}
\E( |\E( G(\vec x, \vec z) | \vec z \in Z )|^2 | \vec x \in X ) = O_{K}(1).
\end{equation}
We expand the left-hand side and use \eqref{fke-better} to estimate this by
$$
\E\left( \prod_{k \in K} \prod_{\omega \in \{0,k\}^{e \backslash (e' \cup \{j\})} \times \{k,k'\}^{\{j\}} \times \{0,1\}^{e'}}
\nu_e( x^{(\omega)}_e ) | \vec x \in X, x^{(k)}_{j}, x^{(k')}_{j} \in V_{j} \hbox{ for all } k \in K\right)$$
where $k \mapsto k'$ is some arbitrary bijection from the label set $K$ to a disjoint label set $K'$ of equal cardinality.
Expanding out $\vec x$, we can factorize this expression as
$$ \E( L(x^{(0)}_{e \backslash \{j\}}, x^{(1)}_{e'})^K
| x^{(0)}_{e \backslash \{j\}} \in V_{e \backslash \{j\}};  x^{(1)}_{e'} \in V_{e'} )$$
where
$$ L(x^{(0)}_{e \backslash \{j\}}, x^{(1)}_{e'})
:= \E\left( \prod_{\omega \in \{0,k\}^{e \backslash (e' \cup \{j\})} \times \{k,k'\}^{\{j\}} \times \{0,1\}^{e'}}
\nu_e( x^{(\omega)}_e ) | x^{(k)}_{e \backslash (e' \cup \{j\}} \in V_{e \backslash (e' \cup \{j\})}; x^{(k)}_{j}, x^{(k')}_{j} \in V_{j} \right)$$
for some arbitrary label $k \in K$ (the exact value of $k$ is irrelevant).  But after relabeling, we have
\begin{align*}
L(x^{(0)}_{e \backslash \{j\}}, x^{(1)}_{e'})
&= \E\left( \prod_{\omega \in \{0,1\}^e}
\nu_e( x^{(\omega)}_e ) | x^{(1)}_{e \backslash (e' \cup \{j\}} \in V_{e \backslash (e' \cup \{j\})}; x^{(0)}_{j}, x^{(1)}_{j} \in V_{j} \right)\\
&= \E\left( M(x^{(0)}_{e \backslash \{j\}}, x^{(1)}_{e \backslash \{j\}}) | x^{(1)}_{e \backslash (e' \cup \{j\}} \in V_{e \backslash (e' \cup \{j\})} \right)
\end{align*}
where
$$ M(x^{(0)}_{e \backslash \{j\}}, x^{(1)}_{e \backslash \{j\}}) := \E( \prod_{\omega \in \{0,1\}^e}
\nu_e( x^{(\omega)}_e ) | x^{(0)}_{j}, x^{(1)}_{j} \in V_{j} ).$$
By Minkowski's inequality (i.e. the triangle inequality in $\ell^K$), we have
\begin{align*}
&\E\left( L(x^{(0)}_{e \backslash \{j\}}, x^{(1)}_{e'})^K
| x^{(0)}_{e \backslash \{j\}} \in V_{e \backslash \{j\}};  x^{(1)}_{e'} \in V_{e'} \right)^{1/K} \\
&\quad\leq \E\left( M(x^{(0)}_{e \backslash \{j\}}, x^{(1)}_{e \backslash \{j\}})^K | x^{(0)}_{e \backslash \{j\}}, x^{(1)}_{e \backslash \{j\}} \in 
V_{e \backslash \{j\}} \right)^{1/K},
\end{align*}
and hence by the correlation condition \eqref{correlation} we obtain \eqref{GB} as required.
\end{proof}

An immediate corollary of Proposition \ref{udp} and the triangle inequality is

\begin{corollary}[Uniform distribution property with respect to polynomials]\label{udp-poly}  Let $V = (J, (V_j)_{j \in J}, d, H)$ be a hypergraph system, let $(\nu_e)_{e \in H}$ be a system of pseudorandom measures, and let $e \in H$.  Let $K$ be a finite set, let $D \geq 0$ be an integer, and let
$P: \R^K \to \R$ be a polynomial of degree $D$ in $K$ variables, with all coefficients bounded by some quantity $M$.  For
 each $k \in K$ let $f_k: V_e \to \R$ be a function such that
\begin{equation}\label{fke-poly}
|f_k(x)| \leq \nu_e(x) + 1 \hbox{ for all } x \in V_e.
\end{equation}
Then we have
\begin{equation}\label{nude-poly}
 |\E\left( (\nu_e(x_e) - 1) P( (\D_e f_k(x_e))_{k \in K} ) \bigl| x_e \in V_e \right)| = o_{N \to \infty; K, D, M}(1).
 \end{equation}
 (Recall we allow our constants to depend implicitly on $J$).
\end{corollary}

\begin{remark} Following \cite{gt-primes}, one could also extend this corollary from polynomials to continuous functions using the Weierstrass approximation theorem (and the bound \eqref{def}), but we will not need to do so here.
\end{remark}

\section{$\sigma$-algebras of dual functions}\label{sigma-sec}

We continue the proof of Theorem \ref{structure}.
As in \cite[Theorem 8.1]{gt-primes}, we will exploit the above uniform distribution property to associate a $\sigma$-algebra to every
dual function.  We first give a minor variant of \cite[Proposition 7.2]{gt-primes}:

\begin{proposition}[Each bounded function generates a $\sigma$-algebra]\label{dual-sigma}  Let $V = (J, (V_j)_{j \in J}, d, H)$ be a hypergraph system, let $(\nu_e)_{e \in H}$ be a system of pseudorandom measures, and let $e \in H$.  Let $0 < \eps < 1$ and $0 < \sigma < 1/2$ be parameters, let
$I$ be an interval in $\R$, and let $G: V_e \to I$ be a function.  Then, if the pseudorandomness parameter $N$ is sufficiently large depending on $\eps,\sigma$, there exists a $\sigma$-algebra $\B_{\eps,\sigma,e}(G)$ on $V_e$
with the following properties:
\begin{itemize}
\item ($G$ lies in its own $\sigma$-algebra)  For any $\sigma$-algebra $\B$ on $V_e$, we have
\begin{equation}\label{G-approx}
|G(x) - \E\left( G | \B_{\eps,\sigma,e}(G) \vee \B \right)(x)| \leq \eps \hbox{ for all } x \in V_e.
\end{equation}
\item (Bounded complexity) $\B_{\eps,\sigma,e}(G)$ is generated by at most $O_{\eps,I}(1)$ atoms.
\item (Approximation by polynomials of $G$) If $A$ is any atom in $\B_{\eps,\sigma,e}(G)$, then there exists a polynomial $P_{A,\eps,\sigma,I}$ of degree
$O_{\eps,\sigma,I}(1)$ and all co-efficients $O_{\eps,\sigma,I}(1)$, such that $P_{A,\eps,\sigma,I}(x) = O(1)$ for all $x \in I$ and
\begin{equation}\label{uapprox}
\E\left( |1_A(x) - P_{A,\eps,\sigma,I}(G(x))| (\nu_e(x) + 1) \bigl| x \in V_e \right) = O(\sigma).
\end{equation}
\end{itemize}
\end{proposition}

\begin{proof}
Observe from Fubini's theorem that
$$ 
\int_0^1
\sum_{n \in \Z} \E\big( {\bf 1}_{G(x) \in [\eps(n-\sigma+\alpha), \eps(n+\sigma+\alpha)]} (\nu_e(x)+1) \; \big| \; x \in V_e \big)\ d\alpha
= 2\sigma \E( \nu_e(x)+1 | x \in V_e ).$$
Since $\nu_e$ is pseudorandom, we have 
\begin{equation}\label{avy}
\E( \nu_e(x)+1 | x \in V_e ) = O(1)
\end{equation}
if $N$ is large enough.  Thus by the pigeonhole principle
we can find $0 \leq \alpha \leq 1$ such that
\begin{equation}\label{summy}
\sum_{n \in \Z} \E\big( {\bf 1}_{G(x) \in [\eps(n-\sigma+\alpha), \eps(n+\sigma+\alpha)]} (\nu_e(x)+1) \; \big| \; x \in V_e \big) = O(\sigma).
\end{equation}
We now set $\B_{\eps,\sigma,e}(G)$ to be the $\sigma$-algebra whose atoms are the sets $G^{-1}([\eps (n+\alpha), \eps(n+1+\alpha)))$ for $n \in \Z + \alpha$ (discarding all the empty atoms, of course).  This is well-defined since the intervals $[\eps (n+\alpha), \eps(n+1+\alpha))$ tile the real line.   Since $G$ takes values in $I$ we see that there are only $O_{\eps,I}$ non-empty atoms.

It is clear that if $\B$ is an arbitrary $\sigma$-algebra on $V_e$, then on any atom of $\B \vee \B_\eps(G)$, the function $G$ takes values
in an interval of diameter $\eps$, which yields \eqref{G-approx}.  
Now we verify the approximation by continuous functions property.  Let $A := G^{-1}([\eps (n+\alpha), \eps(n+1+\alpha)))$ be an atom.
Since $G$ takes values in $I$, we may assume that $n = O_{I,\eps}(1)$, since $A$ is empty otherwise; note that this already establishes
the bounded complexity property.  By combining Urysohn's lemma with the Weierstrass approximation theorem, we can 
find a polynomial $P_{A,\eps,\sigma,I}$ which is equals $1 + O(\sigma)$ on $[\eps(n+\alpha+\sigma), \eps(n+\alpha+1-\sigma)]$, equals
$O(\sigma)$ on $I \backslash [\eps(n+\alpha-\sigma), \eps(n+\alpha+1+\sigma)]$, and equals $O(1)$ on all of $I$.
Furthermore, a simple compactness argument shows that the degree of $P_{A,\eps,\sigma,I}$ can be chosen to be $O_{\eps,\sigma,I}(1)$,
and all the coefficients can also be chosen to be $O_{\eps,\sigma,I}(1)$.  We have the pointwise estimate
$$ |1_A(x) - P_{A,\eps,\sigma,I}(G(x))| = O(\sigma) + \sum_{m=n}^{n+1} O\left( 1_{[\eps(m+\alpha-\sigma), \eps(m+\alpha+\sigma)]}(G(x))\right)$$
so by applying \eqref{summy} and \eqref{avy} we obtain \eqref{uapprox}.
\end{proof}

We specialize this Proposition to functions $G$ which are dual functions, to conclude the following analogue of \cite[Proposition 7.3]{gt-primes}.

\begin{proposition}\label{bohr-sigma} Let $V = (J, (V_j)_{j \in J}, d, H)$ be a hypergraph system, let $(\nu_e)_{e \in H}$ be a system of pseudorandom measures, and let $e \in H$.  Let $K$ be an integer, and for each 
$1 \leq k \leq K$ let $f_k: V_e \to \R$ be a function such that \eqref{fke-poly} holds.
 Let $0 < \eps < 1$ and $0 < \sigma < 1/2$ be parameters, and let $\B_{\eps,\sigma,e}(\D_e f_k)$ for $1 \leq k \leq K$
 be constructed as in Proposition \ref{dual-sigma} (note from \eqref{def} that we can take $I$ to be a fixed interval of width $O(1)$)).  Let $\B := \bigvee_{1 \leq k \leq K} \B_{\eps,\sigma,e}(\D_e f_k)$.
 Then if $\sigma$ is sufficiently small depending on $K, \eps$, and $N$ is sufficiently large depending on $K$, $\eps$, $\sigma$, $J$, $d$, we have
\begin{equation}\label{trivial-specific}
| \D_e f_k(x) - \E( \D_e f_k | \B)(x) | \leq \eps \hbox{ for all } 1 \leq k \leq K, x \in V_e.
\end{equation}
Furthermore there exists a set $\Omega \in \B$ obeying the smallness condition
\begin{equation}\label{mok-small} 
\E( (\nu_e(x) + 1) 1_{\Omega}(x) | x \in V_e ) = O_{K,\eps}(\sigma^{1/2})
\end{equation}
and such that 
\begin{equation}\label{unif-mu}
\E(\nu_e - 1 | \B)(x) = O_{K,\eps}(\sigma^{1/2}) \hbox{ for all } x \in V_e \backslash \Omega.
\end{equation}
\end{proposition}

\begin{proof}
The claim \eqref{trivial-specific} follows immediately from \eqref{G-approx}.  Now we prove \eqref{mok-small} and \eqref{unif-mu}.
Since each of the $\B_{\eps,\sigma,e}(\D_e f_k)$ are generated by $O_\eps(1)$ atoms, we see that $\B$ is generated by $O_{K,\eps}(1)$ atoms.
Call an atom $A$ of $\B$ \emph{small} if $\E( (\nu_e(x) + 1) 1_A(x) | x \in V_e) \leq \sigma^{1/2}$, and let $\Omega$ be the union of all the small atoms.  Then clearly $\Omega$ lies in $\B$ and obeys \eqref{mok-small}.  To prove the remaining claim \eqref{unif-mu}, it suffices to show that
\begin{equation}\label{sigma-smooth}
\frac{ \E( (\nu_e(x) - 1) 1_A(x) | x \in V_e ) }{ \E( 1_A(x) | x \in V_e ) } =o_{N \to \infty; K,\eps,\sigma}(1) + O_{K,\eps}(\sigma^{1/2})
\end{equation}
for all atoms $A$ in $B$ which are not small.  However, by definition of ``small'' we have
$$ \E( (\nu_e(x) - 1) 1_A(x) | x \in V_e ) + 2 \E( 1_A(x) | x \in V_e ) = \E( (\nu_e(x)+1) 1_A(x) | x \in V_e ) \geq \sigma^{1/2}.$$
Thus to complete the proof of \eqref{sigma-smooth} it will suffice (since $\sigma$ is small and $N$ is large) to show that
\begin{equation}\label{nubf}
\E( (\nu_e(x) - 1) 1_A(x) | x \in V_e ) = o_{N \to \infty; K,\eps,\sigma}(1) + O_{K,\eps}(\sigma).
\end{equation}
On the other hand, since $A$ is the intersection of atoms $A_k \in \B_{\eps,\sigma,e}(\D_e f_k)$ for each
$1 \leq k \leq K$, we see from Proposition \ref{dual-sigma} (and H\"older's inequality) that we can find a polynomial
$P: \R^K \to \R$ of degree $O_{\eps,\sigma,K}(1)$ and coefficients $O_{\eps,\sigma,K}(1)$ such that
$$\E\left( (\nu_e(x) + 1) |1_A(x) - P(\D_e f_1(x), \ldots, \D_e f_k(x))| \biggl| x \in V_e \right) = O_K(\sigma),$$
so in particular
$$\E\left( (\nu_e(x) - 1) (1_A(x) - P(\D_e f_1(x), \ldots, \D_e f_k(x)) \biggl| x \in V_e \right) = O_K(\sigma).$$
On the other hand, Corollary \ref{udp-poly} we have
$$\E\left( (\nu_e(x) - 1) P(\D_e f_1(x), \ldots, \D_e f_k(x)) \biggl| x \in V_e \right) = o_{N \to \infty; K,\eps,\sigma}(1).$$
The claim \eqref{nubf} now follows from the triangle inequality.
\end{proof}

\section{A Furstenberg tower, and the proof of Theorem \ref{structure}}\label{structure-sec}

We are now ready to prove Theorem \ref{structure}.  As in \cite{gt-primes}, this theorem shall be proven by a constructing a 
Furstenberg tower of increasingly complex $\sigma$-algebras.

Fix $V$, $e$, $\nu_e$, $f_e$, $\eps$.
We shall need a parameter $0 < \sigma \ll \eps$ which we shall choose 
later, and then we shall assume $N$ is sufficiently large depending on $\sigma$ and $\eps$.

To construct $\B_e$ and $\Omega_e$ we shall iteratively construct a sequence of basic Gowers anti-uniform
functions $\D_e F_{e,1}, \ldots, \D_e F_{e,K}$ on $V_e$, 
exceptional sets $\Omega_{e,0} \subseteq\Omega_{e,1} \subseteq\ldots \subseteq\Omega_{e,K} \subseteq V_e$, and a nested sequence of $\sigma$-algebras 
$\B_{e,0} \subseteq\ldots \subseteq\B_{e,K}$ for some integer $K \geq 0$ as follows.

\begin{itemize}
\item Step 0.  Initialize $K = 0$, and define $\B_{e,0} := \{ \emptyset, V_e\}$ and $\Omega_{e,0} := \emptyset$.  
\item Step 1.  Set $F_{e,K+1} := (1 - 1_{\Omega_{e,K}}) (f_e - \E(f_e|\B_{e,K}))$.
If we have 
$$ \| F_{e,K+1} \|_{\Box^e} \leq \eps^{1/2^{d+1}}$$
then we set $\Omega_e := \Omega_{e,K}$ and $\B_e = \B_{e,K}$, and successfully terminate the algorithm.
\item Step 2.  If instead we have
\begin{equation}\label{non-terminate}
 \| F_{e,K+1} \|_{\Box^e} > \eps^{1/2^{d+1}},
\end{equation}
then we let $\B_{e,K+1} := \B_{e,K} \vee \B_{\eps,\sigma,e}(\D_e F_{e,K+1})$, where $\B_{\eps,\sigma,e}(\D_e F_{e,K+1})$ is as in Proposition \ref{dual-sigma}. 
\item Step 3.  Locate an exceptional set $\Omega_{e,K+1} \supset \Omega_{e,K}$ in $\B_{e,K+1}$ obeying
the smallness condition
\begin{equation}\label{om-small} \E( (\nu_e(x_e) + 1) 1_{\Omega_{e,K+1}}(x_e) | x_e \in V_e) = O_{K,\eps}(\sigma^{1/2})
\end{equation}
and such that we have the bound
\begin{equation}\label{nube}
\E( \nu_e | \B_{e,K+1} )(x_e) = 1 + O_{K,\eps}(\sigma^{1/2}) \hbox{ for all } x_e \in V_e \backslash \Omega_{e,K+1}.
\end{equation}
If such an exceptional set cannot be found, we terminate the algorithm with an error; otherwise, we move
on to Step 4.
\item Step 4.  Increment $K$ to $K+1$, and return to Step 1.
\end{itemize}

Let $K_0$ be a large multiple of $1/\eps$ to be chosen later.
We claim that this algorithm necessarily terminates without error in Step 1 in less than $K_0$
steps (so $K$ always remains smaller than $K_0$), if $N$ is sufficiently large depending on $\eps$ and $\sigma$.  Assuming
this for the moment, then by construction we have \eqref{nu-cube}, as well as
the bounds
$$
\E( (\nu_e(x_e) + 1) 1_{\Omega_{e}}(x_e) | x_e \in V_e) = O_{\eps}(\sigma^{1/2})
$$
and 
$$ \E( \nu_e | \B_{e} )(x_e) - 1 = O_{\eps}(\sigma^{1/2}) \hbox{ for all } x_e \in V_e \backslash \Omega_{e},$$
where we use the hypothesis that $K = O(K_0)= O_\eps(1)$.  If we choose $\sigma$ sufficiently small depending on $\eps$,
and then assume $N$ sufficiently large depending on $\eps$ and $\sigma$, we thus see that the right-hand sides of these bounds
can be made as small as desired, thus obtaining \eqref{omegasmall} and \eqref{nu-distributed}.

It remains to show that the algorithm does indeed terminate without error in less than $K_0$ steps.  We first show that it will not
terminate with error in the first $K_0$ steps.  To see this, observe that we only have to show that Step 3 can be executed without error
whenever $K \leq K_0$.  But observe from \eqref{fnue} and \eqref{nube} for step $K-1$ (if $K \geq 1$) or from the bound \eqref{nue-mean} (if $K=0$) that we have the pointwise bound
\begin{equation}\label{feek}
 |\E(f_e|\B_{e,K})(x_e)| \leq 1 + O_{K,\eps}(\sigma^{1/2}) + o_{N \to \infty}(1) \hbox{ for all } x_e \not \in \Omega_{e,K}
 \end{equation}
and hence by \eqref{fnue} again
\begin{equation}\label{Fek-bound}
|F_{e,K+1}(x_e)| \leq \nu_e(x_e) + 1 + O_{K,\eps}(\sigma^{1/2}) + o_{N \to \infty}(1) \hbox{ for all } x_e \in V_e.
\end{equation}
Applying (a slightly rescaled) version of \eqref{def}, we conclude
\begin{equation}\label{def-2}
|\D_e F_{e,K+1}(x_e)| \leq O(1) + O_{K,\eps}(\sigma^{1/2}) + o_{N \to \infty}(1)  \hbox{ for all } x_e \in V_e.
\end{equation}
The claim now follows by letting $\Omega \in \B_{e,K+1}$ be the set defined in Proposition \ref{bohr-sigma}, using
the family of functions $F_{e,1},\ldots,F_{e,K+1}$ instead of $f_1,\ldots,f_K$ and then setting $\Omega_{e,K+1} := \Omega_{e,K} \cup \Omega$.

The only other remaining possibility to eliminate is that the first $K_0$ loops of the algorithm are executed without error or termination.
We shall show this cannot happen by establishing the energy incrementation inequality
\begin{equation}\label{eje-increment}
\begin{split}
&\E( (1 - \Omega_{e,j}(x_e)) \E( f_e | \B_{e,j} )(x_e)^2 | x_e \in V_e )\\
&\quad\geq \E( (1 - \Omega_{e,j-1}(x_e)) \E( f_e | \B_{e,j-1} )(x_e)^2 | x_e \in V_e ) + c^2 \eps
\end{split}
\end{equation}
for all $1 \leq j \leq K_0$, and some $c > 0$ independent of $\eps$ or $\sigma$.  On the other hand, the quantity $\E( (1 - \Omega_{e,j}(x_e)) \E( f_e | \B_{e,j} )(x_e)^2 | x_e \in V_e )$
is clearly bounded below by zero, and bounded above by
$$ \E( (1 - \Omega_{e,j}(x_e)) \E( \nu_e | \B_{e,j} )(x_e)^2 | x_e \in V_e ) \leq 1 + O_{j,\eps}(\sigma^{1/2})$$
thanks to \eqref{nube}.  The two facts are contradictory by choosing $K_0$ to be a large multiple of $1/\eps$, if $\sigma$ is chosen sufficiently small.

It remains to prove \eqref{eje-increment}.  Since the algorithm successfuly executed the first $K_0$ loops, we have
$$  \| F_{e,j} \|_{\Box^e} > \eps^{1/2^{d+1}}.$$
Raising this to the power $2^{d-1}$, and using \eqref{fdef}, we conclude
$$ \langle F_{e,j}, {\cal D}_e F_{e,j} \rangle > \eps^{1/2}$$
where we are using the usual inner product
$$ \langle f, g \rangle := \E( f(x_e) g(x_e) | x_e \in V_e ).$$
On the other hand, from \eqref{Fek-bound}, \eqref{def-2}, \eqref{om-small} we have
$$ \langle 1_{\Omega_{e,j}} F_{e,j}, {\cal D}_e F_{e,j} \rangle = O_{K,\eps}(\sigma^{1/2}) [ O(1) + O_{j,\eps}(\sigma^{1/2}) + o_{N \to \infty}(1) ]
= O_{\eps}(\sigma^{1/2})$$
(since $j = O(K_0) = O_\eps(1)$)
and hence by the triangle inequality
$$ \langle (1 - 1_{\Omega_{e,j}}) F_{e,j}, {\cal D}_e F_{e,j} \rangle > \eps^{1/2} - O_{\eps}(\sigma^{1/2})$$
On the other hand, from \eqref{trivial-specific} we have
$$ {\cal D}_e F_{e,j}(x_e) = \E( {\cal D}_e F_{e,j} | \B_{e,j} )(x_e) + O(\eps) \hbox{ for all } x_e \in V_e$$
and hence by \eqref{Fek-bound} and \eqref{nue-mean}
$$ \langle (1 - 1_{\Omega_{e,j}}) F_{e,j}, {\cal D}_e F_{e,j} - \E( {\cal D}_e F_{e,j} | \B_{e,j} )  \rangle = O(\eps).$$
We conclude that
$$ \langle (1 - 1_{\Omega_{e,j}}) F_{e,j}, \E( {\cal D}_e F_{e,j} | \B_{e,j}) \rangle > \eps^{1/2} - O_{\eps}(\sigma^{1/2})$$
Since $1_{\Omega_{e,j}}	$ is measurable in $\B_{e,j}$, we obtain
$$ \langle (1 - 1_{\Omega_{e,j}}) \E( F_{e,j} | \B_{e,j} ), \E( {\cal D}_e F_{e,j} | \B_{e,j}) \rangle > \eps^{1/2} - O_{\eps}(\sigma^{1/2}).$$
By \eqref{def} and Cauchy-Schwarz we conclude
$$ \| (1 - 1_{\Omega_{e,j}}) \E( F_{e,j} | \B_{e,j} ) \|_{L^2} > c \eps^{1/2} -  O_{\eps}(\sigma^{1/2})$$
for some $c > 0$ independent of $\eps,\sigma$, and where 
$$\|f\|_{L^2} = \langle f, f\rangle^{1/2} = \E(f^2(x_e)|x_e \in V_e)^{1/2}.$$
Using the definition of $F_{e,j}$, we conclude
$$ \| (1 - 1_{\Omega_{e,j}}) (\E(f_e|\B_{e,j-1}) - \E( f_e | \B_{e,j} )) \|_{L^2} \geq c \eps^{1/2} - O_{\eps}(\sigma^{1/2}).$$
We now use the cosine rule to conclude
\begin{align*}
&\| (1 - 1_{\Omega_{e,j}}) \E( f_e | \B_{e,j} ) \|_{L^2}^2 \geq
\| (1 - 1_{\Omega_{e,j}}) \E( f_e | \B_{e,j-1} ) \|_{L^2}^2 + c^2 \eps - O_{\eps}(\sigma^{1/2}) \\
&\quad + 2 \left\langle 
(1 - 1_{\Omega_{e,j}}) 
[\E( f_e | \B_{e,j} ) - \E(f_e,\B_{e,j-1})], 
(1 - 1_{\Omega_{e,j}}) \E( f_e | \B_{e,j-1} ) 
\right\rangle.
\end{align*}
The inner product here can be rewritten as 
$$
\left\langle \E( f_e | \B_{e,j} ) - \E(f_e|\B_{e,j-1}), (1 - 1_{\Omega_{e,j}}) \E( f_e | \B_{e,j-1} ) \right\rangle.$$
Now observe that the quantity in square brackets has zero conditional expectation with respect to $\B_{e,j-1}$, and
in particular is orthogonal to $(1 - 1_{\Omega_{e,j-1}}) \E( f_e | \B_{e,j-1} )$.  Thus the above inner product can be rewritten as
$$
\left\langle \E( f_e | \B_{e,j} ) - \E(f_e|\B_{e,j-1}), (1_{\Omega_{e,j-1}} - 1_{\Omega_{e,j}}) \E( f_e | \B_{e,j-1} ) \right\rangle.$$
Observe that the second factor is measurable in $\B_j$, and so the inner product can be rewritten again as
$$
\left\langle  f_e - \E(f_e|\B_{e,j-1}), (1_{\Omega_{e,j-1}} - 1_{\Omega_{e,j}}) \E( f_e | \B_{e,j-1} ) \right\rangle.$$
Using \eqref{nube} we have $\E(f_e | \B_{e,j-1} ) = O(1)$ outside of $\Omega_{e,j-1}$, and so from this and \eqref{om-small}, \eqref{fnue} we see
that this inner product is $O_{j,\eps}(\sigma^{1/2}) = O_\eps(\sigma^{1/2})$, since $j = O(K_0) = O_\eps(1)$.  Summarizing all the above
computations, we conclude that
$$ \| (1 - 1_{\Omega_{e,j}}) \E( f_e | \B_{e,j} ) \|_{L^2}^2 \geq
\| (1 - 1_{\Omega_{e,j}}) \E( f_e | \B_{e,j-1} ) \|_{L^2}^2 + c^2 \eps - O_\eps(\sigma^{1/2}).$$
On the other hand, from \eqref{feek}, \eqref{om-small} we have
$$ \| (1 - 1_{\Omega_{e,j}}) \E( f_e | \B_{e,j-1} ) \|_{L^2}^2 = \| (1 - 1_{\Omega_{e,j-1}}) \E( f_e | \B_{e,j-1} ) \|_{L^2}^2
+ O_{\eps}(\sigma^{1/2})$$
and we thus conclude \eqref{eje-increment}, if $\sigma$ is sufficiently small depending on $\eps$.
This concludes the proof of Theorem \ref{structure}.
\endprf

\section{Constellations in the Gaussian primes: preliminaries}\label{Gaussian-sec}

We now begin the proof of Theorem \ref{main}.  In this section we shall reduce matters (via a number of somewhat artificial technical reductions) to the point where we can apply Theorem \ref{sz-multi-alt-relative}, at which point the only remaining task will be to establish that a certain family of measures $(\nu_e)_{e \in H}$ constructed here is pseudorandom.  This will then be achieved in the next section.  

By making the substitution $a \mapsto a - rv_0$ if necessary we may take $v_0 = 0$.  By adding some dummy elements to the $v_j$ if necessary, we may assume the ergodic hypothesis that 
the $v_j$ (and hence their differences $v_i - v_j$) generate $\Z[i]$ as an additive group.  Such a maneuvre is terrible for the 
quantitative bounds, but for the qualitative question of merely establishing infinitely many prime constellations, it is harmless.
In fact by adding a few more dummy elements we can easily impose the following slightly stronger hypothesis:

\begin{hypothesis}[Improved ergodic hypothesis]\label{super-ergodic}  
If $i, j$ are two distinct elements of $J$, then the vectors $\{ v_k - v_j: k \in J \backslash \{i,j\} \}$ span $\Z[i]$ as an
additive group.
\end{hypothesis}

\begin{remark} The above hypothesis is not strictly necessary for our argument, but it does simplify matters slightly and is easy to attain, so we shall take advantage of it.
\end{remark}

Henceforth we allow all implicit constants in the $O()$ and $o()$ notation to depend on $k$ and $v_0,\ldots,v_{k-1}$.  We will also
use $C, c > 0$ to denote various positive constants (possibly depending on the above parameters) which can vary from line to line.

To avoid confusion let us use the terminology \emph{rational prime} to denote a prime in the natural
numbers $\Z^+$, and \emph{Gaussian prime} to denote a prime in $\Z[i]$.  Thus for instance $5 = (2+i)(2-i)$ is a rational prime but not a Gaussian prime.  Similarly we use \emph{rational integer} to denote an element of $\Z$.
We let $\Z[i]^\times := \{1, i, -1, -i\}$ denote the \emph{Gaussian units}, that is the invertible elements in $\Z[i]$.
Let us call two Gaussian non-zero integers \emph{associate} if their quotient is a Gaussian unit, and \emph{non-associate} otherwise.  

Given any non-zero Gaussian integer $z$, we define its \emph{norm} $\N(z)$ to be the quantity $\N(z) := |\Z[i] / z\Z[i]|$; it is easy to verify that $\N(zz') = \N(z) \N(z')$ and $\N(a+bi) = a^2+b^2$.  
As is well known (see e.g. \cite{hardy-wright}), $\N(P[i])$ consists of the number 2, as well as the rational primes equal to 1 modulo 4, and
the squares of the rational primes equal to 3 modulo 4.  Of these three cases, the second case is by far the most prevalent.  As the other two cases cause some minor difficulty\footnote{Specifically, the exceptional primes cause an unwelcome irregularity, namely that the Gaussian primes have an anomalous density on certain lines, such as the real or imaginary axes, which will disrupt the pseudorandomness hypothesis we impose later.  This is ultimately due to the fact that the field $\Z[i]/p\Z[i]$ does not have rational prime order if $p$ is exceptional.},
 we shall remove them by defining the \emph{unexceptional Gaussian primes} $P[i]'$ to be those Gaussian primes $p \in Z[i]$ such that $\N(p)$ is a rational prime equal to 1 modulo 4, and define $\Z[i]'_{sq}$ to be those non-zero square-free Gaussian integers whose prime factorization consists only of unexceptional Gaussian primes (and Gaussian units, of course).  Note that unexceptional Gaussian primes have non-zero real part and non-zero imaginary part.  We define $P[i]'_+ \subset P[i]'$ to be those unexceptional Gaussian primes which lie in the first quadrant; thus every unexceptional Gaussian prime is conjugate to exactly one prime in $P[i]'_+$.  

Clearly, in order to obtain infinitely many constellations in $P[i]$ it suffices to obtain infinitely many constellations in $P[i]'$.
The first main task is to obtain a number-theoretic pseudorandom majorant $\nu$ for the unexceptional Gaussian primes, or more precisely for
a weight function adapted to a variant of the unexceptional Gaussian primes in which all the non-uniformity arising from small divisors has been
eliminated.  
Recall that every rational prime $p$ equal to 1 modulo 4 is the norm of exactly eight unexceptional Gaussian primes (two of which lie in $P[i]'_+$).
From Dirichlet's theorem (in the modulo 4 case) and the prime number theorem we thus have
\begin{equation}\label{primes-square}
 |\{ p \in P[i]': \N(p) \leq N^2 \}| = (2 + o_{N \to \infty}(1)) \frac{N^2}{\log N}.
 \end{equation}

\begin{remark}
This bound is of course consistent with (a very simple case of) the Chebotarev density theorem.
\end{remark}

We now adopt a Gaussian integer version of the ``$W$-trick'' from \cite{gt-primes}, whose purpose is to eliminate non-uniformities in the Gaussian primes which arise from small divisors.  Let $N$ be a large rational prime; we view this as a parameter which will eventually be
sernt to infinity.  Let $w = w(N)$ be a positive rational integer which grow very slowly to infinity as $N \to \infty$, thus
we can write any expression of the form $o_{w \to \infty}(1)$ as $o_{N \to \infty}(1)$, and any specified expression of the form $o_{N \to \infty; w}(1)$ as $o_{N \to \infty}(1)$; we shall frequently  take advantage of these facts in the sequel without further comment.  
We let $W = W(N) := \prod_{p \in P[i]: \N(p) \leq w} \N(p)$ be the product of the norms of all the Gaussian primes of 
norm less than $w$; note that the growth comments about $w$ apply just as well to $W$, thus for instance any specified expression of the form $o_{W \to \infty}(1)$ or $o_{N \to \infty; W}(1)$ can be written as $o_{N \to \infty}(1)$.  
We can partition $\Z[i]$ into cosets $W \cdot \Z[i] + b$, where $b \in [0,W)^2$.  Let $\phi_{\Z[i]}(W)$ denote the 
number of Gaussian integers in $[0,W)^2$ which are coprime to $W$.  We also need a small number $0 < \epsilon < \frac{1}{100}$, depending
on $k, v_1,\ldots,v_k$, to be chosen later.  By \eqref{primes-square} we have
$$ 
|\{ p \in P[i]': \frac{(\epsilon N W)^2}{2} \leq \N(p) \leq (\epsilon N W)^2 \}| \geq c_\epsilon \frac{N^2 W^2}{\log N}$$
for some $c_\epsilon > 0$, if $N$ is sufficiently large depending on $\epsilon$ (and $W$ is slowly growing with respect to $N$).
We caution that the value of $c_\epsilon$ will vary from line to line.
By the pigeonhole principle\footnote{One could also use the Gaussian integer analogue of Dirichlet's theorem at this point, but it is unnecessary for this argument.  One could also replace $P[i]'$ here by any dense subset of $P[i]'$ without difficulty; since $P[i]'$ has density one inside $P[i]$ this ultimately means that we can generalise Theorem \ref{main} to dense subsets of $P[i]$.  We omit the details.} we can find
a $b \in [0,W)^2$ coprime to $W$ such that
\begin{equation}\label{acapback}
|A_b| \geq c_\epsilon \frac{N^2 W^2}{\phi_{\Z[i]}(W) \log N}
\end{equation}
for some slightly different $c_\epsilon > 0$, where $A_b \subset \Z[i]$ is the set
$$ A_b := \{ n \in \Z[i]: \frac{1}{2} \epsilon^2 N^2 \leq \N(p) \leq \epsilon^2 N^2; Wn + b \in P[i]' \}.$$
Fix such a $b$.
It would now suffice to the quantitative estimate
\begin{align*}
\{ (a,r) \in \Z[i] \times &\Z: a + r v_j \in A_b \hbox{ for } 0 \leq j < k \}| \\
&\geq 
(c_\epsilon - o_{N \to \infty;\epsilon}(1)) \frac{N^3 W^{2k}}{\phi_{\Z[i]}(W)^k \log^k N}.
\end{align*}
Note that the contribution of the degenerate cases $r=0$ becomes negligible for $N$ large enough.
  
Let $Z := \Z/N\Z$, and let $\pi: \Z[i] \to Z^2$ be the obvious projection map.
If $\epsilon$ is sufficiently small depending on $v_1,\ldots,v_k$, and $N$ is large enough
depending on $\epsilon, v_1,\ldots,v_k$, it now suffices to show that
\begin{align*}
&\E\left( \prod_{0 \leq j < k} \frac{\phi_{\Z[i]}(W) \log N}{W^2}  1_{\pi(A_b)}(a + r v_j) \bigg| a \in Z^2; r \in Z \right)\\
&\quad\geq c_\epsilon - o_{N \to \infty; \epsilon}(1)
\end{align*}
for some $c_\epsilon > 0$.

\begin{remark} This bound is consistent with the analogue of the Hardy-Littlewood prime tuples conjecture for Gaussian primes; see \cite{zoo}.  Of course, our work here does not make any serious progress towards that conjecture.
\end{remark}

From \eqref{acapback} we have
\begin{equation}\label{acap-2}
 \E\left( \frac{\phi_{\Z[i]}(W) \log N}{W^2}  1_{\pi(A_b)}(x) \bigg| x \in Z^2 \right) \geq c_\epsilon > 0.
 \end{equation}

The next step is to construct a suitable pseudorandom measure on $Z^2$ so that we may invoke Theorem \ref{sz-multi-alt-relative}. 
One could modify the truncated divisor sums of Goldston and Y{\i}ld{\i}r{\i}m (as used in \cite{gt-primes} for the rational primes)
directly.  However we take advantage of a slight simplification to their approach introduced in \cite{tao-gy} which uses less information on the Gaussian integer $\zeta$-function (in particular, using only the very crude zero-free region in the vicinity of the pole at $s=1$) to obtain a qualitatively similar result in a slightly more elementary fashion.

Define the \emph{M\"obius function} $\mu_{\Z[i]}: \Z[i] \to \R$ for the Gaussian integers by setting $\mu_{\Z[i]}(n) := (-1)^m$ when $n \in \Z[i]$ is the product of 
$m$ pairwise non-associate Gaussian primes, and zero otherwise.  Similarly, we define the \emph{von Mangoldt function} $\Lambda_{\Z[i]}: \Z[i] \to \R^+$ for the Gaussian integers by setting $\Lambda_{\Z[i]}(n) := \log \N(p)$ if $n$ is associate to a power of a Gaussian prime $p$, and equal to zero otherwise.  From unique factorization in $\Z[i]$, one easily verifies the identities
\begin{align*}
\log\N(n) &= \frac{1}{4} \sum_{d \in \Z[i] \backslash \{0\}: d|n} \Lambda_{\Z[i]}(d)\\ 
\Lambda_{\Z[i]}(n) &= \frac{1}{4} \sum_{d \in \Z[i] \backslash \{0\}: d|n} \mu_{\Z[i]}(d) \log \N(\frac{n}{d}) 
\end{align*}
for all $n \in \Z[i] \backslash \{0\}$; the factor of $\frac{1}{4}$ is due to the four Gaussian units $1, -1, i, -i$.  

We now smoothly truncate the above formula for $\Lambda_{\Z[i]}(n)$ to obtain a truncated divisor sum of Goldston-Y{\i}ld{\i}r{\i}m type, and also restrict to the unexceptional Gaussian integers $\Z[i]'_{sq}$.
Let $R := N^{c}$ for some small $c = c_k > 0$ to be chosen later (e.g. $c_k = 2^{-100k}$ would suffice).
Let $\varphi: \R \to \R^+$ be a smooth bump function\footnote{Any standard bump function will do here.  The actual
truncated divisor sum corresponding to Goldston-Y{\i}ld{\i}r{\i}m corresponds to the choice $\varphi(x) := \max(1-|x|,0)$, which offers the advantage that all integrals involving $\varphi$ can (in principle) be worked out explicitly, but has only a limited amount of regularity which necessitates knowledge of the zero-free region on the axis $\Re s = 1$ in order to proceed.}
 supported on $[-1,1]$ which equals 1 at 0 (any standard bump function would do here), and define
\begin{equation}\label{trunc-def}
\Lambda_{\Z[i]'_{sq},R,\varphi}(n) := \frac{1}{4} \log \N(R) \sum_{d \in \Z[i]'_{sq}: d|n} \mu_{\Z[i]}(d) \varphi\left( \frac{\log \N(d)}{\log \N(R)} \right),
\end{equation}
One observes that $\Lambda_{\Z[i],R,\varphi}(n) = \log \N(R)$ whenever $n$ is an unexceptional Gaussian prime with $\N(n) > \N(R)$; in particular,
this is true whenever $\frac{\epsilon^2 N^2}{2} \leq \N(n) \leq \epsilon^2 N^2$.
We now define the function $\nu: Z^2 \to \R^+$ by
\begin{equation}\label{nudef}
 \nu(n) \; := \; \left\{
\begin{array}{ll}
C_\varphi \frac{\phi_{\Z[i]}(W)}{\N(W)} \frac{\Lambda_{\Z[i]'_{sq},R,\varphi}(W\pi^{-1}(n) + b)^2}{\log \N(R)} & \hbox{ when } \N(\pi^{-1}(n)) \leq \epsilon^2 N^2\\
1 & \hbox{ otherwise}
\end{array}\right. 
\end{equation}
where $C_\varphi > 0$ is a normalization factor depending only on $\varphi$ to be chosen later (it is the constant which ensures that $\nu$ has mean close to 1), and $\pi^{-1}: Z^2 \to (-N/2,N/2)^2$ is the inverse of $\pi$ taking values in the fundamental domain $(-N/2,N/2)^2$.  
By construction we see that $\nu$ is non-negative and
$$ \nu(n) = C_\varphi  \frac{\phi_{\Z[i]}(W) \log \N(R)}{\N(W)} $$
whenever $n \in \pi(A_b)$.  In particular, from \eqref{acap-2} we have
$$ \E( \nu 1_{\pi(A_b)}(a) | a \in Z^2 ) \geq c_{\epsilon,\varphi} > 0.$$
Our task is now to show that
\begin{align*}
 \E( \prod_{0 \leq j < k} &(\nu 1_{\pi(A_b)})(a + r v_j) | a \in Z^2; r \in Z )\\
&\geq c_{\epsilon, \varphi} 
- o_{N \to \infty; \epsilon,\varphi}(1)
\end{align*}
for some $c_{\epsilon,\varphi} > 0$.
Since the $v_j$ were assumed to contain zero and span $\Z^2$, they will also span $Z^2$ if the prime $N$
is sufficiently large.  We thus see that the maps $\phi_j(r) := v_j r$ obey the ergodicity hypothesis
in Theorem \ref{sz-multi-alt-relative}.  We can invoke that theorem (in the contrapositive)
and be done as soon as we establish

\begin{proposition}[Existence of a system of pseudorandom majorants]\label{majorant} Consider the hypergraph system  $(J, (Z)_{j\in J}, d, H)$, where
$J := \{0,\ldots,k-1\}$, $d := k-2$, $H := {J \choose d}$, and for each  $e = J \backslash \{i\} \in H$, 
define the function $\nu_e: Z^e \to \R^+$ by
\begin{equation}\label{nue-def2}
\nu_e( (x_j)_{j \in e} ) := \nu( \sum_{j \in e} (v_j-v_i) x_j ).
\end{equation}
Then, the constant $C_\varphi$ is chosen properly, and if $\epsilon$ is sufficiently small depending on $k,v_1,\ldots,v_k$,
the system $(\nu_e)_{e \in H}$ is a 
pseudorandom system of measures, i.e. it obeys the dual function condition, the linear forms condition, and the correlation condition.
Note that we take $N$ to be the pseudorandomness parameter, and allow our bounds to depend on $\eps$, $k$, $v_1,\ldots,v_k$, $\varphi$.
\end{proposition}

It remains to prove the above proposition.  We shall do this in stages.  First we reduce matters from controlling various estimates involving
$\nu_e$ to estimates involving $\nu$.  More precisely, we will deduce Proposition \ref{majorant} from the following two propositions, whose
proof we shall give in later sections.  We first need some notation.

\begin{definition}[Gaussian $t$-tuples]  Let $T$ be a finite set.
If $\L = (L_t)_{t \in T} \in \Z[i]^T$ is a $T$-tuple of Gaussian integers,
and $\x = (x_t)_{t \in T} \in Z^T$ is a $T$-tuple of elements of $Z$, we define the quantity
$\L \cdot \x \in Z^2$ to be the quantity
$$ \L \cdot \x := \left(\sum_{t \in T} \Re(L_t) x_t, \sum_{j \in T} \Im(L_t) x_t\right).$$
We say that two $T$-tuples $\L$, $\L'$ are \emph{incommensurate} if they are both not identically zero, and we have
$\L \neq q \L'$ and $\L \neq q \overline{\L'}$ for any Gaussian rational $q \in \Q[i]$.  We say that $\L$ is \emph{self-incommensurate} if $\L \neq q \overline{\L}$ for any Gaussian rational $q \in \Q[i]$.
\end{definition}

The basic point here is that if $\L$, $\L'$ are non-degenerate and incommensurate then for any fixed $b, b' \in Z^2$ and some unknown $\x \in Z^T$, there is no obvious correlation between $\L \cdot \x + b$ being a Gaussian prime (or almost prime)
and between $\L' \cdot \x + b'$ being a Gaussian prime (or almost prime), other than those arising from small divisors (which
have already been eliminated through the $W$-trick).  The self-incommensurate hypothesis is needed to prevent the components of $\L$ from lying in a subspace of $\C$ (e.g. on the real axis), which would constrain $\L \cdot \x + b'$ to a line.

We now formalize the above heuristics.

\begin{proposition}[Linear forms condition for {$\Z[i]$}]\label{linear-forms-condition-gauss}
Let $S$ be a finite set of cardinality $|S| \leq k 2^k$, and $T$ be a finite set of cardinality $|T| \leq 2k$.
For each $s \in S$, let $\L_s \in \Z[i]^T$ be a $T$-tuple,
with any two $\L_s, \L_{s'}$ with $s \neq s'$ being incommensurate, and all $\L_s$ being self-incommensurate.
Then, if the exponent $c_{k}$ used to define $R$ is sufficiently small depending on $k$, $W$ is sufficiently large depending on $(\L_s)_{s \in S}$,  $C_\varphi$ in \eqref{nudef} is chosen correctly (depending only on $\varphi$), and $\epsilon$ is sufficiently small depending on
$(\L_s)_{s \in S}$, we have
\begin{equation}\label{lfc-gauss}
 \E( \prod_{s \in S} \nu(\L_s \cdot \x + b_s) | \x \in Z^T ) = 1 + o_{N \to \infty; \varphi,\epsilon, (\L_s)_{s \in S}}(1)
\end{equation}
uniformly for all choices of $(b_s)_{s \in S} \in (Z^2)^S$.
\end{proposition}

\begin{proposition}[Correlation condition for {$\Z[i]$}]\label{correlation-condition-gauss}  Let $m \leq 2^k$ and $v \in \Z[i] \backslash \{0\}$ 
be arbitrary.  Then, if $\epsilon$ is sufficiently small depending on $m,v$, there exists
functions $\tau^{(l)} = \tau^{(l)}_{v,m}: Z^2 \to \R^+$ for $l=1,2,3$ which are even (i.e. $\tau^{(l)}(-x) = \tau^{(l)}(x)$) which
obey the moment conditions
\begin{equation}\label{corr1}
\E( \tau^{(l)}(x)^q | x \in Z^2 ) = O_{q,v}(1) \hbox{ for } l=1,2
\end{equation}
and
\begin{equation}\label{corr0}
\E( \tau^{(3)}(0,x)^q | x \in Z ) = O_{q,v}(1) 
\end{equation}
for all integers $1 \leq q < \infty$, and furthermore will obey the moment conditions
\begin{equation}\label{corr2}
\E( \tau^{(l)}(v' \cdot x)^q | x \in Z ) = O_{q,v,v'}(1) \hbox{ for } l=1,2
\end{equation}
for all integers $1 \leq q < \infty$ and $w \in \Z[i] \backslash \{0\}$, if $\epsilon$ is sufficiently small depending on $m,v,v'$.
Furthermore we have the correlation estimate
\begin{equation}\label{corr-cond}
\begin{split}
& \E( \nu(v \cdot x + h_1) \ldots \nu(v \cdot x + h_m) | x \in Z ) \\
&\quad \leq \sum_{1 \leq i < j \leq m} \tau^{(1)}(h_i - h_j) + \tau^{(2)}(\overline{v} h_i - v \overline{h_j}) + \sum_{1 \leq j \leq m} \tau^{(3)}(\overline{v} h_j - v \overline{h_j})
\end{split}
\end{equation}
for all $h_1,\ldots,h_m \in Z^2$ (not necessarily distinct), where the conjugation operation $h \mapsto \overline{h}$ and the scalar multiplication
operation $h \mapsto vh$ on $Z^2$ are inherited from the corresponding operations on $\Z[i]$ in the obvious manner.
\end{proposition}

The $\tau^{(1)}$ term in \eqref{corr-cond} appeared in \cite{gt-primes}.  The $\tau^{(2)}$ term is new
and reflects the unavoidable fact that $\nu(x)$ and $\nu(\overline{x})$ will be very strongly correlated\footnote{At least, this is the case if $b$ is real.  If $b$ is complex then one has to shift $x$ or $\overline{x}$ by a fixed factor.}, since if $p$ is a Gaussian prime or almost prime then $\overline{p}$ will be also.  Similarly, the $\tau^{(3)}$ term is new and reflects the facts that $\nu$ will have an anomalous density on the real line (or on multiples of that line by $v$).  Note that while $\tau^{(3)}$ is ostensibly defined on $Z^2$, only its values on $0 \times Z$ are relevant, since this is where $\overline{v} h_j - v \overline{h_j}$ takes its values.

\begin{proof}[of Proposition \ref{majorant} assuming Proposition \ref{linear-forms-condition-gauss} and Proposition \ref{correlation-condition-gauss}]  We have to verify that $(\nu_e)_{e \in H}$ obeys the dual function condition (Definition \ref{dual-c}), linear forms condition (Definition \ref{lfc}) and the correlation condition (Definition \ref{corr-c}).  We begin with the dual function condition.  Fix $e = J \backslash \{i\}$ and $x^{(0)}_e \in V_e$. Using \eqref{defdef} to expand out $\D_e(\nu_e+1)$, and then using \eqref{nue-def2}, it suffices to show that
$$
\E( \prod_{\omega \in \Omega} \nu( \sum_{j \in e} (v_j - v_i) x^{(\omega_j)}_j) | x^{(1)}_e \in V_e ) = O(1)$$
for all $\Omega \subseteq \{0,1\}^e \backslash 0^e$.  But this follows from Proposition \ref{linear-forms-condition-gauss}; note that as the
$v_j$ are all distinct, each of the linear forms $\sum_{j \in e} (v_j - v_i) x^{(\omega_j)}_j$ utilizes a distinct non-empty subset of the variables
in $x^{(1)}_e$ and so the hypotheses of that Proposition are easily verified.

We now verify the linear forms condition. By \eqref{nue-def2}, it suffices to show that
\begin{equation}\label{lint}
\E( \prod_{(e,\omega) \in S} \nu(\sum_{j \in e} (v_j-v_i) x^{(\omega_j)}_j ) 
| x_J^{(0)}, x_J^{(1)} \in V_J ) = 1 + o_{W \to \infty}(1)
\end{equation}
for any finite set $S$ of pairs $(e,\omega)$ such that $e \in H$,
$\omega \in \{0,1\}^e$.  We can parameterize the averaging variables
by $(x_t)_{t \in T} \in Z^T$, where $T$ is the finite set
$T = J \times \{0,1\}$.
We can thus write the left-hand side of \eqref{lint} as
$$\E\left( \prod_{s \in S} \nu(\L_s \cdot \x + b_s) \bigl| \x \in Z^T \right)$$
where for any $(e,\omega) \in S$ with $e = J \backslash \{i\}$, we have
$$ \L_{e,\omega} := ( 1_{\omega_j = a} (v_j-v_i) )_{(j,a) \in T}; \quad b_{e,\omega} := \sum_{j \in e \cap J': \omega_j = 0} (v_j-v_i) x_j.$$
The hypothesis that the $v_j$ are all distinct ensures
that the $\L_{e,\omega}$ are non-zero.  In fact they are all pairwise incommensurate, because each $\L_{e,\omega}$ has a different set of non-zero
co-ordinates.  Because the $v_j-v_i$ span $\Z[i]$, we also see that each $\L_{e,\omega}$ is self-incommensurate.
Thus \eqref{lint} follows from Proposition \ref{linear-forms-condition-gauss}, if $\epsilon$ is sufficiently small depending
on the $\L_{e,\omega}$, which in turn depend only on the $k, v_1,\ldots,v_k$.

We now turn to the correlation condition.  Fix $e = J \backslash \{i\} \in H$, $j \in e$, $K \geq 0$, and $n_{e,\omega} \in \{0,1\}$.
The left-hand side of \eqref{correlation} can be expanded as
$$ \E\left( \prod_{a=0}^1 \E\left( \prod_{\omega \in A_a} \nu( (v_j-v_i) x_j + h_\omega ) \bigr| x_j \in Z\right)^K \biggr| x_{e \backslash \{j\}}^{(0)}, x_{e \backslash \{j\}}^{(1)} \in Z^{e \backslash \{j\}} \right)$$
where $A_a := \{ \omega \in \{0,1\}^e: n_{e,\omega} = 1; \omega_j = a\}$ and
$$ h_\omega := \sum_{j' \in e \backslash \{j\}} (v_{j'}-v_i) x^{(\omega_{j'})}_{j'} .$$
By Cauchy-Schwarz and symmetry it suffices to show that
$$ \E( \E( \prod_{\omega \in A_0} \nu( (v_j-v_i) x_j + h_\omega ) | x_j \in Z)^{2K} | x_{e \backslash \{j\}}^{(0)}, x_{e \backslash \{j\}}^{(1)} \in Z^{e \backslash \{j\}} ) = O_K(1).$$
Applying Proposition \ref{correlation-condition-gauss} with $v := v_j - v_i \neq 0$, we have
\begin{align*}
&\E( \prod_{\omega \in A_0} \nu( (v_j-v_i) x_j + h_\omega ) | x_j \in Z)\\ 
&\quad \leq \sum_{\omega, \omega' \in A_0: \omega \neq \omega'} \tau^{(1)}(h_\omega - h_{\omega'}) + 
\tau^{(2)}(\overline{v} h_\omega - v \overline{h_{\omega'}}) + \sum_{\omega \in A_0} \tau^{(3)}(\overline{v} h_\omega - v \overline{h_\omega})
\end{align*}	
where $\tau^{(l)} = \tau^{(l)}_{|A_0|,v}$, assuming of course that $\epsilon$ is sufficiently small depending on $k,v_1,\ldots,v_k$.  
By the triangle inequality, we can thus bound the left-hand side of \eqref{correlation} by
\begin{align*}
 O_K\biggl( \sum_{\omega, \omega' \in A_0: \omega \neq \omega'} &\E\bigl( \tau^{(1)}(h_{\omega} - h_{\omega'})^{2K} 
+
\tau^{(2)}(\overline{v} h_\omega - v \overline{h_{\omega'}})^{2K}\\
&+ \sum_{\omega \in A_0} \tau^{(3)}(\overline{v} h_\omega - v \overline{h_\omega})^{2K}
\bigr|
x_{e \backslash \{j\}}^{(0)}, x_{e \backslash \{j\}}^{(1)} \in Z^{e \backslash \{j\}} \bigr) \biggr).
\end{align*}
Since $|A_0| = O(1)$, it thus suffices to show that
$$\E( \tau^{(1)}(h_{\omega} - h_{\omega'})^{2K} +
\tau^{(2)}(\overline{v} h_\omega - v \overline{h_{\omega'}})^{2K}
+ \tau^{(3)}(\overline{v} h_\omega - v \overline{h_\omega})^{2K}
|
x_{e \backslash \{j\}}^{(0)}, x_{e \backslash \{j\}}^{(1)} \in Z^{e \backslash \{j\}} ) = O_K(1)$$
for any distinct $\omega, \omega'$ in $A_0$.

Fix $\omega,\omega' \in A_0$.  Let us first deal with the
$\tau^{(1)}(h_{\omega} - h_{\omega'})^{2K}$ term.
Observe that the map from $(x_{e \backslash \{j\}}^{(0)}, x_{e \backslash \{j\}}^{(1)})$ to
$h_{\omega} - h_{\omega'}$ is a group homomorphism from $Z^{e \backslash \{j\}} \times Z^{e \backslash \{j\}}$ to $Z^2$.
Since $Z$ is a cyclic group of prime order, we thus see that the image of this homomorphism is either $\{0\}$, $Z^2$, or a line of the form
$\{ v' \cdot x: x \in Z \}$ for some $v' \in \Z[i]$.  Also, all the fibers of this group homomorphism have the same cardinality (they are all cosets
of the same kernel).
Since all the $v_j$ are distinct and $\omega \neq \omega'$, we see that the image is
not zero.  If it is $Z^2$ then the claim now follows from \eqref{corr1}.  If the image is a line, then the Gaussian integer $v' \in \Z[i]$ depends only
on $k,v_1,\ldots,v_k,\omega,\omega'$.  Since the number of values of $\omega, \omega'$ is $O(1)$, we thus see that if $\epsilon$ is small
enough depending on $k,v_1,\ldots,v_k$, the claim will now follow from \eqref{corr2}.

The contribution of the $\tau^{(2)}(\overline{v} h_\omega - v \overline{h_{\omega'}})^{2K}$ term is dealt with similarly; note
that the map from $(x_{e \backslash \{j\}}^{(0)}, x_{e \backslash \{j\}}^{(1)})$ to
$\overline{v} h_\omega - v \overline{h_{\omega'}}$ is still a group homomorphism whose image is not identically zero.

Finally, we control the contribution of $\tau^{(3)}$.  Here we use the improved ergodic hypothesis, Hypothesis \ref{super-ergodic}.
This implies that the map from $(x_{e \backslash \{j\}}^{(0)}, x_{e \backslash \{j\}}^{(1)})$ to
$h_\omega$ is a \emph{surjective} group homomorphism, and hence the map to $\overline{v} h_\omega - v \overline{h_\omega}$
has image $0 \times Z$. Thus the contribution of $\tau^{(3)}$ can be controlled purely by \eqref{corr0}.
\end{proof}

\section{Reduction to a number-theoretic estimates}

To conclude the proof of Theorem \ref{main}, we have to verify Proposition \ref{linear-forms-condition-gauss} and
Proposition \ref{correlation-condition-gauss}.  These propositions are estimates on the function $\nu$, which was defined
in \eqref{nudef}, partly in terms of the truncated divisor sum $\Lambda_{\Z[i]'_{sq},R,\varphi}$ and partly in terms of the constant
function 1.  In this section we reduce matters purely to estimation of the truncated divisor sum.  In particular we reduce Proposition \ref{linear-forms-condition-gauss} to the following estimate.

\begin{proposition}[First Goldston-Y{\i}ld{\i}r{\i}m correlation estimate for {$\Z[i]'_{sq}$}]\label{GY-Gauss}  
Let $S, T$ be finite sets. For each $s \in S$, let $\L_s \in \Z[i]^T$ be a $T$-tuple,
with any two $\L_s, \L_{s'}$ with $s \neq s'$ being incommensurate, and all $\L_s$ being self-incommensurate.  For each $s \in S$, let $a_s \in \Z[i]$ be a Gaussian integer coprime to $W$.
Let $B \subset \Z^T$ is a product $B = \prod_{t \in T} I_t$ of $t$ intervals $I_t \subseteq \Z$, each of length at least $R^{10|S|}$.
Then, if $W$ is sufficiently large depending on $(\L_s)_{s \in S}$, we have
\begin{equation}\label{gy-cor}
\begin{split}
 \E&\left( \prod_{s \in S}  \Lambda_{\Z[i]'_{sq},R,\varphi}(W (\L_s \cdot \x) + a_s)^2 | \x \in B \right) \\
 &= 
(1 + o_{R \to \infty; |S|,|T|,\varphi,W, (\L_s)_{s \in S}}(1) + o_{W \to \infty; |S|,|T|,\varphi, (\L_s)_{s \in S}}(1)) \\
&\quad \left(c_\varphi \frac{\N(W) \log \N(R)}{\phi_{\Z[i]}(W)}\right)^{|S|}
\end{split}
\end{equation}
for an explicit quantity $c_\varphi > 0$ depending only on $\varphi$.
\end{proposition}

Similarly, we will reduce Proposition \ref{correlation-condition-gauss} to the following estimate.

\begin{proposition}[Second Goldston-Y{\i}ld{\i}r{\i}m correlation estimate for {$\Z[i]'_{sq}$}]\label{GY-Gauss-2} 
Let $m$ be a positive integer, let $b$ be a Gaussian integer coprime to
$W$, and let $v$ be a Gaussian integer.  Let $h_1,\ldots,h_m$ be Gaussian integers such that the quantity
\begin{equation}\label{Delta-def}
 \Delta := \prod_{1 \leq i < j \leq m} \N(h_i - h_j) 
 \prod_{1 \leq i \leq j \leq m} \N(W(h_i \overline{v} - \overline{h_j} v) - b \overline{v} + \overline{b} v).
 \end{equation}
is non-zero.
Let $I \subset \Z$ be an interval of length at least $R^{10m}$.  Then, if $W$ is
sufficiently large depending on $v$, we have
\begin{equation}\label{parade}
\begin{split}
& \E\left( \prod_{j=1}^m \Lambda_{\Z[i]'_{sq},R,\varphi}(W(h_j + nv) + b)^2 \bigg| n \in I \right)
\\
&\quad \leq O_{m,v}\left( 
\left(\frac{\N(W) \log \N(R)}{\phi_{\Z[i]}(W)}\right)^m 
\prod_{p \in \P[i]'_{+,W}: p | \Delta} (1 + O_m( \N(p)^{-1/2} )) \right)
\end{split}
\end{equation}
where $\P[i]'_{+,W}$ are those primes in $\P[i]'_+$ which are coprime to $W$.
\end{proposition}

The presence of the rather unusual expression $W(h_i \overline{v} - \overline{h_j} v) - b \overline{v} + \overline{b} v$ in
\eqref{Delta-def} can be partially explained
by the following observation: if $W(h_i \overline{v} - \overline{h_j} v) - b \overline{v} + \overline{b} v = 0$,
then we have
$$ \overline{v} (W(h_i + nv) + b) = v \overline{(W(h_j + nv) + b)}$$
for all $n$.  Thus there is likely to be a strong correlation between $W(h_i + nv) + b$ being prime or almost prime, and
$W(h_j + nv) + b$ being prime or almost prime.  This correlation also occurs in the diagonal case $i=j$, reflecting the fact that $\Lambda_{\Z[i]'_{sq},R,\varphi}$ is substantially larger on the real line (and on multiples of the real line by gaussian rationals of small height)
than in general.  Thus we expect the left-hand side of \eqref{parade} to be abnormally large
when $\Delta$ is zero; it turns out that it can also be large when $\Delta$ is very smooth (has many small prime factors).

We now show how these propositions imply Propositions \ref{linear-forms-condition-gauss} and \ref{correlation-condition-gauss}.

\begin{proof}[Proof of Proposition \ref{linear-forms-condition-gauss} from Proposition \ref{GY-Gauss}]
This shall follow the proof of \cite[Proposition 9.8]{gt-primes}; the main idea is to discretize the domain $Z^T$ to the point where the boundary
effects caused by the constraint $\N(n) \leq \epsilon^2 N^2$ in \eqref{nudef} are negligible.  In this proof we allow all constants to depend on $|S|$ and $|T|$, $\varphi$, and $(\L_s)_{s \in S}$.

Let us view $Z^T$ as the discrete cube $(-N/2,N/2)^T$.  If the constant $c_k$ used to define $R$ is sufficiently small, we can
find an integer $Q = Q(N)$ such that $N/Q \geq 2 R^{10|S|}$ and $1/Q = o_{N \to \infty}(1)$ (thus $Q$ grows slowly with $N$).
We partition $(-N/2,N/2)^T$ into $Q^{|T|}$ boxes $(B_\alpha)_{\alpha \in A}$, each of sidelength $N/Q (1 + o_{N \to \infty}(1))$.
Then, up to multiplicative errors of $1 + o_{N \to \infty}(1)$, the left-hand side of \eqref{lfc-gauss} is equal to
$$
 \E\left( \E( \prod_{s \in S} \nu(\L_s \cdot \x + b_s) | \x \in B_\alpha ) \bigr| \alpha \in A \right).
$$
It thus suffices to show that
$$
 \E\left( \E( \prod_{s \in S} \nu(\L_s \cdot \x + b_s) | \x \in B_\alpha ) - 1 \bigr| \alpha \in A \right) = o_{N \to \infty}(1).
$$
Writing $\nu = 1 + (\nu-1)$ and expanding, it suffices to show that
\begin{equation}\label{nusub}
 \E\left( \E( \prod_{s \in S'} (\nu(\L_s \cdot \x + b_s)-1) | \x \in B_\alpha ) \bigr| \alpha \in A \right) = o_{N \to \infty}(1)
\end{equation}
for all non-empty $S' \subseteq S$.  

Fix $S'$.  Let $D$ denote the disk $D := \{ n \in \Z[i]: \N(n) \leq \epsilon^2 N^2\}$. 
We divide the boxes $B_\alpha$ into three categories.  We say that a box $B_\alpha$ is \emph{interior} if $\L_s \cdot \x + b_s \in \pi(D)$ for all $s \in S'$ and $\x \in B_\alpha$.  We say that a box $B_\alpha$ is \emph{exterior} if there exists an $s \in S'$
such that $\L_s \cdot \x + b_s \not \in \pi(D)$ for all $\x \in B_\alpha$.  We say that a box $B_\alpha$ is \emph{borderline of type $s$}
for some $s \in S'$ if $\L_s \cdot \x + b_s \in \pi(D)$ for at least one $\x \in B_\alpha$, and
$\L_s \cdot {\mathbf y} + b_s \not \in O$ for at least one ${\mathbf y} \in B_\alpha$.  Clearly every box is either interior, exterior, or borderline for some $s \in S'$.
From \eqref{nudef}, the exterior boxes give a zero contribution to \eqref{nusub}.  Now consider an interior box $B_\alpha$.  For these boxes we claim that
$$ \E( \prod_{s \in S'} (\nu(\L_s \cdot \x + b_s)-1) | \x \in B_\alpha ) = o_{N \to \infty}(1).$$
Expanding out the product using the binomial formula, it suffices to show that
$$ \E( \prod_{s \in S''} \nu(\L_s \cdot \x + b_s) | \x \in B_\alpha ) = 1 + o_{N \to \infty}(1)$$
for all $S'' \subseteq S'$.  At this point we need to make a technical remark concerning the identification between elements of $\Z[i]$
and elements of $Z^2$, and between $Z^T$ and $(-N/2,N/2)^T$.  Currently, $\x$ is viewed as an element of $Z^T$, and $b_s$ is an element
of $Z^2$, and so $\L_s \cdot \x + b_s$ is also an element of $Z^2$.  But using $\pi$, this element of $Z^2$ is then considered
to be an element of $\Z[i]$, which in fact lies in the disk $D$.

We now change this perspective, viewing $\x$ now as an element of $(-N/2,N/2)^T$ (and $B_\alpha$ as a box of sidelengths $\approx N/Q$ inside
$(-N/2,N/2)^T$).  This makes $\L_s \cdot \x$ an element of $\Z[i]$ rather than $Z^2$, although the dimensions of the box $B_\alpha$ will keep
$\L_s \cdot \x$ constrained to a ball of radius $O_{\L_s}(N/Q)$.  We now wish to view $b_s$ as an element of $\Z[i]$ also, but one has the freedom to
modify $b_s$ by an element of $N \Z[i]$ in doing so.  However, only one of these ``lifts'' of $b_s$ will place $\L_s \cdot \x$ to lie in $D$.  
Indeed, since $\L_s \cdot \x$ is constrained to a ball of radius much less than $N$, there is
a unique lift of $b_s$ in $\Z[i]$ (which by abuse of notation we shall continue to call $b_s$), 
independent of the choice of $\x$, for which $\L_s \cdot x + b_s$ lies in $D$, now viewed as a 
subset of $\Z[i]$ rather than $Z^2$.  Applying \eqref{nudef}, we can now write the left-hand side as
$$ \frac{C_\varphi \phi_{\Z[i]}(W)}{\N(W) \log \N(R)} 
\E( \prod_{s \in S''} \Lambda_{\Z[i]'_{sq},R,\varphi}(W \L_s \cdot \x + a_s) | \x \in B_\alpha )$$
where $a_s := Wb_s + b$.  Applying Proposition \ref{GY-Gauss} and choosing $C_\varphi := 1/c_\varphi$, this expression is equal to
$$ 1 + o_{R \to \infty; W}(1) + o_{W \to \infty}(1)$$
which is acceptable since $R$ is a small power of $N$, and $W$ is chosen to grow extremely slowly in $R$.

It remains to control the contribution of the borderline boxes of type $s_0$ for some $s_0 \in S'$.  Since $|S'| = O(1)$, we may fix $s_0$.  
Bounding $\nu-1$ in absolute value by $\nu+1$, and using several applications of Proposition \ref{GY-Gauss}, we can control
$$ \E( \prod_{s \in S'} (\nu(\L_s \cdot \x + b_s)-1) | \x \in B_\alpha ) $$
crudely by $O(1)$.  To conclude the proof it suffices to show that the number of borderline boxes of type $s_0$ is small, in the sense
that it is $o_{N\to \infty}(1)$ times the total number $Q^{|T|}$ of boxes.

Observe that the set $\{ \L_{s_0} \cdot \x + b_{s_0}: \x \in B_\alpha \}$ 
has a diameter of $O(N/Q)$, where the metric on $Z^2$ is the quotient metric inherited from $\Z[i]$.  
Since $B_\alpha$ is borderline of type $s_0$, we conclude that
$$ \{ \L_{s_0} \cdot \x + b_{s_0}: \x \in B_\alpha \} \subseteq \{ n \in \Z[i]: |n| = N + O(N/Q) \},$$
where the annulus on the right-hand side is thought of as a subset of $Z^2$.  Next, observe that the map
$\x \mapsto \L_{s_0} \cdot \x + b_{s_0}$ is an affine homomorhpism from $Z^T$ to $Z^2$, and thus
(since $Z$ has prime order) the image is an affine subspace of $Z^2$, with the fibers at each point of this image having equal cardinality.  Since $\L_{s_0}$ is not identically zero, the image is either an affine line in $Z^2$ (with a ``slope'' determined entirely by $\L_{s_0}$)
or is all of $\Z^2$.  In either case, we see from elementary geometry that the proportion of points in this image which lie
in the annulus $ \{ n \in \Z[i]: |n| = N + O(N/Q) \}$ is $o_{Q \to \infty}(1)$ (indeed one can obtain the more precise
bound of $O(Q^{-1/2})$, because the circle $\{ |n| = N \}$ has non-vanishing curvature, though we will not need that improved
bound here).  Thus the proportion of boxes which are borderline of type $s_0$ is $o_{Q \to \infty}(1)$, which is acceptable since
$Q$ is growing with $N$.
\end{proof}

\begin{proof}[of Proposition \ref{correlation-condition-gauss} from Proposition \ref{GY-Gauss-2}]
The arguments here are somewhat similar to the derivation of Proposition \ref{linear-forms-condition-gauss} from Proposition \ref{GY-Gauss}
but are simpler because we are only seeking upper bounds rather than asymptotics.  On the other hand, some number theory is required
to control the expressions arising from Proposition \ref{GY-Gauss-2}.

Fix $m,v$.
We first observe that we may remove the requirement that $\tau^{(l)}$ is even, since we may simply replace $\tau^{(l)}(x)$ by 
$\tau^{(l)}(x) + \tau^{(l)}(-x)$ if necessary.  

We begin by establishing the very crude estimate
\begin{equation}\label{crudo}
 \E( \prod_{j=1}^m \nu(v \cdot x + h_j) | x \in Z ) =
O_{m,\varphi} \left( (\log^m N) \sup_{n \in \Z[i]: |n| \leq NW} d_{\Z[i]}(n)^{2m} \right),
\end{equation}
where $d_{\Z[i]}(n)$ is the Gaussian divisor function
$$ d_{\Z[i]}(n) := \sum_{d \in \Z[i]: d | n} 1.$$
To see this, we first use H\"older's inequality to bound the left-hand side of
\eqref{corr-cond} very crudely by
$$ \left(\sup_{x \in Z^2} \nu(x)\right)^m.$$
By \eqref{nudef} and \eqref{trunc-def}, this can in turn be crudely estimated by
$$ O_{m,\varphi}\left(1 + \left(\frac{\phi_{\Z[i]}(W) \log \N(R)}{\N(W)}\right)^m \sup_{n \in \Z[i]: \N(n) \leq \epsilon^2 N^2} \left[\sum_{d \in \Z[i]'_{sq}: d|Wn+b} 1\right]^{2m} \right).$$
The claim \eqref{crudo} follows.  Next, observe that $d_{\Z[i]}(n) = O_\eps( \N(n)^\eps )$
whenever $n$ is a power of a Gaussian prime $p$, and can in fact improve this to $d_{\Z[i]}(n) \leq \N(n)^\eps$ if $\N(p)$ is sufficiently large depending on $\eps$.  Using the multiplicativity of $d_{\Z[i]}(n)$ we conclude that $d_{\Z[i]}(n) = O_\eps( \N(n)^\eps )$ for all $n$, and so
we see that the right-hand side of \eqref{crudo} is $O_{m,\varphi,\eps}(N^\eps)$ for any $\eps > 0$.

In light of \eqref{crudo}, we will define $\tau^{(1)}(0)$, $\tau^{(2)}(W^{-1}(b \overline{v} - \overline{b} v) )$ and
$\tau^{(3)}(W^{-1}(b \overline{v} - \overline{b} v) )$ to equal the
right-hand side of \eqref{crudo}, and observe from the preceding discussion that this will not significantly affect \eqref{corr1}, \eqref{corr0} or
\eqref{corr2}. 

It now remains to treat the cases when $h_i - h_j \neq 0$
for all $1 \leq i < j \leq m$ and $W(h_i \overline{v} - \overline{h_j} v) - b \overline{v} + \overline{b} v \neq 0$ for all $1 \leq i \leq j \leq m$.
In other words, we are left with the case where the quantity $\Delta$ defined in \eqref{Delta-def} is non-zero.
We now use \eqref{nudef} to crudely estimate
$$ \nu(x) \leq 1 + C_\varphi \frac{\phi_{\Z[i]}(W)}{\N(W)} \frac{\Lambda_{\Z[i]'_{sq},R,\varphi}(W\pi^{-1}(n) + b)^2}{\log \N(R)} 1_{\N(\pi^{-1}(n)) \leq \epsilon^2 N^2}
$$
and expand terms, to reduce to establishing an estimate of the form
\begin{align*}
&\E\left( \prod_{j=1}^m \Lambda_{\Z[i]'_{sq},R,\varphi}(W\pi^{-1}(v \cdot x + h_j) + b)^2 1_{\N(\pi^{-1}(v \cdot x + h_j)) \leq \epsilon^2 N^2} \bigr| x \in Z \right)\\
&\quad \leq 
\left(\frac{\N(W) \log \N(R)}{\phi_{\Z[i]}(W)}\right)^m 
\left[\sum_{1 \leq i < j \leq m} \tau^{(1)}(h_i - h_j) + \tau^{(2)}( h_i \overline{v} - \overline{h_j} v )
+ \sum_{1 \leq i \leq m} \tau^{(3)}(h_i \overline{v} - \overline{h_i} v)\right]
\end{align*}
for functions $\tau^{(1)},\tau^{(2)}, \tau^{(3)}$ obeying \eqref{corr1}, \eqref{corr0}, \eqref{corr2}, in the case $\Delta \neq 0$.

Using the identity
$$ \E( f(x) | x \in Z ) = \E( \E( f(x + n) | 1 \leq n \leq N^{1/2} ) | x \in Z )$$ 
we see it suffices to obtain an estimate of the form
\begin{align*}
& \E\left( \prod_{j=1}^m \Lambda_{\Z[i]'_{sq},R,\varphi}(W\pi^{-1}(v \cdot (x+n) + h_j) + b)^2 1_{\N(\pi^{-1}(v \cdot (x+n) + h_j)) \leq \epsilon^2 N^2} \bigr| 
1 \leq n \leq N^{1/2} \right)\\
&\quad \leq 
\left(\frac{\N(W) \log \N(R)}{\phi_{\Z[i]}(W)}\right)^m 
\left[\sum_{1 \leq i < j \leq m} \tau^{(1)}(h_i - h_j) + \tau^{(2)}( h_i \overline{v} - \overline{h_j} v )
+ \sum_{1 \leq i \leq m} \tau^{(3)}(h_i \overline{v} - \overline{h_i} v)
\right]
\end{align*}
uniformly in $x$.  By absorbing $v \cdot x$ into the $h_j$ term we may take $x=0$.  We then observe that this sum is zero unless
there exists an $n$ for which $|\pi^{-1}(v \cdot n + h_j)| \leq \epsilon N$ for all $j$.  In particular this forces
$$ |\pi^{-1}(h_j)| \leq \epsilon N + O( |v| N^{1/2} ) = 2 \epsilon N$$
if $N$ is sufficiently large depending on $\epsilon$ and $|v|$.  In particular, by the triangle inequality,
$\tau$ only needs to be defined on the region $\{ x \in Z^2: |\pi^{-1}(x)| \leq 4 \epsilon N \}$.
Now observe that $\pi^{-1}(v \cdot n + h_j) = \pi^{-1}(h_j) + nv$ (if $\epsilon$ is sufficiently small to avoid wraparound issues).  Setting
$\tilde h_j := \pi^{-1}(h_j)$, we thus
reduce to showing that
\begin{equation}\label{lrw}
\begin{split}
& \E\left( \prod_{j=1}^m \Lambda_{\Z[i]'_{sq},R,\varphi}(W(\tilde h_j + nv) + b)^2 
\bigr| 1 \leq n \leq N^{1/2} \right)\\
&\quad \leq 
\left(\frac{\N(W) \log \N(R)}{\phi_{\Z[i]}(W)}\right)^m 
\left[\sum_{1 \leq i < j \leq m} \tilde \tau_1(\tilde h_i - \tilde h_j) + \tilde \tau_2( \tilde h_i \overline{v} - \overline{\tilde h_j} v )
+ \sum_{1 \leq i \leq m} \tilde \tau_3(\tilde h_i \overline{v} - \overline{\tilde h_i} v) \right]
\end{split}
\end{equation}
for all distinct $\tilde h_j \in \{ x \in \Z[i]: |x| \leq 2 \epsilon N \}$, and for functions $\tilde \tau_1, \tilde \tau_2, \tilde \tau_3: 
\Z[i]\backslash \{0\} \to \R^+$ supported
on the punctured disk $D := \{ x \in \Z[i]: 0 < |x| \leq 4 \epsilon N \}$, such that the functions $\tau^{(l)} := \tilde \tau_l \circ \pi^{-1}$ obeys
\eqref{corr1}, \eqref{corr0}, \eqref{corr2}.

Applying Proposition \ref{GY-Gauss-2}, and assuming that $W$ is large enough depending on $v$, and the exponent $c_k$ used to define $R$
is sufficiently small, we can bound the left-hand side of \eqref{lrw} by
\begin{align*}
O_{m,v}\biggl( &\left(\frac{\N(W) \log \N(R)}{\phi_{\Z[i]}(W)}\right)^m \\
&\prod_{1 \leq i < j \leq m} \prod_{p \in \P[i]'_{+,W}: p | \N(\tilde h_i - \tilde h_j) 
\N(W(\tilde h_i \overline{v} - \overline{\tilde h_j} v) - b \overline{v} + \overline{b} v)
} (1 + O_m( \N(p)^{-1/2} )) \\
&\prod_{1 \leq i \leq m} \prod_{p \in \P[i]'_{+,W}: p |
\N(W(\tilde h_i \overline{v} - \overline{\tilde h_i} v) - b \overline{v} + \overline{b} v)
} (1 + O_m( \N(p)^{-1/2} )) \biggr)
\end{align*}
and so by the arithmetic mean-geometric mean inequality we will be able to satisfy \eqref{lrw} by setting
$$ \tilde \tau_1(x) := 1_D(x) O_m( \prod_{p \in \P[i]'_{+,W}: p | \N(x)} (1 + O_m( \N(p)^{-1/2} ))^{m^2} )$$
and
$$ \tilde \tau_2(x) = \tilde \tau_3(x) := 1_D(x) O_m( \prod_{p \in \P[i]'_{+,W}: p | \N(Wx - b \overline{v} + \overline{b} v)} (1 + O_m( \N(p)^{-1/2} ))^{m^2} )$$

We now need to verify \eqref{corr1}, \eqref{corr0}, and \eqref{corr2}.  Let us first verify \eqref{corr1} for $\tau^{(1)}$.  We need to show that
$$ \sum_{x \in D} \prod_{p \in \P[i]'_+: p | \N(x)} (1 + O_m( \N(p)^{-1/2} ))^{m^2q} = O_{m,q}(N^2).$$
Estimating
\begin{equation}\label{pax}
\begin{split}
 \prod_{p \in \P[i]'_{+,W}: p | \N(x)} (1 + O_m( \N(p)^{-1/2} ))^{m^2 q}
&\leq \prod_{p \in \P[i]'_{+,W}: p | \N(x)} (1 + O_{m,q}( \N(p)^{-1/2} ) \\
&\leq O_{m,q}( \prod_{p \in \P[i]'_{+,W}: p | \N(x)} (1 + \N(p)^{-1/4} ) )\\
&\leq O_{m,q}( \sum_{d \in \Z[i]'_{sq,W}: d | \N(x)} \N(d)^{-1/4} )
\end{split}
\end{equation}
where $\Z[i]'_{sq,W}$ are those elements of $\Z[i]'_{sq}$ which are coprime to $W$.
We then see that
\begin{equation}\label{pax-2}
\begin{split}
\sum_{x \in D} \prod_{p \in \P[i]'_+: p | \N(x)} (1 + O_m( \N(p)^{-1/2} ))^{m^2 q} 
&\leq O_{m,q}( \sum_{x \in D} \sum_{d \in \Z[i]'_{sq,W}: d|\N(x)} \N(d)^{-1/4} )\\
&= O_{m,q}( \sum_{d \in \Z[i]'_{sq,W}} \N(d)^{-1/4} |\{ x \in D: d | \N(x) \}|).
\end{split}
\end{equation}
Since $d \in \Z[i]'_{sq,W}$, we have $d = p_1 \ldots p_k$ for some distinct (non-associate) $p_1,\ldots,p_k \in \P[i]'$.
we see that if $d|\N(x)$, then $d'|\N(x)$, where $d' = p'_1 \ldots p'_k$ and each $p'_j$ is either associate to $p_j$ or to $\overline{p_j}$.
There are $O(2^k)$ possible values of $d'$, and they all have the same norm as $d$.  We thus see that
$$ |\{ x \in D: d | \N(x)\}| = O( 2^k |D| / \N(d) ) = O( 2^k N^2 / \N(d) ).$$
Since there are only finitely many Gaussian primes in any given bounded set, we see that $2^k = O( \N(d)^{1/8} )$ (for instance).
Thus we have
\begin{align*}
\sum_{x \in D} \prod_{p \in \P[i]'_+: p | \N(x)} (1 + O_m( \N(p)^{-1/2} ))^{mq} 
&= O_{m,q}( \sum_{d \in \Z[i]'_{sq,W}} \N(d)^{-1/4} \N(d)^{1/8} N^2 / \N(d) )\\
&= O_{m,q}( N^2 \sum_{d \in \Z[i] \backslash \{0\}} \N(d)^{-9/8} )\\
&= O_{m,q}(N^2)
\end{align*}
as desired.

Now we verify \eqref{corr2} for $\tau^{(1)}$.  If $v'$ is fixed,
and $W$, $N$ are large with respect to $v'$, then the set $\{ x \in Z: \pi^{-1}(v' \cdot x) \in D \}$ is essentially a union of
$O_{v'}(1)$ intervals of length $O(N)$, on which $\pi^{-1}(v' \cdot x)$ is an arithmetic progression of step $r := \pi^{-1}(v')$.  It thus suffices
to show that
$$ \sum_{j \in \Z: j=O(N), a+jr \neq 0} \prod_{p \in \P[i]'_+: p | \N(a+jr)} (1 + O_m( \N(p)^{-1/2} ))^{m^2 q} = O_{m,q,r}(N)$$
where $a = O(N)$ is a Gaussian integer.  Applying \eqref{pax} again, we bound the left-hand side by
\begin{align*}
&O_{m,q}\left( \sum_{j \in \Z: j=O(N), a+jr \neq 0} \sum_{d \in \Z[i]'_{sq}: d | \N(a+jr)} \N(d)^{-1/4} \right)\\
= &O_{m,q}\left( \sum_{d \in \Z[i]'_{sq}} \N(d)^{-1/4} |\{ j \in \Z: j = O(N), d|\N(a+jr), a+jr \neq 0 \}| \right). 
\end{align*}
We can assume that $d = O_r(N)$ since the larger values $d$ give a zero contribution (recall that $a = O(N)$).
As before, we can replace the constraint $d | \N(a+jr)$ by $d' | a+jr$, where $d'$ ranges over $O(2^k) = O(\N(d)^{1/8})$ possible values,
all with norm equal to $\N(d)$.
The Gaussian integer $d'$ is (up to Gaussian units) the product of primes $p$ in $P[i]_+$, with no prime appearing at most once.
For each of these primes, the group $\Z[i]/p\Z[i]$ is a cyclic group of prime order, thus has no proper subgroups.  From this fact
and the Chinese remainder theorem for Gaussian integers, we see that $|\{ j \in \Z: d'|j \}| = \N(d) \Z$.  We can thus estimate the previous
expression by
$$ O_{m,q}\left( \sum_{d \in \Z[i]'_{sq}: d = O_r(N)} \N(d)^{-1/4} \N(d)^{1/8} O( 1 + \frac{N}{\N(d)} ) \right) = O_{m,q,r}(N)$$
as desired.

Now we verify \eqref{corr1} for $\tau^{(2)}$.  We need to show that
$$ \sum_{x \in D} \prod_{p \in \P[i]'_+: p | \N(Wx - b \overline{v} + \overline{b} v) \neq 0} (1 + O_m( \N(p)^{-1/2} ))^{m(m-1)q} = O_{m,q}(N^2).$$
By modifying the computations in \eqref{pax}, \eqref{pax-2} we see the left-hand side is
$$ O_{m,v}(\sum_{d \in \Z[i]'_{sq,W}} \N(d)^{-1/4} |\{ x \in D: d | \N(Wx - b \overline{v} + \overline{b} v \neq 0\}| ).$$
As before, we see that if $d$ divides $\N(Wx-b\overline{v} + \overline{b} v)$, then $d'$ divides $Wx - b \overline{v} + \overline{b} v$,
where $d'$ ranges over $O(\N(d)^{1/8})$ Gaussian integers in $\Z[i]'_{sq,W}$ with the same norm as $d$.  Since the radius of $D$ is large
compared with $d$, $W$, $b$, or $v$, we see (using the Chinese remainder theorem, since $d'$ and $W$ are coprime)
that for any fixed $d'$, the number of elements $x$ of $D$ for which $d'$ divides $Wx - b \overline{v} + \overline{b} v$ is
$O( N^2 / \N(d') ) = O( N^2 / \N(d) )$.  The proof of \eqref{corr1} then proceeds as with $\tau^{(1)}$.   For similar reasons we can adapt
the proof of \eqref{corr2} for $\tau^{(1)}$ to also give a proof for $\tau^{(2)}$, which then also implies \eqref{corr0}.
\end{proof}

\section{Proof of Proposition \ref{GY-Gauss}}\label{Gaussian-gy}

To conclude the proof of Theorem \ref{main}, we need to prove Proposition \ref{GY-Gauss} and Proposition \ref{GY-Gauss-2}.  
This is the purpose of this section and the next.
As in \cite[Section 10]{gt-primes} or in the earlier work of Goldston-Y{\i}ld{\i}r{\i}m, the idea is to first use the Chinese remainder theorem 
to essentially replace $\Z[i]$ with the product of more local objects such as $\Z[i]/p\Z[i]$.  We then use the non-degeneracy hypotheses on the $L_{aj}$ to compute the contribution of each local object, leaving us with an Euler product over Gaussian primes, which we will estimate by using the pole and residue of the modified Gaussian integer zeta function $\zeta_{\Z'[i]}(s)$ (which also has an Euler product representation) at $s=1$.

We turn to the details, starting with the proof of Proposition \ref{GY-Gauss}.  
We begin by eliminating the role of the box $B$.  Using \eqref{trunc-def}, we can write the left-hand side of
\eqref{gy-cor} as
\begin{align*}
\frac{\log^{2|S|} \N(R)}{16^{|S|}} \E\biggl( &\prod_{s \in S} \sum_{d_s, d'_s \in \Z'[i]: d_s,d'_s | W(\L_s \cdot \x) + a_s} \\
&\mu_{\Z[i]}(d_s) \mu_{\Z[i]}(d'_s)
\varphi\left( \frac{\log \N(d_s)}{\log \N(R)} \right) \varphi\left( \frac{\log \N(d'_s)}{\log \N(R)} \right) \bigr| \x \in B \biggr).
\end{align*}
From the support of $\varphi$, we may restrict the $d_s, d'_s$ summations to the range where $\N(d_s), \N(d'_s) \leq \N(\R)$.  We can
thus rearrange the above expression as
\begin{equation}\label{rearrange}
\begin{split}
&\frac{\log^{2|S|} \N(R)}{16^{|S|}}  \sum_{d_s,d'_s \in \Z'[i]: \N(d_s), \N(d'_s) \leq \N(R) \forall s \in S}\\
&\quad \prod_{s \in S} \mu_{\Z[i]}(d_s) \mu_{\Z[i]}(d'_s) \varphi\left( \frac{\log \N(d_s)}{\log \N(R)} \right) \varphi\left( \frac{\log \N(d'_s)}{\log \N(R)} \right) \\
&\quad\E\left( \prod_{s \in S} 1_{[d_s, d'_s] | W(\L_s \cdot \x) + a_s} \bigr| \x \in B \right),
\end{split}
\end{equation}
where $[d,d']$ is the least common multiple of $d$ and $d'$ (this is only defined up to association).
Let $D = D(([d_s,d'_s])_{s \in S}) \in \Z^+$ be the smallest positive
rational integer which is a multiple of all of the $d_s, d'_s$ for $s \in S$.
Since each of the $d_s, d'_s$ have norm at most $\N(R) = R^2$, we have
$$D \leq \prod_{s \in S} \N(d_s) \N(d'_1) = R^{4|S|}$$  
On the other hand, $B$ has sidelength at least $R^{10|S|}$.  Since the solutions $\x$ to the system
$$ [d_s, d'_s] | W(\L_s \cdot \x) + a_s \hbox{ for all } s \in S$$
are periodic of period $D$ (of course, the period could in fact be smaller) in each component of $\x$, we thus conclude that
$$ \E( \prod_{s \in S} 1_{[d_s, d'_s] | W(\L_s \cdot \x) + a_s} | \x \in B ) =
\omega(([d_s,d'_s])_{s \in S}) + O_{|T|}(R^{-6|S|})$$
where $\omega$ is the expression
\begin{equation}\label{omega-def}
 \omega((d_s)_{s \in S}) := \E( \prod_{s \in S} 1_{d_s | W(\L_s \cdot \x) + a_s} | \x \in (\Z/D\Z)^T ).
 \end{equation}
Note that $\omega$ is unchanged if one of its arguments is replaced by an associate, so we may legitimately use expressions such as
$[d_s, d'_s]$ in the arguments of $\omega$. The contribution of the error term $O_{|T|}(R^{-6|S|})$ to \eqref{rearrange} is at most
$$ O_{|S|,|T|}( (\log^{|S|} \N(R)) \N(R)^{2|S|} R^{-6|S|} ) $$
which is certainly acceptable.  Thus it suffices to show that
\begin{equation}\label{nanda}
\begin{split} 
&\sum_{d_s, d'_s \in \Z'[i]: \N(d_s), \N(d'_s) \leq \N(R) \forall s \in S}\\
&\prod_{s \in S} \mu_{\Z[i]}(d_s) \mu_{\Z[i]}(d'_s) \varphi\left( \frac{\log \N(d_s)}{\log \N(R)} \right) \varphi\left( \frac{\log \N(d'_s)}{\log \N(R)} \right)
\omega(([d_s,d'_s])_{s \in S})\\
&=
(1 + o_{R \to \infty; |S|,|T|,\varphi,W, (\L_s)_{s \in S}}(1) + o_{W \to \infty; |S|,|T|,\varphi, (\L_s)_{s \in S}}(1)) 
\left(\frac{c'_\varphi}{\log \N(R)} \frac{\N(W)}{\phi_{\Z[i]}(W)}\right)^{|S|}
\end{split}
\end{equation}
for some $c'_\varphi > 0$ depending only on $\varphi$.
Now observe that since the $a_s$ are coprime to $W$, so are $W(\L_s \cdot \x) + a_s$.  Thus the above summand vanishes if
any of the $d_s, d'_s$ share a common factor with $W$.  Thus we reduce to showing that
\begin{equation}\label{nanda-2}
\begin{split} 
&\sum_{d_s,d'_s\in \Z[i]'_{sq,W} \hbox{ for all } s \in S}
\prod_{s \in S} \mu_{\Z[i]}(d_s) \mu_{\Z[i]}(d'_s) \varphi\left( \frac{\log \N(d_s)}{\log \N(R)} \right) \varphi\left( \frac{\log \N(d'_s)}{\log \N(R)} \right)
\omega(([d_s,d'_s])_{s \in S})\\
&=
(1 + o_{R \to \infty; |S|,|T|,\varphi,W, (\L_s)_{s \in S}}(1) + o_{W \to \infty; |S|,|T|,\varphi, (\L_s)_{s \in S}}(1)) \left(\frac{c'_\varphi}{\log \N(R)} \frac{\N(W)}{\phi_{\Z[i]}(W)}\right)^{|S|}
\end{split}
\end{equation}
where $\Z[i]'_{sq,W} := \{ n \in \Z'[i]_W: (n,W) = 1 \}$.  Also, we 
have taken advantage of the supports of the $\varphi$ to drop the restrictions $\N(d_s), \N(d'_s) \leq \N(R)$.

To proceed further, we need to understand the quantity $\omega( (d_s)_{s \in S} )$.  This quantity clearly ranges between 0 and 1, but
much better estimates are possible.  Firstly, we observe that $\omega$ is partially multiplicative:

\begin{lemma}\label{mult}  If $d_s \in \Z[i]'_{sq,W}$ for all $s \in S$, we have
$$ \omega((d_s)_{s \in S}) = \prod_{n \in \N(P'[i]): n \geq w} \omega(((d_s,n))_{s \in S}) $$
where $(d,n)$ is the greatest common divisor of $d$ and $n$ in the Gaussian primes (defined up to association).
\end{lemma}

\begin{proof}  Observe from unique factorization (and the hypothesis $d_s \in \Z[i]'_{sq,W}$) that solving the linear system
$$d_s | W( \L_s \cdot \x ) + a_s \hbox{ for all } s \in S$$
is the same as solving the linear systems
$$(d_a,n) | W ( \L_s \cdot \x ) + a_s \hbox{ for all } s \in S$$
simultaneously for each $n \in \N(P'[i])$ with $n \geq w$.  Note that
each individual linear system is then periodic with period $n \in \N(P'[i])$.  Since all the elements of $\N(P'[i])$ are rational primes,
the claim then follows from the Chinese remainder theorem.
\end{proof}

The above lemma splits $\omega$ into local expressions at a single value of $n \in \N(P'[i])$.  We now estimate each of these local terms; it is here
 that we must use the various non-degeneracy hypotheses we have placed on the $L_{aj}$.

\begin{lemma}[No significant local correlations]\label{local}  Let $n \in \N(P'[i])$ be such that $n \geq w$, and for each $s$ let
 $d_s\in Z[i]'_{sq,W}$ be such that $d_s | n$ (thus for fixed $n$ there are only four possible values of $d_s$, up to association).  Suppose that $w$ is sufficiently large depending on the linear forms $(\L_s)_{s \in S}$.  Then
$\omega((d_s)_{s \in S}) = 1$ if $\prod_{s \in S} d_s$ is a Gaussian unit, 
$\omega((d_s)_{s \in S}) = 1/n$ if $\prod_{s \in S} d_s$ is a Gaussian prime, 
and $\omega((d_s)_{s \in S}) = O(1/n^2)$ otherwise.
\end{lemma}

\begin{proof}  The claim is trivial when $\prod_{s \in S} d_S$ is a Gaussian unit.  Now suppose that $\prod_{s \in S} d_S$ is a 
Gaussian prime $p$, which is necessarily unexceptional since $d_1,\ldots,d_m \in Z[i]'_{sq,W}$.  We thus have
$\N(p) = n$, and one of the $d_s$ is associate to $p$, with the remaining $d_s$ being Gaussian units.  By \eqref{omega-def}, it suffices to show that
$$ \E( 1_{p | W( \L_s \cdot \x ) + a_s } | \x \in (\Z/n\Z)^T ) = 1/n$$
for each $s \in S$.  Since $\N(p)=n \geq w$, we see that $W$ is invertible in $\Z[i]/p\Z[i]$, and so the map
$x \mapsto Wx + a_s$ is a bijection on $\Z[i]/p\Z[i]$.  It will thus suffice to show that the homomorphism
from $(\Z/n\Z)^T$ to $\Z[i]/p\Z[i]$ induced by the map $\x \mapsto \L_s \cdot \x$ is
surjective.  But since $p$ is unexceptional, $\Z[i]/p\Z[i]$ is a cyclic group of prime order.  Since the linear part $\L_s$ of $\psi_s$ is not identically zero, the claim follows if $w$ is assumed sufficiently large.

Now suppose $\prod_{s \in S} d_s$ is not a unit or a Gaussian prime, then there exist Gaussian primes $p, p'$ with norm $\N(p)=\N(p')=n$,
and indices $s,s'$, with either $s \neq s'$ or $p$ not associate to $p'$, such that $d_s$ is a multiple of $p$ and 
$d_{s'}$ is a multiple of $p'$.  It thus suffices to show that
$$ \E( 1_{p | W( \L_s \cdot \x ) + a_s } 1_{p' | W( \L_{s'} \cdot \x ) + a_{s'} } 
| {\mathbf{x}} \in (\Z/n\Z)^T ) \leq 1/n^2.$$
Observe that $n^2$ is the cardinality of $\Z[i]/p\Z[i] \times \Z[i]/p'\Z[i]$.  Again, since $\N(p) = \N(p') > w$,
the map $(x,y) \mapsto (Wx + a_s, Wy + a_{s'})$ is a bijection on $\Z[i]/p\Z[i] \times \Z[i]/p'\Z[i]$.  It thus
suffices to show that the homomorphism $\Phi$ from $(\Z/n\Z)^T$ to $\Z[i]/p\Z[i] \times \Z[i]/p'\Z[i]$
induced by $\x \mapsto (\L_s \cdot \x, \L_{s'} \cdot \x)$ is surjective.

Suppose first that $s \neq s'$.  Observe that as $p$ and $p'$ are Gaussian primes with the same norm $n$, they are either 
associate to each other, or else $p$ is associate to the complex conjugate of $p'$.  In the latter case we may replace $\L_s, a_s$ with their complex conjugates $\overline{\L_s}$, $\overline{a_s}$; note that this does not affect the hypotheses we have placed on the $L_{aj}$ or $c_a$.  Thus up to association we may assume that $p=p'$.

Suppose for contradiction that $\Phi$ is not surjective, then its image is a proper subgroup of $\Z[i]/p\Z[i] \times \Z[i]/p\Z[i]$, i.e.
a line or the origin (note that $\Z[i]/p\Z[i]$ is a finite field of rational prime order, since $p \in P[i]'$ is unexceptional).  
Since the $\L_s$ are non-zero, the latter option is ruled out (if $W$ is large enough depending on the $\L_s$).
Thus the image is a line.  This forces $\L_s$ and $\L_{s'}$ to be concurrent in the finite field geometry $(\Z/n\Z)^T$.
But this implies that $L_{st} L_{s't'} - L_{st'} L_{s't'}$ is divisible by $n$ for all $t,t' \in T$.  If $w$ and hence $n$ is sufficiently large depending on the $L_{st}$, we conclude that $L_{st} L_{s't} - L_{st'} L_{s't} = 0$ for all $t,t' \in T$, but this forces
$\L_s$ and $\L_{s'}$ to be $\Q[i]$-multiples of each other, contradicting the incommensurability hypothesis.

It remains to consider the case when $s=s'$, which forces $p'$ to be associate to a conjugate of $p$.  Performing the conjugation, it suffices to show that the homomorphism $(\Z/n\Z)^t$ to $\Z[i]/p\Z[i])^2$ induced by $\x \mapsto (\L_s \cdot \x, \overline{\L_s} \cdot x)$
is not surjective.  But this follows by arguing as before (using the hypothesis that $\L_s$ is not self-incommensurate).
\end{proof}

As a particular corollary we obtain the following crude estimate:

\begin{lemma}\label{crude}  If $d_s \in P[i]'_W$ for all $s \in S$, we have
$$ \omega((d_s)_{s \in S}) \leq \frac{1}{[(\N(d_s))_{s \in S}]}$$
where $[(\N(d_s))_{s \in S}]$ is the least common multiple of the $\N(d_s)$.
\end{lemma}

\begin{proof}
Using Lemma \ref{mult} it suffices to verify this when $d_s$ all divide $n$ for some $n \in \N(P[i])$ with $n \geq w$.
But then this follows from Lemma \ref{local}, just by using the crude bound $\omega((d_s)_{s \in S}) \leq 1/n$ whenever $\prod_{s \in S} d_s$ is not
a Gaussian unit.
\end{proof}

With these estimates in hand, we can now return to proving \eqref{nanda-2}.  We would like to take advantage of the multiplicativity of
$\omega$ to obtain a Euler factorization of the left-hand side, but we must first deal with the non-multiplicative factors $\varphi$.
This we shall do by Fourier expansion\footnote{One could also express $\varphi$ as a contour integral, which amounts to much the same thing.}.  Since $\varphi(x)$ is smooth and compactly supported, so is $e^x \varphi(x)$, and so we have an expansion
\begin{equation}\label{exvarphi}
e^x \varphi(x) = \int_{-\infty}^\infty \psi(t) e^{-ixt}\ dt
\end{equation}
for some function $\psi$ depending on $\varphi$ which is rapidly decreasing in the sense that $\psi(t) = O_A( (1 + |t|)^{-A} )$ for all $A > 0$.
In particular $\psi$ is absolutely integrable and there there will be no difficulty justifying interchange of sums and integrals in what follows.
We can now expand
$$ \varphi\left( \frac{\log \N(d)}{\log \N(R)} \right) = \int_{-\infty}^\infty \N(d)^{-(1+it)/\log \N(R)} \psi(t)\ dt.$$
We could substitute this into \eqref{nanda-2}, which is essentially what is done in \cite{gt-primes} (and in the earlier work of Goldston and Y{\i}ld{\i}r{\i}m
in \cite{goldston-yildirim}, \cite{goldston-yildirim-old1}, \cite{goldston-yildirim-old2}).  However, one would then eventually 
need to estimate expressions for large $t$ which would require knowledge of a zero-free region of the zeta function for $\Z[i]$ around the axis $s = 1 + it$.  While this is 
certainly possible, one can avoid any dependence on a zero-free region (other than that near $s=1$) by truncating $t$ at this stage of the argument,
thus making the argument slightly more elementary.
More precisely, let $I$ be the interval $\{ t \in \R: |t| \leq \log^{1/2} \N(R)\}$, and exploit the rapid decrease of $\psi$ to now write
$$  \int_I \N(d)^{-(1+it)/\log \N(R)} \psi(t)\ dt
= \varphi\left( \frac{\log \N(d)}{\log \N(R)} \right) + O_{A,\varphi}( d^{-1/\log \N(R)} \log^{-A} \N(R)).$$
for any $A > 0$.  Multiplying this out (and taking advantage of the fact that the $\varphi$ terms are supported
on the region where $\N(d) \leq \N(R)$), we obtain
\begin{align*}
&\int_I \dots \int_I\\
&\prod_{s \in S} \N(d_s)^{-(1+it_s)/\log \N(R)} \N(d'_s)^{-(1+it'_s)/\log \N(R)}\ \psi(t_s) \psi(t'_s) dt_s dt'_s \\
&\quad = \prod_{s \in S} \varphi\left( \frac{\log \N(d_s)}{\log \N(R)} \right) \varphi\left( \frac{\log \N(d'_s)}{\log \N(R)} \right)\\
&\quad\quad
+ O_{A,\varphi,|S|}( (\prod_{s \in S} d_s d'_s)^{-1/\log \N(R)} \log^{-A} \N(R)).
\end{align*}
This allows us to write the left-hand side of \eqref{nanda-2} as
\begin{equation}\label{nanda-3}
\begin{split}
\int_I &\dots \int_I
\sum_{d_s, d'_s \in \Z[i]'_{sq,W} \forall s \in S} \omega(([d_s,d'_s])_{s \in S})\\
&\prod_{s \in S} \mu_{\Z[i]}(d_s) \mu_{\Z[i]}(d'_s) 
\N(d_s)^{-(1+it_s)/\log \N(R)} \N(d'_s)^{-(1+it'_s)/\log \N(R)}\ \psi(t_s) \psi(t'_s) dt_s dt'_s 
\end{split}
\end{equation}
plus an error term
\begin{equation}\label{nanda-error}
\sum_{d_s, d'_s \in \Z[i]'_{sq,W} \forall s \in S}
O_{A,\varphi,|S|}( \frac{\omega(([d_s,d'_s])_{s \in S})}{|\prod_{s \in S} d_s d'_s|^{1/\log \N(R)} \log^{A} \N(R)} ).
\end{equation}

Let us first dispose of the error term.  By Lemma \ref{crude} this expression is bounded by
$$ O_{A,\varphi,|S|}( \log^{-A} \N(R) ) \sum_{d_s, d'_s \in \Z[i]'_{sq,W} \forall s \in S} 
\frac{|\prod_{s \in S} d_s d'_s|^{-1/\log \N(R)}}{[(\N([d_s,d'_s]))_{s\in S}]}  $$
which has an Euler factorization
$$ O_{A,\varphi,|S|}( \log^{-A} \N(R) ) \prod_{n \in \N(P[i]'): n \geq w}
\sum_{d_s,d'_s \in \Z[i]_+^{(n)} \forall s \in S} 
\frac{|\prod_{s \in S} d_s d'_s|^{-1/\log \N(R)}}{[(\N([d_s,d'_s]))_{s \in S}]}  $$
where $\Z[i]_+^{(n)} := \{ a+bi \in \Z[i]_W: a,b > 0; a+bi | n \}$; note this set consists of two elements for every
$n \in \N(P[i]')$. Direct calculation shows that
\begin{align*}
 \sum_{d_s,d'_s \in \Z[i]_+^{(n)} \forall s \in S} 
\frac{|\prod_{s \in S} d_s d'_s|^{-1/\log \N(R)}}{[(\N([d_s,d'_s]))_{s \in S}])]} 
&= 1 + O_{|S|}( n^{-1-1/2\log \N(R)} ) \\
&\leq (1 + n^{-1-1/2\log \N(R)})^{O_{|S|}(1)}.
\end{align*}
Thus we can bound \eqref{nanda-error} by
$$ O_{A,\varphi,|S|}( \log^{-A} \N(R) ) \prod_{n \in P} (1 + n^{-1-1/2\log \N(R)})^{O_{|S|}(1)}$$
where $P = \{2,3,5,\ldots\}$ is the set of rational primes.
Expanding out the Euler product, this can be bounded by
$$ O_{A,\varphi,|S|}\left( \log^{-A} \N(R) \zeta( 1 + \frac{1}{2 \log \N(R)})^{O_{|S|}(1)} \right)$$
where $\zeta(\sigma+it) = \sum_{n=1}^\infty \frac{1}{n^{\sigma+it}} = \prod_{q \in P} (1 - q^{-\sigma-it})^{-1}$ is the usual Riemann zeta function.
Using the crude bound $\zeta(\sigma+it) = O( 1 + 1/|\sigma-1|)$ for $\sigma > 1$ coming from the integral test, we obtain the upper bound
$$ O_{A,\varphi,|S|} ( \log^{-A + O_{|S|}(1)} \N(R)).$$
The contribution of this to \eqref{nanda-2} will be acceptable if $A$ is chosen sufficiently large depending on $|S|$.

It remains to show that the main term \eqref{nanda-3} is equal to
$$
(1 + o_{R \to \infty; |S|,|T|,\varphi,W,(\L_s)_{s \in S}}(1) + o_{W \to \infty; |S|,|T|,\varphi, (\L_s)_{s \in S}}(1)) \left(\frac{c'_\varphi}{\log \N(R)} \frac{\N(W)}{\phi_{\Z[i]}(W)}\right)^{|S|}.$$
Using Lemma \ref{mult}, we can factorize the integrand, writing \eqref{nanda-3} as
\begin{equation}\label{sextet}
16^{|S|}
\int_I \dots \int_I
K((t_s, t'_s)_{s \in S})
\ \prod_{s \in S} \psi(t_s) \psi(t'_s) dt_s dt'_s
\end{equation}
where
\begin{align*}
K(&(t_s,t'_s)_{s \in S}) := \prod_{n \in \N(P[i]'): n \geq w}
\sum_{d_s,d'_s \in \Z[i]_+^{(n)} \forall s \in S} \\
&\omega(([d_s,d'_s])_{s \in S})
\prod_{s \in S} \mu_{\Z[i]}(d_s) \mu_{\Z[i]}(d'_s) \\
&\quad\quad
\N(d_s)^{-(1+it_s)/\log \N(R)} \N(d'_s)^{-(1+it'_s)/\log \N(R)};
\end{align*}
the factor of $16^{|S|}$ comes from the freedom to multiply each of $d_s, d'_s$ by one of the four Gaussian units.

Now we control the local factor.

\begin{lemma}\label{local-control}  Let $n \in \N(P[i]')$ be such that $n \geq w$.  Then the expression
\begin{equation}\label{doofus}
\begin{split}
&\sum_{d_s, d'_s \in \Z[i]_+^{(n)} \forall s \in S} 
\omega(([d_s,d'_s])_{s \in S})\\
&\quad\prod_{s \in S} \mu_{\Z[i]}(d_s) \mu_{\Z[i]}(d'_s) 
\N(d_s)^{-(1+it_s)/\log \N(R)} \N(d'_s)^{-(1+it'_s)/\log \N(R)}
\end{split}
\end{equation}
is equal to
$$ (1 + O_{|S|}( \frac{1}{n^2} )) \prod_{s \in S} \prod_{p \in P[i]'_+: \N(p) = n} \frac{(1 - \N(p)^{-1-(1+it_s)/\log \N(R)}) (1 - \N(p)^{-1-(1+it'_s)/\log \N(R)})}{1 - \N(p)^{-1-(2+it_s+it'_s)/\log \N(R)}}.$$
\end{lemma}

\begin{proof}  By Lemma \ref{local}, all the terms in which $\prod_{s \in S} [d_s,d'_s]$ contain more than one Gaussian prime will give
a net contribution of $O_{|S|}(1/n^2)$.
We are left with those terms in which all but at most one of the expressions $[d_s, d'_s]$ are equal to 1,
with the remaining expression $[d_s, d'_s]$ equal to either $1$ or a Gaussian prime in $P[i]'_+$ with norm $n$.  We thus can write \eqref{doofus} as
\begin{align*}
&1 + \sum_{s \in S} -\N(p)^{-1-(1+it_s)/\log \N(R)} -\N(p)^{-1-(1+it'_s)/\log \N(R)} \\
&\quad\quad+ \N(p)^{-1-(2+it_s+it'_s)/\log \N(R)}\\
&+ O_{|S|}(\frac{1}{n^2})
\end{align*}
and the claim follows.
\end{proof}

From the convergence of the infinite product $\prod_{n \geq 1} 1 + O_{|S|}(1/n^2)$, we see that
$$ \prod_{n \geq w} 1 + O_{|S|}(1/n^2)
 = 1 + o_{W \to \infty;|S|}(1).$$
We thus have
\begin{align*}
&K((t_s,t'_s)_{s \in S}) = (1 + o_{W \to \infty;|S|}(1)) \\
&\quad \prod_{s \in S} \frac{\zeta_{\Z[i]'_{sq},W}(1 + \frac{2+it_s+it'_s}{\log \N(R)})}
{\zeta_{\Z[i]'_{sq},W}(1 + \frac{1+it_s}{\log \N(R)}) \zeta_{\Z[i]'_{sq},W}(1 + \frac{1+it'_s}{\log \N(R)})}
\end{align*}
where $\zeta_{\Z[i]'_{sq},W}$ is the truncated Gaussian integer zeta function
\begin{equation}\label{trunczeta-def}
 \zeta_{\Z[i]'_{sq},W}(\sigma+it) := \prod_{p \in P[i]'_+: \N(p) \geq w} \frac{1}{1 - \N(p)^{-\sigma+it}}.
 \end{equation}

Next, we obtain a crude estimate on this zeta function.

\begin{lemma}\label{trunczeta-lemma}  If $\sigma > 1$ and $|\sigma+it-1| \leq c$ for some absolute constant $c>0$, we have
$$ \zeta_{\Z[i]'_{sq},W}(\sigma+it) = (c_0 + o_{\sigma+it \to 1; W}(1)) \frac{\phi_{\Z[i]}(W)}{\N(W)} \frac{1}{\sigma+it-1}$$
for some absolute constant $c_0 > 0$ (which does not depend on any parameter).
\end{lemma}
  
\begin{proof}  Observe that
\begin{align*}
\zeta_{\Z[i]'_{sq},W}(\sigma+it) &= \frac{1}{4} \zeta_{\Z[i]}(\sigma+it) \frac{\prod_{p \in P[i]'_+: \N(p) < w} (1 - \N(p)^{-\sigma-it})}{\prod_{p \in P[i]_+ \backslash P[i]'_+} (1 - \N(p)^{-\sigma-it})}\\
&= (1 + o_{\sigma+it \to 1;W}(1)) \frac{1}{4} \zeta_{\Z[i]'_{sq}}(\sigma+it) 
\frac{\prod_{p \in P[i]'_+: \N(p) < w} (1 - \N(p)^{-1})}{ \prod_{p \in P[i]_+ \backslash P[i]'_+} (1 - \N(p)^{-1})}
\end{align*}
where $P[i]_+$ denotes those Gaussian primes in the first quadrant $\{a+bi: a > 0; b \geq 0\}$ and
$$ \zeta_{\Z[i]}(\sigma+it) := 4 \prod_{p \in P[i]_+} \frac{1}{1 - \N(p)^{-\sigma-it}}.$$
On the other hand, from the Chinese remainder theorem and the definition of $W$ we see that
$$ \prod_{p \in P[i]_+: \N(p) < w} (1 - \N(p)^{-1}) = \frac{\phi_{\Z[i]}(W)}{\N(W)}.$$
Also, since $P[i]_+ \backslash P[i]'_+$ consists of 2 and the rational primes equal to 3 modulo 4, we see that
$$ \prod_{p \in P[i]_+ \backslash P[i]'_+} (1 - \N(p)^{-1}) = c_1$$
for some absolute constant $c_1 > 0$.  To conclude the claim (for $s$ sufficiently close to 1), it will suffice to show that
$$ \zeta_{\Z[i]}(\sigma+it) = (1 + o_{\sigma+it \to 1}(1)) \frac{\pi}{\sigma+it-1}.$$
But by the unique factorization of the Gaussian integers (and the fact that there are exactly 4 Gaussian units) we have
$$ \zeta_{\Z[i]}(\sigma+it) = \sum_{n \in \Z[i] \backslash 0} \frac{1}{\N(n)^{\sigma+it}}.$$
By the integral test we can estimate
$$ \zeta_{\Z[i]}(\sigma+it) = \int_{a^2 + b^2 \geq 1} \frac{da db}{(a^2+b^2)^{\sigma+it}} + O(1)$$
which after polar co-ordinates becomes
$$ \zeta_{\Z[i]}(\sigma+it) = \frac{\pi}{\sigma+it-1} + O(1)$$
and the claim follows.
\end{proof}

Applying this lemma and recalling that $t_s,t'_s \in I$ and hence $t_s, t'_s = O(\log^{1/2} \N(R))$, we conclude that
$$ K((t_s, t'_s)_{s \in S}) = \frac{1}{c_0^{|S|}}(1 + o_{W \to \infty;|S|}(1) + o_{R \to \infty; W,|S|}(1))
\prod_{s \in S} \frac{(1 + it_s) (1 + it'_s)}{2 + it_s + it'_s}$$
and so we can write \eqref{sextet} as
\begin{align*}
&(\frac{16}{c_0} \frac{\N(W)}{\phi_{\Z[i]}(W) \log \N(R)})^{|S|} 
\int_I \dots \int_I\\
&\quad
(1 + o_{W \to \infty;|S|}(1) + o_{R \to \infty; W,|S|}(1))
\prod_{s \in S} \frac{(1 + it_s) (1 + it'_s)}{2 + it_s + it'_s}
\ \psi(t_s) \psi(t'_s) dt_s dt'_s.
\end{align*}
The contributions of the error terms $o_{W \to \infty;|S|}(1) + o_{R \to \infty; W,|S|}(1))$ will be acceptable, thanks to the rapid
decay of the $\psi$ factors (and the at most polynomial growth of the $\frac{(1 + it_s) (1 + it'_s)}{2 + it_s + it'_s}$ factors), so it suffices to estimate the main term, which factorizes as
$$
\left[\frac{16}{c_0} \frac{\N(W)}{\phi_{\Z[i]}(W) \log \N(R)} 
\int_I \int_I
\frac{(1 + it) (1 + it')}{2 + it + it'}
\ \psi(t) \psi(t') dt dt' \right]^{|S|}.$$
Using the rapid decay of the $\psi$, we can write this as
$$ \left[ (c'_\varphi + o_{R \to \infty}(1)) \frac{\N(W)}{\phi_{\Z[i]}(W) \log \N(R)}\right]^{|S|}$$
where
$$ c'_\varphi := \frac{16}{c_0} \int_{-\infty}^\infty \int_{-\infty}^\infty
\frac{(1 + it) (1 + it')}{2 + it + it'}
\ \psi(t) \psi(t') dt dt'.$$
It thus suffices to show that $c'_\varphi$ is real and positive.  We remark that this can be shown indirectly, by observing that
the left-hand side of \eqref{gy-cor} is necessarily non-negative, and when $|S|=|T|=1$ one can show using \eqref{primes-square} and a pigeonholing argument
that this left-hand side is at least $C_{\varphi,W}^{-1} \log \N(R)$ for some $C_{\varphi, W} > 0$, and all $R$ sufficiently large depending on $\varphi, W$; by choosing $W$ appropriately we obtain the positivity of $c'_\varphi$.  However, we can also argue directly via the following
Fourier-analytic argument\footnote{One can of course also use contour integration as a substitute for Fourier analysis here; the two approaches are essentially equivalent.}.
Making the change of variables $s :=t+t'$, we have
$$ \int_{-\infty}^\infty \int_{-\infty}^\infty
\frac{(1 + it) (1 + it')}{2 + it + it'}
\ \psi(t) \psi(t') dt dt' = \int_{-\infty}^\infty \frac{1}{2+is} F(s)\ ds$$
where $F$ is the convolution
$$ F(s) := \int_{-\infty}^\infty (1+it) \psi(t) (1+i(s-t)) \psi(s-t)\ dt.$$
Observe that for any real number $x$, the Fourier transform of $F$ can be computed as
\begin{align*}
\int_{-\infty}^\infty F(s) e^{-ixs}\ ds
&= \int_{-\infty}^\infty (1+it) \psi(t) e^{-ixt} (1+i(s-t)) \psi(s-t)e^{-ix(s-t)}\ dt ds \\
&= (\int_{-\infty}^\infty (1+it) \psi(t) e^{-ixt}\ dt)^2 \\
&= [(1 - \frac{d}{dx}) \int_{-\infty}^\infty \psi(t) e^{-ixt}\ dt]^2 \\
&= [e^x \varphi'(x)]^2
\end{align*}
where we have used the rapid decrease of the $\psi$ to justify all the swapping of integrals, and \eqref{exvarphi} in the last line.
Now we write $\frac{1}{2+is} = \int_0^\infty e^{-2x} e^{-ixs}\ dx$ and interchange integrals again (using the rapid decay of $F$) to conclude
$$ \int_{-\infty}^\infty \frac{1}{2+is} F(s)\ ds = \int_0^\infty e^{-2x} \int_{-\infty}^\infty F(s) e^{-ixs}\ ds dx =
\int_0^\infty [\varphi'(x)]^2\ dx$$
and hence
$$ c'_\varphi = \frac{64}{\pi} \int_0^\infty [\varphi'(x)]^2\ dx  > 0$$
as desired.  This concludes the proof of Proposition \ref{GY-Gauss}.
\endprf

\section{Proof of Proposition \ref{GY-Gauss-2}}

Now we turn to Proposition \ref{GY-Gauss-2}.  This will be similar to the proof of Proposition \ref{GY-Gauss} in the
preceding section, but with a number of differences.  It is a little simpler because there is only one parameter $n$ to sum over rather
than $|T|$ parameters, and also we only seek an upper bound rather than an asymptotic.  As such we shall move more rapidly with this proof as compared
with the similar but more complicated proof from the previous section.

We begin by eliminating the role of the interval $I$.  Using \eqref{trunc-def}, we can rewrite the left-hand side of
\eqref{parade} as
\begin{align*}
\frac{\log^{2m} \N(R)}{16^m}
\sum_{d_1,d'_1,\ldots,d_m,d'_m \in Z[i]'_{sq}} &\prod_{j=1}^m \varphi( \frac{\log \N(d_j)}{\log \N(R)} ) \varphi( \frac{\log \N(d_j)}{\log \N(R)} )
\mu_{\Z[i]}(d_j) \mu_{\Z[i]}(d'_j) \\
&\E( \prod_{j=1}^m 1_{d_j, d'_j | W(h_j+nv)+b} | n \in I ).
\end{align*}
Due to the support of the $\varphi$, we can restrict the $d_j$ and $d'_j$ to the region $\N(d_j), \N(d'_j) \leq \N(R)$.  Now from the Chinese
remainder theorem we have
$$ \E( \prod_{j=1}^m 1_{d_j, d'_j | W(h_j+nv)+b} | n \in I ) = 
\omega([d_1,d'_1],\ldots,[d_m,d'_m]) + O_v( R^{4m} / |I| )$$
where
$$ \omega(q_1,\ldots,q_m) := \E( \prod_{j=1}^m 1_{q_j | W(h_j+nv)+b} | n \in \Z/D\Z ).$$
and $D = D(q_1,\ldots,q_m)$ is the smallest positive rational integer which is a multiple of all
the $d_1,\ldots,d_m$. In our situation we have the crude estimate $D = O(R^{4m})$.  Since $|I| \geq R^{10m}$, it is easy to see
that the contribution of the error term $O_v(R^{4m}/|I|)$ is acceptable (if $W$ is large enough depending on $v$, but is sufficiently slowly growing in $N$).  Thus it suffices to show that
\begin{equation}\label{nanda-corr}
\begin{split}
&\sum_{d_1,d'_1,\ldots,d_m,d'_m \in Z[i]'_{sq}} \prod_{j=1}^m \varphi( \frac{\log \N(d_j)}{\log \N(R)} ) \varphi( \frac{\log \N(d_j)}{\log \N(R)} )
\mu_{\Z[i]}(d_j) \mu_{\Z[i]}(d'_j) \omega([d_1,d'_1],\ldots,[d_m,d'_m]) \\
&= 
O_{m,v}\left( \left(\frac{\N(W)}{\phi_{\Z[i]}(W) \log \N(R)} \right)^m 
\prod_{p \in \P[i]'_+: p | \Delta} (1 + O_m( \N(p)^{-1/2} )) \right)
\end{split}
\end{equation}
Here we have used the support of $\varphi$ to drop the constraints $\N(d_j), \N(d'_j) \leq \N(R)$ again.
Now observe that $\omega$ vanishes if any one of the $d_j$ or $d'_j$ shares a common factor with $W$, since $b$ is coprime to $W$.  Thus without loss of generality we may restrict $d_1,\ldots,d_m,d'_1,\ldots,d'_m$ to $\Z[i]'_{sq,W}$.

Now we must obtain analogues to Lemmas \ref{mult}, \ref{local}, \ref{crude}.  By repeating the proof of Lemma \ref{mult} with
only trivial changes, we have

\begin{lemma}\label{mult-2}  If $d_1,\ldots,d_m,d'_1,\ldots,d'_m \in \Z[i]'_{sq,W}$, we have
$$ \omega(q_1,\ldots,q_m) = \prod_{n \in \N(P'[i]): n \geq w} \omega((q_1,n),\ldots,(q_m,n)).$$
\end{lemma}

Now we give the analogue of Lemma \ref{local}.

\begin{lemma}[No significant local correlations]\label{local-2}  Let $n \in \N(P'[i])$ be such that $n \geq w$, and let
$q_1,\ldots,q_m \in Z[i]'_{sq,W}$ divide $n$.  Suppose that $w$ is sufficiently large depending on $v$.
Then
$\omega(q_1,\ldots,q_m) = 1$ if $q_1 \ldots q_m$ is a Gaussian unit, 
and $\omega(q_1,\ldots,q_m) = 1/n$ if $q_1 \ldots q_m$ is a Gaussian prime.
In all other cases, we have $\omega(q_1,\ldots,q_m) = O(1/n) 1_{n|\Delta}$,
where $\Delta$ was defined in \eqref{Delta-def}.
\end{lemma}

\begin{proof}  In the first two cases (when $q_1\ldots q_m$ is a Gaussian unit or a Gaussian prime), the
claim follows just as in Lemma \ref{local}, noting that $\Z[i]/p\Z[i]$ is cyclic of prime order $n$ whenever $\N(p)=n$, and that $W$ and $v$
are invertible in $\Z[i]/p\Z[i]$.  

Now suppose that $q_1\ldots q_m$ is the product of at least two primes.  Then from the preceding
discussion we certainly have $\omega(q_1,\ldots,q_m) \leq 1/n$, by discarding all but one of the constraints
$q_j | W(h_j+nv)+b$.  This settles the claim when $n$ divides $\Delta$, so now suppose that $n$ does not divide $\Delta$.
This implies in particular that the $h_j$ are all distinct in $\Z[i]/p\Z[i]$ for any Gaussian prime $p$ dividing $n$.  Since $W$ and $v$ are also invertible in $\Z[i]/p\Z[i]$, this means any two constraints of the form $p | W(h_j + nv) + b$ and $p | W(h_{j'} + nv) + b$ cannot simultaneously
be true for any distinct $j,j'$.  In a similar spirit, since $\Delta$ is non-zero in $\Z[i]/p\Z[i]$ for any $p$ dividing $n$, we see that
$W(h_j \overline{v} - \overline{h_{j'}} v) - b \overline{v} + \overline{b} v$ is similarly non-zero for any $1 \leq j \leq j' \leq m$.  
A little algebra then shows that the constraints $p | W(h_j + nv)+b$ and $\overline{p} | W(h_{j'} + nv) + b$ cannot simultaneously be true.
Combining all these facts together, we see that the constraints $q_j | W(h_j+nv) + b$ cannot be simultaneously satisfied for $1 \leq j \leq m$,
anmd $\omega(q_1,\ldots,q_m)$ vanishes as claimed.
\end{proof}

As a particular corollary we obtain the analogue of Lemma \ref{crude}:

\begin{lemma}\label{crude-2}  If $d_1,\ldots,d_m,d'_1,\ldots,d'_m \in P[i]'_W$, then
$$ \omega(d_1,\ldots,d_m,d'_1,\ldots,d'_m) \leq \frac{1}{[\N(d_1),\ldots,\N(d_m),\N(d'_1),\ldots,\N(d'_m)]}.$$
\end{lemma}

The proof is the same as that of Lemma \ref{crude} and is omitted.

We return to the proof of \eqref{nanda-corr}.  Once again, we use the expansion \eqref{exvarphi} of $e^x \varphi(x)$, and obtain the expansion
\begin{align*}
\int_I \dots \int_I &\prod_{j=1}^m \N(d_j)^{-(1+it_j)/\log \N(R)} \N(d'_j)^{-(1+it'_j)/\log \N(R)}\ \psi(t_j) \psi(t'_j) dt_j dt'_j \\
&\quad = \prod_{j=1}^m \varphi( \frac{\log \N(d_j)}{\log \N(R)} ) \varphi( \frac{\log \N(d'_j)}{\log \N(R)} )\\
&\quad\quad
+ O_{A,\varphi,m}\left( (\prod_{j=1}^m d_j d'_j)^{-1/\log \N(R)} \log^{-A} \N(R)\right)
\end{align*}
where $I$ is the interval $I := \{ t \in \R: |t| \leq \log^{1/2} \N(R)\}$, $\psi$ is rapidly decreasing, and $A > 0$ is arbitrary.
This allows us to write the left-hand side of \eqref{nanda-corr} as a main term
\begin{equation}\label{nanda-main-corr}
\begin{split}
\int_I \ldots \int_I \sum_{d_1,d'_1,\ldots,d_m,d'_m \in Z[i]'_{sq}} \prod_{j=1}^m 
&\mu_{\Z[i]}(d_j) \mu_{\Z[i]}(d'_j) \omega([d_1,d'_1],\ldots,[d_m,d'_m]) \\
&\N(d_j)^{-(1+it_j)/\log \N(R)} \N(d'_j)^{-(1+it'_j)/\log \N(R)}\ \psi(t_j) \psi(t'_j) dt_j dt'_j 
\end{split}
\end{equation}
plus an error term
$$
\sum_{d_1,\ldots,d_m,d'_1,\ldots,d'_m \in \Z[i]'_{sq,W}}
O_{A,\varphi,m}\left( \frac{\omega([d_1,d'_1],\ldots,[d_m,d'_m])}{|d_1\ldots d_m d'_1 \ldots d'_m|^{1/\log \N(R)} \log^{A} \N(R)} \right).
$$

The error term is treated exactly as with \eqref{nanda-error}, so we turn to treating the main term
\eqref{nanda-main-corr}.  Our task is to estimate this term by
\begin{equation}\label{movit}
O_{m,v}\left( \left(\frac{\N(W)}{\phi_{\Z[i]}(W) \log \N(R)} \right)^m 
\prod_{p \in \P[i]'_+: p | \Delta} (1 + O_m( \N(p)^{-1/2} )) \right).
\end{equation}
Using Lemma \ref{mult-2}, we can rewrite \eqref{nanda-main-corr} as
\begin{equation}\label{sextet-2}
16^m \int_I \dots \int_I
K(t_1,\ldots,t_m,t'_1,\ldots,t'_m)
\ \prod_{j=1}^m \psi(t_j) \psi(t'_j) dt_j dt'_j
\end{equation}
where
\begin{align*}
K(&t_1,\ldots,t_m,t'_1,\ldots,t'_m) := \prod_{n \in \N(P[i]'): n \geq w}
\sum_{d_1,\ldots,d_m,d'_1,\ldots,d'_m \in \Z[i]_+^{(n)}} \\
&\omega([d_1,d'_1],\ldots,[d_m,d'_m])
\prod_{j=1}^m \mu_{\Z[i]}(d_j) \mu_{\Z[i]}(d'_j) \\
&\quad\quad
\N(d_j)^{-(1+it_j)/\log \N(R)} \N(d'_j)^{-(1+it'_j)/\log \N(R)}.
\end{align*}

Now we control the local factor, in complete analogy with Lemma \ref{local-control}.

\begin{lemma}  Let $n \in \N(P[i]')$ be such that $n \geq w$.  Then the expression
\begin{equation}\label{doofus-2}
\begin{split}
&\sum_{d_1,\ldots,d_m,d'_1,\ldots,d'_m \in \Z[i]_+^{(n)} \forall s \in S} 
\omega(([d_1,d'_1],\ldots,[d_m,d'_m])_{s \in S})\\
&\quad\prod_{j=1}^m \mu_{\Z[i]}(d_j) \mu_{\Z[i]}(d'_j) 
\N(d_j)^{-(1+it_j)/\log \N(R)} \N(d'_s)^{-(1+it'_j)/\log \N(R)}
\end{split}
\end{equation}
is equal to $1 + O_m(1/n)$ if $n$ divides $\Delta$, and is equal to
$$ (1 + O_{m}( \frac{1}{n^2} )) \prod_{j=1}^m \prod_{p \in P[i]'_+: \N(p) = n} \frac{(1 - \N(p)^{-1-(1+it_j)/\log \N(R)}) (1 - \N(p)^{-1-(1+it'_j)/\log \N(R)})}{1 - \N(p)^{-1-(2+it_j+it'_j)/\log \N(R)}}$$
otherwise,.
\end{lemma}

\begin{proof}  If $n$ divides $\Delta$, then the claim follows from Lemma \ref{crude-2}, so suppose that $n$ does not divide $\Delta$.
But then the claim follows by exact repetition of the proof of Lemma \ref{local-control}.
\end{proof}

From the above lemma we see that
\begin{align*}
K(t_1,\ldots,t_m,t'_1,\ldots,t'_m) = &O_m\biggl( \left[\prod_{n \in \N(P[i]'): n \geq w, n | \Delta} (1 + O_m(1/n))\right]  \\
&\prod_{j=1}^m 
\frac{\zeta_{\Z[i]'_{sq},W}(1 + (2+it_j+it'_j)/\log \N(R))}
{\zeta_{\Z[i]'_{sq},W}(1 + (1+it_j)/\log \N(R)) \zeta_{\Z[i]'_{sq},W}(1 + (1+it'_j)/\log \N(R))} \biggr)
\end{align*}
where $\zeta_{\Z[i]'_{sq},W}$ was defined in \eqref{trunczeta-def}.  Applying Lemma \ref{trunczeta-lemma}, we conclude
\begin{align*}
K(t_1,\ldots,t_m,t'_1,\ldots,t'_m) &= O_m\biggl( \left[\prod_{n \in \N(P[i]'): n \geq w, n | \Delta} (1 + O_m(1/n))\right] \\ 
&\frac{\N(W)^m}{\phi_{\Z[i]}(W)^m \log^m \N(R)} 
\prod_{j=1}^m \frac{(1+|t_j|)(1+|t'_j|)}{1+|t_j+t'_j|} \biggr).
\end{align*}
Inserting this into \eqref{sextet-2} and using the rapid decay of $\psi$, we can thus bound \eqref{sextet-2} by
$$ O_m( [\prod_{n \in \N(P[i]'): n \geq w, n | \Delta} (1 + O_m(1/n))]  
\frac{\N(W)^m}{\phi_{\Z[i]}(W)^m \log^m \N(R)} )$$
which is bounded by \eqref{movit} as desired.
This concludes the proof of Proposition \ref{GY-Gauss} and hence Theorem \ref{main}.
\endprf

\section{Discussion}\label{discussion-sec}

The proof of Theorem \ref{main} also gives a little bit more, namely that any subset of the Gaussian primes $P[i]$ of positive \emph{relative} density
will contain infinitely many constellations of a prescribed shape, but we have chosen not to give this generalization in order to 
simplify the exposition slightly.

Our method is also likely to extend to other number fields than the Gaussian integers, at least if one has unique factorization, a finite Galois group, and only finitely many units (though these constraints certainly place severe restrictions on which number fields are available!).  A good ``litmus test'' seems to be whether the number field supports a reasonable notion of an \emph{almost prime}, and whether sieve theory type techniques can easily produce a large number of constellations amongst these almost primes.  If this is the case, then there is a good chance that the methods here will then extend to primes (or irreducibles), and dense subsets thereof.  It is also likely
that a relative version of Theorem \ref{multiszemeredi} exists, in which the set $A$ is a dense subset of the set $P^d$ - the set of lattice points in $\Z^d$ with prime coefficients - rather than $\Z^d$.  However, a technical problem arises when working with $P^d$, namely that $P^d$ (or any majorant
of $P^d$) generates significant correlations between certain elements $a+rv_j$ of the constellation, even after removing obstructions coming from small divisors.  For instance, if $a + r (1,0)$ and $a + r(0,1)$ both lie in $P^2$, then $a$ itself necessarily also lies in $P^2$.  This issue means
that the entire approach to this problem, based on viewing $P^d$ as a dense subset of some suitably pseudorandom set, needs to be somehow modified,
unless one is working in a case when these correlations do not appear (for instance, if the $i^{th}$ co-ordinate of the $v_j$ are distinct in $j$
for each $i$).  However even in such a model case there appear to be some non-trivial technical difficulties, most notably
in obtaining the dual function condition (Definition \ref{dual-c}).  We will not pursue these matters here.


\begin{thebibliography}{10}

\bibitem{frankl02}
P. Frankl, V. R\"odl, \emph{Extremal problems on set systems}, Random Struct. Algorithms \textbf{20} (2002), no. 2, 131-164. 

\bibitem{fk}
H. Furstenberg, Y. Katznelson, \emph{An ergodic Szemer\'edi theorem for commuting transformations}. J. Analyse Math. \textbf{34} (1978), 275--291.

\bibitem{fk2}
H. Furstenberg, Y. Katznelson, \emph{A density version of the Hales-Jewett theorem}, J. d'Analyse Math. \textbf{57} (1991), 64--119.

\bibitem{goldston-yildirim-old1} D. Goldston and C.Y. Y{\i}ld{\i}r{\i}m, \emph{Higher correlations of divisor sums related to primes, I: Triple correlations,} Integers \textbf{3} (2003) A5, 66pp.

\bibitem{goldston-yildirim-old2} D. Goldston and C.Y. Y{\i}ld{\i}r{\i}m, \emph{Higher correlations of divisor sums related to primes, III: $k$-correlations,} preprint (available at AIM preprints)

\bibitem{goldston-yildirim} D. Goldston and C.Y. Y{\i}ld{\i}r{\i}m, \emph{Small gaps between primes, I,} preprint.

\bibitem{gowers}
T. Gowers, \emph{A new proof of Szemeredi's theorem}, GAFA \textbf{11} (2001), 465--588.

\bibitem{gowers-hyper}
T. Gowers, \emph{Hypergraph regularity and the multidimensional Szemer\'edi theorem}, preprint.

\bibitem{gt-primes}
B. Green, T. Tao, \emph{The primes contain arbitrarily long arithmetic progressions}, preprint.

\bibitem{hardy-wright}
G.H. Hardy, E.M. Wright, \emph{An introduction to the theory of numbers}, 5th Ed., Oxford, Clarendon Press, 1979.

\bibitem{krs-hyper}
Y. Kohayakawa, V. R\"odl, J. Skokan,
\emph{Hypergraphs, quasi-randomness, and conditions for regularity},
J. Combin. Theory Ser. A \textbf{97} (2002), no. 2, 307--352.

\bibitem{nrs}
B. Nagle, V. R\"odl, J. Skokan, \emph{The counting lemma for regular $k$-uniform hypergraphs}, to appear, Random Structures and Algorithms.

\bibitem{zoo}
J. Renze, S. Wagon, B. Wick, \emph{The Gaussian Zoo}, Experimental Math. \textbf{10} (2001), 161--173.

\bibitem{rodl}
V. R\"odl, J. Skokan, \emph{Regularity lemma for $k$-uniform hypergraphs}, to appear, Random Structures and Algorithms.

\bibitem{rodl2}
V. R\"odl, J. Skokan, \emph{Applications of the regularity lemma for uniform hypergraphs}, preprint.

\bibitem{rsz}
I. Ruzsa, E. Szemer\'edi, \emph{Triple systems with no six points carrying three triangles}, Colloq. Math. Soc. J. Bolyai \textbf{18} (1978), 939--945.

\bibitem{soly-roth}
J. Solymosi, \emph{Note on a generalization of Roth's theorem}, Discrete and computational geometry, 825--827, Algorithms Combin. \textbf{25}, Springer Verlag, 2003.

\bibitem{soly-2}
J. Solymosi, \emph{A note on a question of Erd\"os and Graham}, Combinatorics, Probability and Computing \textbf{13} (2004), 263--267.

\bibitem{szemeredi}
E. Szemer\'edi, \emph{On sets of integers containing no four elements in arithmetic progression},
Acta Math. Acad. Sci. Hungar. \textbf{20} (1969), 89--104.

\bibitem{szemeredi-4}
E. Szemer\'edi, \emph{On sets of integers containing no $k$ elements in arithmetic progression},
Acta Arith. \textbf{27} (1975), 299--345.

\bibitem{tao:ergodic}
T. Tao, \emph{A quantitative ergodic theory proof of Szemer\'edi's theorem}, preprint.

\bibitem{tao:hyper}
T. Tao, \emph{A variant of the hypergraph removal lemma}, preprint.

\bibitem{tao:survey}
T. Tao, \emph{Arithmetic progressions in the primes}, El Escorial conference proceedings.

\bibitem{tao-gy}
T. Tao, \emph{A remark on Goldston-Y{\i}ld{\i}r{\i}m correlation estimates}, unpublished.

\end{thebibliography}
\end{document}